\documentclass{amsart}
\usepackage{amssymb, amsmath, amsthm}
\usepackage{epsfig}
\usepackage{epstopdf}
\usepackage{graphicx}
\usepackage{color}

\theoremstyle{plain}
\newtheorem{theorem}{Theorem}[section]

\newtheorem{proposition}[theorem]{Proposition}
\newtheorem{corollary}[theorem]{Corollary}
\newtheorem{lemma}[theorem]{Lemma}

\theoremstyle{definition}

\newtheorem*{question}{Question}
\newtheorem*{definition}{Definition}
\newtheorem*{remark}{Remark}
\newtheorem*{example}{Example}

\newtheorem*{thmA}{Theorem A}
\newtheorem*{thmB}{Theorem B}
\newtheorem*{thmC}{Theorem C}
\newtheorem*{thm6.4}{Theorem 6.4}
\newtheorem*{prop6.5}{Proposition 6.5}
\newtheorem*{prop6.6}{Proposition 6.6}
\newtheorem*{prop6.7}{Proposition 6.7}
\newtheorem*{prop6.8}{Proposition 6.8}

\newcommand{\R}{{\mathbb R}}

\newcommand{\Z}{\mathbb Z}
\newcommand{\N}{\mathbb N}

\newcommand{\nil}{\varnothing}

\newcommand{\defn}[1]{\textbf{#1}}
\newcommand{\boundary}{\partial}

\newcommand{\mc}[1]{\mathcal{#1}}
\newcommand{\ob}[1]{\overline{#1}}
\newcommand{\fr}{\operatorname{fr}} 
\newcommand{\cl}{\operatorname{cl}} 
\newcommand{\interior}{\operatorname{int}} 
\newcommand{\slide}[1]{{\curvearrowright \atop #1}}

      \makeatletter
      \def\@setcopyright{}
      \def\serieslogo@{}
      \makeatother
      
\begin{document}


   \title[On Non-compact Heegaard Splittings]{On Non-compact Heegaard Splittings}
   \author{Scott Taylor}
   \address{Mathematics Department, University of California, Santa Barbara}
   \email{staylor@math.ucsb.edu}



 \maketitle
 
 \begin{abstract}
A Heegaard splitting of an open 3-manifold is the partition of the manifold into two non-compact handlebodies which intersect on their common boundary.  This paper proves several non-compact analogues of theorems about compact Heegaard splittings.  The main theorem is: if N is a compact, connected, orientable 3-manifold with non-empty boundary, with no $S^2$ components, and if $M$ is obtained from $N$ by removing the boundary then any two Heegaard splittings of $M$ are properly ambient isotopic.  This is a non-compact analogue of the classifications of splittings of $\text{(closed surface)} \times I$ and $\text{(closed surface)} \times S^1$ by Scharlemann-Thompson and Schultens.  Work of Frohman-Meeks and a non-compact analogue of the Casson-Gordon theorem on weakly reducible Heegaard splittings are key tools.
 \end{abstract}
 
 
\section{Introduction}

Non-compact 3-manifolds vary widely in the degree to which they are similar to compact 3-manifolds.  The most tractable are the deleted boundary 3-manifolds which are obtained by removing boundary components from compact 3-manifolds.  At the other end of the spectrum are those which are not the connect sums of prime 3-manifolds.  Heegaard splittings play an important role in compact 3-manifolds, so it is interesting to ask about the ways in which this structure can be extended to non-compact manifolds.  This paper studies Heegaard splittings of an important class of 3-manifolds: eventually end-irreducible 3-manifolds.  One of the main results is that eventually end-irreducible 3-manifolds have exhausting sequences which interact nicely with a given Heegaard splitting.  This result is applied to classify Heegaard splittings of ``most" deleted boundary 3-manifolds.  The only deleted boundary 3-manifolds which have splittings that cannot be classified are those which are obtained from a compact 3-manifold by removing finitely many closed 3-balls. \newline

The study of non-compact Heegaard splittings was initiated by Frohman and Meeks in \cite{FrMe97}.  They showed that any two infinite genus Heegaard surfaces in $\R^3$ are properly ambient isotopic and used this result to give a topological classification of complete 1-ended minimal surfaces in $\R^3$.  In this paper, Frohman and Meek's result is extended to non-compact 3-manifolds which are obtained by removing boundary components from a compact 3-manifold.  The hope is that a better understanding of the Heegaard splittings of these ``household name" manifolds will aid in understanding how these manifolds compare to their more exotic cousins.  As Heegaard splittings lift to covering spaces, there may also be some hope of using these results in the study of compact Heegaard splittings.   We will eventually need to distinguish between two types of Heegaard splittings: relative and absolute, but for the initial statement of our results, we make the following preliminary definitions here.

\begin{definition}
Let $\mc{H}$ be the disjoint union of finitely or countably many 3-balls.  A \defn{handlebody} $H$ is formed by attaching finitely or infinitely many 1-handles $D^2 \times I$ to $\mc{H}$ so that a component of $D^2 \times \boundary I$ is attached to the boundary of a 3-ball.  We allow only finitely many 1-handles to be attached to each component of $\mc{H}$.    
\end{definition}

Handlebodies are characterized by the existence of a properly embedded collection of pairwise disjoint discs with boundary on $\boundary H$ which cut $H$ into 3-balls.  We will use the existence of these discs in many of the arguments in this paper. 

\begin{definition}
A \defn{Heegaard surface} in a 3-manifold $M$ with empty boundary is a properly embedded surface $S \subset M$ such that the closure of the complement of $S$ in $M$ consists of two handlebodies $U$ and $V$ with $S = \boundary U = \boundary V$.  We say that $M = U \cup_S V$ is a \defn{Heegaard splitting} of $M$.
\end{definition}

\begin{remark}
If $M$ has compact boundary, a Heegaard splitting will be a division of $M$ into two \defn{compressionbodies}.  We defer the definition of ``compressionbody" until after we have stated the main theorems of the paper.
\end{remark}

As in the theory of compact Heegaard splittings, stabilizations of Heegaard splittings play an important role.  The following are the most basic definitions.  More details will be given in Section \ref{definitions} and Section \ref{Exh. Seq.}.

\begin{definition}
A Heegaard splitting $M = U \cup_S V$ is \defn{stabilized} if there is an embedded 3-ball $B$ in $M$ such that $S \cap B$ is a properly embedded unknotted once-punctured torus.  The ball $B$ is a \defn{reducing ball} for $S$.
\end{definition}

\begin{definition}
A non-compact Heegaard splitting $M = U \cup_S V$ is \defn{end-stabilized} if for every compact set $C \subset M$ and for every non-compact component $W$ of the closure of the complement of $C$ in $M$ there is a reducing ball for $S$ which is contained in $W$.
\end{definition}

\subsection{Main Results}
We may now give simplified versions of the main results of this paper.  Throughout the paper we assume that $M$ is a non-compact orientable 3-manifold with no $S^2$ boundary components.  The first result is an elementary extension of a theorem of Frohman and Meeks \cite{FrMe97}.  It may be viewed as an analogue of the Reidemeister-Singer theorem for compact manifolds which says that, after some finite number of stabilizations, any two Heegaard splittings with the same partition of the boundary are ambient isotopic.  The proof of this extension is contained in the Appendix.

\begin{thmA}
Any two end-stabilized Heegaard splittings of $M$ which have the same partition of $\boundary M$ are properly ambient isotopic.
\end{thmA}

A (compact) Heegaard splitting is called \defn{weakly reducible} if there are disjoint compressing discs for the Heegaard surface which are contained in different compressionbodies.  In \cite{CaGo87}, Casson and Gordon prove that if a closed 3-manifold has a Heegaard splitting which is weakly reducible but not reducible then the manifold contains an incompressible surface of positive genus.  Since every non-compact Heegaard splitting (of non-zero genus) is weakly reducible, we cannot hope for a direct extension of the Casson and Gordon theorem to non-compact 3-manifolds.  The following result can, however, be viewed as a partial extension to the non-compact case.  The proof is based on a proof of the Casson-Gordon theorem.  Here is a simplified statement of the result.

\begin{definition}
A non-compact 3-manifold $M$ is \defn{end-irreducible} (rel $C$) for $C$ a compact subset if there is an exhausting sequence $\{K_i\}$ of nested compact connected submanifolds whose union is $M$ such that $C \subset \interior(K_1)$ and the frontier of each $K_i$ is incompressible in $M - C$.  $\{K_i\}$ is called a \defn{frontier-incompressible} (rel $C$) exhausting sequence for $M$.
\end{definition}

\begin{thmB}
Suppose that $M$ is end-irreducible (rel $C$) for some compact set $C$ containing $\boundary M$.  Let $M = U \cup_S V$ be a Heegaard splitting of $M$.  Then there is a frontier-incompressible rel $C$ exhausting sequence where each compact submanifold in the sequence ``inherits" a Heegaard splitting from $S$.
\end{thmB}

\begin{remark}
The notion of a compact submanifold ``inheriting" a Heegaard splitting from the non-compact 3-manifold $M$ will be made precise by using the idea of a ``relative Heegaard splitting".
\end{remark}

As an application of these two theorems we can classify Heegaard splittings of ``most" deleted-boundary 3-manifolds.  Here is an example:

\begin{thmC}
Let $\ob{M}$ be a compact, orientable 3-manifold with non-empty boundary, no component of which is a sphere.  Let $M$ be a 3-manifold obtained from $\ob{M}$ by removing one or more boundary components of $M$.  Then any two Heegaard splittings of $M$ which have the same partition of $\boundary M$ are equivalent up to proper ambient isotopy.  In particular, every Heegaard splitting of $M$ is of infinite genus and is end-stabilized.
\end{thmC}

\subsection{Acknowledgements}
Martin Scharlemann has provided many helpful comments and suggestions in the research leading to and on early drafts of this paper.  I am especially grateful for his sustained patience and encouragement.  Thanks also to Maggy Tomova, Ben Benoy, and Kelly Delp for our many conversations.  This research was partially supported by an NSF grant.

\subsection{Outline}

Section \ref{examples} contains several examples of non-compact Heegaard splittings and proves that the inclusion of a Heegaard surface into a non-compact 3-manifold induces a homeomorphism of ends. \newline

Sections \ref{slide-moves} and \ref{rel HS} provide preliminary work that is necessary for understanding the theorems which are the detailed versions of Theorems A, B, and C.  Section \ref{slide-moves} defines and studies handle-slides of boundary-reducing discs in compressionbodies.  Section \ref{rel HS} introduces a certain type of submanifold (due to Frohman-Meeks) which is ``balanced" on a non-compact Heegaard surface.  We discuss the type of Heegaard splittings (called ``relative Heegaard splittings") which these submanifolds inherit from the splitting of the manifold.  Both balanced submanifolds and relative Heegaard splittings are central in the work of Frohman and Meeks. \newline

Section \ref{Exh. Seq.} proves a more detailed version of Theorem B.  The proof is based on the outline of Casson and Gordon's theorem given in \cite{Sc02}. \newline

Theorem C is proved in Section \ref{Deleted Boundary}.  The key ideas are applications of Theorem A, Theorem B, and the classification of Heegaard splittings of $\text{(closed surface)} \times I$ by Scharlemann and Thompson \cite{ScTh93}. \newline

The Appendix proves Theorem A.  This is a slight extension of Frohman and Meeks' \cite[Theorem 2.1]{FrMe97} to manifolds with more than one end.  We also correct a misstatement\footnote{The error occurs in the last sentence of Prop. 2.2.  After including the collars of $J_i - L_i$ and $L_i - J_i$ you have arranged for $K_i$ to have a relative Heegaard splitting, but $K_{i+1} - K_i$ may not.  For example, the frontier of $K_i \cap H_1$ may not be incompressible in $H_1 \cap \cl(K_{i+1} - K_i)$.  This error affects the proof of Proposition 2.3.} in the proof of that theorem.  The correction is not difficult, but does require some work.

\subsection{History}\label{history}
Scharlemann, in his survey paper \cite{Sc02}, gives an overview of the history of Heegaard splittings of compact 3-manifolds.  We rely on that paper for our historical treatment of compact Heegaard splittings.  \newline

Very few types of compact 3-manifolds are known to have unique Heegaard splittings of a given genus and partition of the boundary.  Waldhausen \cite{Wa68b} proved that Heegaard splittings of a given genus of $S^3$ and of $S^2 \times S^1$ are unique up to ambient isotopy.  Scharlemann and Thompson in \cite{ScTh94a} have provided another proof of the classification for $S^3$.  Bonahon and Otal \cite{BoOt83} proved that lens spaces have unique splittings.  The 3-torus and $T^2 \times I$ also have unique splittings of a given genus and boundary partition.  This was proved by Boileau and Otal \cite{BoiOt90}.  Scharlemann and Thompson proved in \cite{ScTh93} that any two Heegaard splittings of $\text{(closed surface)} \times I$ with the same genus and partition of the boundary are isotopic.  In particular, any Heegaard splitting of $\text{(closed surface)} \times I$ is stabilized if it is of genus more than twice the genus of the surface and both boundary components are contained in a single compressionbody.  We will use their classification in Section \ref{Deleted Boundary}.  Schultens in \cite{Sch93} classified splittings of $\text{(compact surface)} \times S^1$.  She proved that any two splittings of $\text{(closed surface)} \times S^1$ of the same genus are isotopic.  The uniqueness of splittings for $\text{(closed surface)} \times I$ and $\text{(closed surface)} \times S^1$ provided the inspiration for Theorem C of this paper.  \newline

As previously mentioned, Frohman and Meeks \cite{FrMe97} in their work on Heegaard splittings of $\R^3$ defined the notion of noncompact Heegaard splitting used in this paper.  Pitts and Rubinstein \cite{PiRu86} have also considered Heegaard splittings of non-compact 3-manifolds.  They, however, consider only deleted boundary 3-manifolds and compact Heegaard surfaces which split the manifold into two ``hollow handlebodies".  For Pitts and Rubinstein, a hollow handlebody is simply a compact compressionbody with $\boundary_-$ removed.  Frohman and Meeks also use the term ``hollow handlebody", but they refer to what we call ``relative compressionbodies".  The term ``relative compressionbody" was introduced by Canary and McCullough \cite{CaMc04}.  Other authors have picked this term up and it has become standard.  This paper uses ``relative compressionbody"; this will, hopefully, avoid confusion with Pitts and Rubinstein's use of ``hollow handlebody".  Both Frohman-Meeks and Pitts-Rubinstein use Heegaard surfaces in non-compact 3-manifolds to study minimal surfaces from a topological point of view.  The main appearances of non-compact Heegaard splittings and non-compact handlebodies have been in minimal surface theory, for example \cite{Fr90,FrMe97, F92, MR05}.\newline

This paper focuses on ``eventually end-irreducible" 3-manifolds.  These manifolds were first studied by Brown and Tucker \cite{BrTu74}.  They are an important class of 3-manifolds since some questions about arbitrary non-compact 3-manifolds can be reduced to questions about eventually end-irreducible 3-manifolds \cite{BrTh87b}. 

\subsection{Definitions}\label{definitions}
\subsubsection*{Notation}
If $X$ is a subcomplex of a complex $Y$, then $\eta(X)$ denotes a closed regular neighborhood of $X$ in $Y$.  The term ``submanifold" will be reserved for codimension 0 submanifolds.  If $X$ is a submanifold of a manifold $Y$ then $\cl(X)$ indicates the closure of $X$ in $Y$ and $\interior(X)$ indicates the interior of $X$ in $Y$.  The number of components of a complex $X$ is denoted $|X|$.  The spaces $[0,1]$ and $[0,\infty)$ are denoted by $I$ and $\R_+$ respectively.  $\R^n$ denotes $n$-dimensional Euclidean space and $S^n$ denotes the sphere of dimension $n$.  The closed unit disc in $\R^2$ is denoted by $D^2$.  The integers and natural numbers are indicated by $\Z$ and $\N$ respectively.  All homology groups use $\Z$ coefficients.

\subsubsection*{3-Manifold Topology} We work in the PL category and use, with a few exceptions, standard terminology from 3-manifold theory (see \cite{He04,Ja80}).  All 3-manifolds and surfaces are assumed to be orientable.  A map $\rho: X \rightarrow Y$ between complexes is \defn{proper} if the preimage of each compact set is compact.  If $X$ is a surface and $Y$ is 3-manifold, $\rho$ is a \defn{proper embedding} if, in addition to being proper and an embedding, $\rho^{-1}(\boundary Y) = \boundary X$.  To say that a graph is properly embedded in a 3-manifold means that the inclusion map is proper and an embedding.  The vagaries of language being what they are, the terminology for surfaces and graphs in 3-manifolds is different, but this should not cause confusion in practice.  A homotopy $\rho: X \times I \rightarrow Y$ is \defn{proper} if it is proper as a map. If $X$ is a surface and $Y$ is a 3-manifold we also require that $\rho^{-1}(\boundary Y) = \boundary(X \times I)$. The homotopy $\rho$ is \defn{ambient}, if $X \subset Y$ and there is an extension of $\rho$ to a proper homotopy $\rho: Y \times I \rightarrow Y$.  An isotopy $\rho: X \times I \rightarrow Y$ is a homotopy where for each $t \in I$, $\rho(\cdot,t): X \rightarrow Y$ is an embedding.  An ambient isotopy $\rho: Y \times I \rightarrow Y$ is required to be a homeomorphism at each time $t \in I$.  To say that a homotopy $\rho$ is \defn{fixed} on a set $C$ means that $\rho$ restricted to $C \times \{t\}$ is the identity map for all $t \in I$.  \newline

A loop on a surface is \defn{essential} if it is not null-homotopic.  If the loop is embedded, this is equivalent to saying that it does not bound a disc on the surface. An embedded 2-sphere in a 3-manifold is \defn{essential} if it does not bound a 3-ball.  A \defn{compressing disc} for a surface $F$ in a 3-manifold is an embedded disc $D$ for which $D \cap F = \boundary D$ and $\boundary D$ is an essential loop on $F$.  A surface $F$ properly embedded in $M$ is \defn{incompressible} if there are no compressing discs for $F$ in $M$.  By the loop theorem and Dehn's lemma this is equivalent (since, for us, all surfaces and manifolds are orientable) to the inclusion map of $F$ into $M$ inducing an injection of fundamental groups.  Note that our definition considers an inessential 2-sphere to be an incompressible surface.  This is slightly non-standard, but it makes the statements and proofs of some of the results easier.  We will emphasize places where this observation matters.  \newline

If $S \subset M$ is a surface embedded in a 3-manifold and if $\Delta$ is the union of pairwise disjoint compressing discs for $S$ then $\sigma(S;\Delta)$ will denote the surface obtained from $F$ by compressing along $\Delta$.  If $R \subset S$ is a properly embedded subsurface (i.e. $\cl(R) = R$) with each component of $\boundary R$ either contained in or disjoint from $\boundary \Delta$ then $\sigma(R;\Delta)$ will denote the surface obtained from $R$ by compressing along those discs of $\Delta$ with boundary in $R$.  If $S \subset \boundary M$ then the manifold obtained by boundary-reducing $M$ along $\Delta$ is denoted $\sigma(M;\Delta)$.  As it will always be clear when we have a surface and when we have a 3-manifold this should not cause confusion.  \newline

A manifold (2 or 3 dimensional) is \defn{open} if it is non-compact and without boundary.  It is \defn{closed} if it is compact and without boundary.  A 3-manifold is \defn{irreducible} if every embedded 2-sphere bounds an embedded 3-ball.  As much as possible, we do not assume irreducibility.  A submanifold of a 3-manifold is a \defn{product region} if it is homeomorphic to $F \times I$ where $F$ is a surface.  A \defn{fiber} of $F \times I$ is $\{x\} \times I$ where $x \in F$.  A set $X \subset F \times I$ is \defn{vertical} if it is the union of fibers.

\subsubsection*{Heegaard Splittings}  For background on Heegaard splittings of compact 3-manifolds the reader is referred to the survey article \cite{Sc02}.  Since we are interested in splittings of non-compact 3-manifolds, some of our definitions differ from conventions in the compact setting.  Let $F$ be either a compact, orientable surface (possibly disconnected) or the empty set.  A \defn{compressionbody} $H$ is formed by taking the disjoint union of $F \times I$ and countably (finitely or infinitely) many disjoint 3-balls and then attaching 1-handles.  1-handles are attached to $F \times I$ on the interior of $F \times \{1\}$ and to the boundaries of the 3-balls.  Only finitely many 1-handles are to be attached to each 3-ball and only finitely many may be attached to $F \times \{1\}$.  We usually require that the result be connected.  The surface $F \times \{0\}$ is denoted $\boundary_- H$ and the surface $\cl(\boundary H - \boundary_- H)$ is denoted $\boundary_+ H$ and is called the \defn{preferred surface} of $H$.  If $F$ is a closed surface then $H$ is an \defn{absolute compressionbody}; if $F$ has non-empty boundary then $H$ is a \defn{relative compressionbody}.  If $H = F \times I$ then $H$ is a \defn{trivial compressionbody}.  If $F$ is empty, then $H$ is a \defn{handlebody}.  We will generally require that $F$ contain no $S^2$ components, as then $H$ is irreducible.  At one point, in section \ref{Exh. Seq.} we will need to allow $S^2$ components.  This will be explicitly pointed out.  A \defn{subcompressionbody} $A$ of $H$ is a submanifold of $H$ whose frontier in $H$ consists of properly embedded discs.  (We do not require these discs to be essential.  Thus, for example, upper half space, which is a handlebody, has an exhausting sequence consisting of subcompressionbodies.)  We denote $\boundary A \cap \boundary_+ H$ by $\boundary_{\boundary_+ H} A$.  There is a proper strong deformation retraction of a compressionbody $H$ onto $\boundary_- H \cup \Sigma$ where $\Sigma$ is a properly embedded graph in $H$ attached at valence one vertices to $\boundary_- H$.  $\Sigma \cup \boundary_- H$ is called the \defn{spine} of the compressionbody.  \newline

A compressionbody $H$ is determined by a properly embedded collection of disjoint discs $\Delta$ with $\boundary \Delta \subset \boundary_+ H$ such that $\sigma(H;\Delta)$ consists of 3-balls and $\boundary_- H \times I$.  Such a collection of discs is called a \defn{defining set of discs}.  A defining set of discs may contain discs which are not compressing discs.  An example is the handlebody which is a closed regular neighborhood of the positive $z$-axis in $\R^3$.  Since its boundary is homeomorphic to $\R^2$ there are no compressing discs, but there is clearly a defining set of discs.  A properly embedded collection of disjoint discs in $H$ with boundary on $\boundary_+ H$ will be called a \defn{disc set} for $H$ or for $\boundary_+ H$.  A disc set $\delta$ is \defn{collaring} if the union of some components of $\sigma(H;\delta)$ is $\boundary_- H \times I$.\newline

A \defn{Heegaard splitting} of a 3-manifold $M$ is a decomposition of $M$ into two compressionbodies $U$ and $V$ glued along $\boundary_+ U = \boundary_+ V$ = S.  If $U$ and $V$ are absolute compressionbodies the splitting is an \defn{absolute Heegaard splitting}.  If $U$ and $V$ are relative compressionbodies the splitting is a \defn{relative Heegaard splitting}.  The surface $S$ is called the \defn{Heegaard surface}.  We write $M = U \cup_S V$.  If the term ``Heegaard splitting" is used without either the adjective ``absolute" or ``relative", we will mean ``absolute Heegaard splitting".  Usually, relative Heegaard splittings will be of compact submanifolds of a non-compact 3-manifold. \newline

A Heegaard splitting of a manifold $M = U \cup_S V$ is \defn{reducible} if there is an essential simple closed curve on $S$ which bounds embedded discs in $U$ and $V$.  To \defn{stabilize} a Heegaard surface, push the interior of an embedded arc on the surface into one of the compressionbodies and include a regular neighborhood of the arc into the other compressionbody. A Heegaard splitting has been stabilized if there is, in $M$, an embedded 3-ball which intersects the Heegaard surface in a properly embedded, unknotted, once-punctured torus.  Such a ball is called a \defn{reducing ball}.  If a splitting other than the genus 1 splitting of $S^3$ is stabilized, it is reducible.  A Heegaard splitting is \defn{weakly reducible} (see \cite{CaGo87}) if there are disjoint compressing discs for the Heegaard surface, one in each compressionbody, which have boundaries that are essential on the Heegaard surface.  Any Heegaard splitting, other than the genus 0 splitting of $\R^3$, of a non-compact 3-manifold must have a properly embedded infinite collection of discs in each compressionbody which have boundaries which are essential on the Heegaard surface.  Thus, apart from the genus 0 splitting of $\R^3$, every non-compact Heegaard splitting is weakly reducible.  A Heegaard splitting $M = U \cup_S V$ is \defn{end-stabilized} if for every compact set $C \subset M$ and every non-compact component $W$ of $\cl(M - C)$ there is a reducing ball for $S$ entirely contained in $W$.   

\subsubsection*{Non-Compact 3-Manifolds}  An \defn{exhausting sequence} for a non-compact 3-manifold $M$ is a sequence $\{K_i\}$ of compact, connected 3-submanifolds such that $K_i \subset \interior(K_{i+1})$ and $M = \cup_i K_i$.  A 3-manifold $M$ is \defn{end-irreducible (rel $C$)} for a compact subset $C$ if there is an exhausting sequence for $M$ such that the frontier of each element of the exhausting sequence is incompressible in $M - C$.  If $C$ can be taken to be the empty set, then $M$ is simply \defn{end-irreducible}.  If $M$ is end-irreducible (rel $C$) for some $C$ then $M$ is \defn{eventually end-irreducible}. If a non-compact 3-manifold is obtained by removing at least one boundary component from a compact 3-manifold then the non-compact 3-manifold is a \defn{deleted boundary 3-manifold}.  Deleted boundary 3-manifolds are eventually end-irreducible.  All 3-manifolds (excluding compressionbodies) considered in this paper will have compact boundary.  When the manifold is end-irreducible (rel $C$) we will assume that $C$ contains $\boundary M$.  

\section{Examples}\label{examples}
Some examples of non-compact Heegaard splittings are in order.  When thinking about non-compact Heegaard splittings, keep in mind that an absolute handlebody is the closed regular neighborhood of a properly embedded, locally finite graph in $\R^3$.  Frohman and Meeks \cite{FrMe97} give an example of a non-compact 3-manifold whose interior is an open handlebody but the closure of the interior is not a handlebody.  Handlebodies have properly embedded discs which cut them into 3-balls.  Another observation, which may help the reader's intuition, is that no essential loop in an absolute compressionbody can be homotoped out of every compact set.  This is easily proved using the proper deformation retraction of the compressionbody to its spine.  This implies, for example, that $\text{(compact surface)} \times \R$ is not a handlebody.\newline

As mentioned in the opening paragraphs of the paper, Heegaard splittings of open manifolds can be obtained from triangulations or by lifting Heegaard splittings from manifolds which are covered by the manifold in question.

\subsection{Heegaard Splittings of $\R^3$}Heegaard splittings of $\R^3$ are easy to construct.  Since the upper and lower half spaces are homeomorphic to closed regular neighborhoods of the positive $z$-axis, $\R^3$ has a genus zero Heegaard splitting.  Obviously, this splitting can be stabilized any given (finite) number of times.  By choosing an infinite, properly embedded collection of arcs in the surface, it can also be stabilized an infinite number of times simultaneously to give an infinite genus Heegaard surface.  Frohman and Meeks prove that this is, up to proper ambient isotopy, the only infinite genus Heegaard splitting of $\R^3$.

\subsection{Finite Genus Heegaard Splittings} Let $\ob{M}$ be a compact 3-manifold with Heegaard splitting $\ob{U} \cup_{\ob{S}} \ob{V}$.  Let $B$ be an embedded closed 3-ball in $M$ which intersects $\ob{S}$ in a properly embedded disc.  Let $X = \ob{X} - B$ for $X = M, U, S, V$.  Then $M = U \cup_S V$ is a finite genus Heegaard splitting of the deleted boundary manifold $M$. (Infinitely many discs parallel to $\boundary B \cap \ob{U}$ ($\boundary B \cap \ob{V})$ are in any defining set of discs for $U$ ($V$).)  Classifying such Heegaard splittings would be equivalent to classifying all Heegaard splittings of compact manifolds.  No such simple classification is to be hoped for, and so our classification of Heegaard splittings for deleted boundary 3-manifolds does not address such examples.  Fortunately, this is the only type not covered by our classification.

\subsection{Amalgamating Heegaard Splittings}An easy way to create Heegaard splittings of non-compact manifolds is to amalgamate splittings of compact submanifolds.  We describe a way to do this, beginning with a description of amalgamation.  See \cite{Sch93} for the definition of amalgamation. Let $N_0$ and $N_1$ be two compact 3-manifolds with absolute Heegaard splittings $N_0 = U_0 \cup_{S_0} V_0$ and $N_1 = U_1 \cup_{S_1} V_1$ and homeomorphic boundary components $F_0 \subset \boundary_- V_0$ and $F_1 \subset \boundary_- V_1$.  Choose a homeomorphism $h:F_1 \rightarrow F_0$.  We can form a Heegaard splitting of the \defn{amalgamated manifold} $N = N_0 \cup_{h} N_1$ by \defn{amalgamation} in the following way: \newline

In $V_i$ there are collaring discs $\delta_i$ which cut off a product region $F_i \times I$ contained in $V_i$. Choose labels so that $F_i = F_i \times \{0\}$.  Let $P$ denote the product region $(F_1 \times I) \cup (F_0 \times I)$ in $N$. Think of $P$ as homeomorphic to $F_1 \times [0,2]$.  Note that it is contained in $V_0 \cup V_1$.  Perform an isotopy of $N_1$ so that, in $P$, $\delta_1 \times [0,2]$ is disjoint from $\delta_0$.   Let $A_1$ be $\delta_1 \times [0,2]$ in $P$.  Let $U = U_0 \cup (V_1 - (F_1 \times I)) \cup A_1$, $V = \cl(N - U)$ and $S = V \cap U$.  Then $M = U \cup_S V$ is a Heegaard splitting with genus equal to $\operatorname{genus}(S_0) + \operatorname{genus}(S_1) - \operatorname{genus}(F_1)$. Note that there there are disjoint discs $\delta_1 \subset U$ and $\delta_0 \subset V$ which, when we compress $S$ along them, leave us with a surface parallel to $F_0 = F_1$ in $N$.  \newline

Here is a method of producing an infinite genus Heegaard splitting of a non-compact 3-manifold $M$.  Let $\{K_i\}$ be an exhausting sequence for $M$ with the properties that $\boundary M \subset K_1$, that no component of $\cl(M - K_i)$ is compact for any $i$, and that for each $i$ and for each component $J$ of $\cl(K_{i+1} - K_i)$ the intersection $J \cap K_i$ is connected.  For each $i$, let $L_i = \cl(K_{i+1} - K_i)$ and $F_i = L_i \cap K_i$. $K_{i+1}$ is formed by amalgamating $K_i$ and each component of $L_i$ along a single component of the surface $F_i$. \newline

We now carefully choose absolute Heegaard splittings of $K_1$ and each component of $L_i$ for each $i \geq 1$.  Choose a Heegaard splitting $K_1 = U_1 \cup_{S_1} V_1$ of $K_1$ so that every boundary component of $K_1$ is contained in $V_1$.  Let $\delta'_1$ be a set of collaring discs for $V_1$.  Now for each component of $L_i$ choose a Heegaard splitting so that $L_i = X_i \cup_{T_i} Y_i$.  Choose the splitting so that $\boundary L_i \subset Y_i$.  Since each component of $L_i$ has at least two boundary components neither $X_i$ nor $Y_i$ has a component which is a stabilization of a trivial Heegaard splitting.  Inductively, form a Heegaard splitting of $K_n = U_n \cup_{S_n} V_n$ for $n \geq 2$ by amalgamating the Heegaard splittings of $K_{n-1}$ and $L_{n-1}$.  Let $V_n$ be the compressionbody which contains $\delta'_1$ and let $U_n$ be the other.  \newline

Recall from the definition of amalgamation that if $F_n \subset U_n$ then $U_{n+1} \cap U_n$ can be created by removing 1-handles in $U_n$ which join $F_i$ to $S_i$ and are vertical in the product structure of $U_n$ compressed along a defining set of discs.  Denote these 1-handles by $A_n$.  The surface $F_n$ is contained in $U_n$ whenever $n$ is even (by our choice of Heegaard splitting for $L_n$).  If $n$ is odd then $F_n$ is not in $U_n$, so for odd $n$, let $A_n = \nil$.  If $n$ is even then $U_n \subset U_{n+1}$.  Define $U'_n = \cl(U_n - A_n)$.  Since for each $n$, $\boundary L_n \subset Y_n$, $U'_n \subset U'_{n+1}$ for all $n$.  In particular, when we extend the 1-handles from $Y_n$ into $K_{n-1}$ they do not not need to reach into $L_{n-2}$.  Let $U = \cup_{\N} U'_n$. \newline

We desire to show that $U$ is an absolute compressionbody.  Since $\boundary M \subset V$, $U$ will be an absolute handlebody.  To prove this we will produce a properly embedded collection of discs in $U$ which cut $U$ into compact handlebodies.  Let $\delta_n$ be a collaring set of discs contained in $L_{n-1}$ for $U_n$ for each even $n$.  We may assume that $\delta_n$ is disjoint from $A_n$ and so $\delta_n$ is  a properly embedded finite collection of discs in $U$, for each even $n$.  Furthermore, since $\delta_n \subset L_{n-1}$ the infinite collection of discs $\delta = \cup \delta_n$ is properly embedded in $M$. The discs $\delta_n$ cut off a compact submanifold $U'_n - (\boundary K_{n+2} \times I)$.  As $U = \cup U'_n$ every component of $\sigma(U;\delta)$ is compact.  Let $H$ be a component of $\sigma(U;\delta)$.  Choose an even $n$ large enough so that $H \subset U'_{n-2}$.  $H$ is thus a component of $\sigma(U'_n;\delta)$ which is not contained in $\boundary U'_n \times I$.  As such, it must be a handlebody as $U'_n$ for $n$ even is an absolute compressionbody.  Hence $U$ is a handlebody. \newline

Let $V = \cl(M - U)$.  The argument to show that $V$ is an absolute compressionbody is similar, except that the disc set $\delta$ will cut $V$ into compact handlebodies and, if $\boundary M \neq \nil$, a compact absolute compressionbody $H$ with $\boundary_- H = \boundary M$.  Letting $S = U \cap V$, we have shown that $U \cup_S V$ is an absolute Heegaard splitting of $M$. \newline

It is instructive to examine this construction in the case when $M$ is a deleted boundary 3-manifold.  Let $M_0$ be a compact, orientable 3-manifold with non-empty boundary component $\boundary_1 M_0 = F \neq S^2$.  Let $M_0 = U_0 \cup_{S_0} V_0$ be a Heegaard splitting of $M_0$ with $F \subset V_0$.  Let $M_i$ for $i \geq 1$ be homeomorphic to $F \times I$ and choose a Heegaard splitting $M_i = U_i \cup_{S_i} V_i$ of $M_i$ which is obtained by tubing together two copies of $F$ in $M_i$.  Such a Heegaard splitting has both boundary components, $\boundary_0 M_i$ and $\boundary_1 M_i$, contained in $V_i$ and has genus which is twice the genus of $F$. (Heegaard splittings of $F \times I$ are classified by Scharlemann and Thompson in \cite{ScTh93}.  This classification will be important for our work in Section \ref{Deleted Boundary}.)  Build a 3-manifold $M$, homeomorphic to $M_0 - F$ by glueing $\boundary_0 M_i$ to $\boundary_1 M_{i-1}$ for $i \geq 1$.  At stage $n$ of the glueing process we can obtain a Heegaard splitting of the new manifold by amalgamating the splittings of the previously constructed manifold and $M_n$.  The new Heegaard splitting will have genus equal to $\operatorname{genus}(S_0) + n\cdot\operatorname{genus}(F)$.  This produces an infinite genus splitting of $M$.  It is easy to verify that the splitting is end-stabilized.  The content of Proposition \ref{no spheres} is that, up to proper ambient isotopy, this is the only Heegaard splitting of $M$. 

\subsection{Infinite Genus Splittings which are not End-stabilized}
A consequence of Theorem C is that all infinite genus splittings of one-ended deleted boundary 3-manifolds are end-stabilized. It is then natural to ask:

\begin{question}
Are there examples of one-ended, irreducible 3-manifolds which have infinite genus Heegaard splittings that are not stabilized?  Are there such examples where the manifold has finitely generated fundamental group?  What if we simply require that the splitting not be end-stabilized?
\end{question}

In this subsection, we give two examples of splittings which are not end-stabilized.  The first example is a non-stabilized splitting of a one-ended, irreducible 3-manifold $M$ with infinitely generated fundamental group.  The second example, which is obtained from the first, is a splitting of the Whitehead manifold $W$ which is not end-stabilized, but which is stabilized and cannot be destabilized finitely many times.  The key point is that, though there are infinitely many ``inequivalent" reducing balls, they are not properly embedded in $W$.  In this sense, this example is similar in spirit to Peter Scott's construction of a simply connected 3-manifold which is not the connect sum of prime 3-manifolds \cite{Sc77}.  I do not know of a one-ended, irreducible manifold with finitely generated fundamental group which has an infinite genus non-stabilized splitting.\newline

We begin by constructing the splitting of $M$.  Let $W_0$ be the exterior of the Whitehead link in $S^3$.  $W_0$ is a compact 3-manifold which contains no essential annuli or essential tori\footnote{This is easy to prove directly, or see \cite{EudUch96}.}.  $W_0$ is hyperbolic (Example 3.3.9 of \cite{Thur97}).  As the Whitehead link is a 2-bridge link, it has tunnel number one, and therefore $W_0$ has a genus 2 Heegaard splitting which does not separate $\boundary W_0$.  Let $\boundary_0 W_0$ and $\boundary_1 W_0$ be the two boundary components of $W_0$.  Let $\lambda_j$ and $\mu_j$ be the longitude and meridian of $\boundary_j W_0$ (for $j = 0,1$).  The choice should be made so that $\lambda_j$ and $\mu_j$ correspond to the longitude and meridian of the corresponding component of the Whitehead link in $S^3$.  In particular, $\lambda_0$ and $\lambda_1$ are homologically trivial in $W_0$ and $\mu_0$ and $\mu_1$ included into $W_0$ generate the first homology of $W_0$.  Let $f:\boundary_0 W_0 \rightarrow \boundary_1 W_{0}$ be a homeomorphism which takes $\lambda_0$ to $\mu_1$ and $\mu_0$ to $\lambda_1$. \newline

For each $i \in \N$ let $W_i$ be a copy of $W_0$.  Denote the boundary components of $W_i$ by $\boundary_0 W_i$ and $\boundary_1 W_i$ in such a way that the labelling corresponds to the labelling of the boundary components of $W_0$.  Let $S_i$ be a genus 2 Heegaard surface for $W_i$ which does not separate the boundary components. Let $f_i: \boundary_0 W_i \rightarrow \boundary_1 W_{i-1}$ be the map $f$.  Let $M_1 = W_1$ and, inductively, let $M_n = M_{n-1} \cup_{f_{n}} W_{n}$, $\boundary_0 M_n = \boundary_0 W_1$, and $\boundary_1 M_n = \boundary_1 W_n$ for $n \geq 2$.  Let $S'_n$ be the Heegaard surface of $M_n$ and $S$ the Heegaard surface of $M$ obtained by amalgamating the surfaces $S_i$, as described previously. \newline

We now show that $S$ is not stabilized.  If it were, then some $S'_n$ would be stabilized, as reducing balls are compact.  Without loss of generality, we may assume that $n$ is odd, so that $S'_n$ does not separate the boundary components of $M_n$.  It will be beneficial to work with a closed 3-manifold: glue a copy of $W_0$ to $M_{n}$ to obtain a closed 3-manifold $M'$.  Use the glueing maps $f:\boundary_0 W_0 \rightarrow \boundary_1 M_{n}$ and $f^{-1}:\boundary_1 W_0 \rightarrow \boundary_0 M_n$.  We may form a Heegaard splitting of $M'$ by amalgamating a genus 2 splitting, which does not separate $\boundary W_0$, of $W_0$ to $S'_{n}$ across $\boundary M_{n}$ to obtain a Heegaard surface $T$.  As neither splitting separates the boundary components of the respective manifolds, this operation gives a well-defined Heegaard splitting $T$ of $M'$, a closed 3-manifold.  The genus of $S'_{n}$ is $(2n - (n - 1)) = n + 1$.  The splitting given by $T$ is obtained from $S'_n$ by adding a single one-handle to the handlebody in the splitting of $M_n$.  Thus, the genus of $T$ is one more than the genus of $S'_n$; that is, the genus of $T$ is $n+2$.  By assumption, $S'_n$ is stabilized, and so $T$ is, as well.  Thus, $M'$ has an irreducible Heegaard splitting of genus $g \leq n+1$. \newline

We now appeal to a theorem of Scharlemann and Schultens.  A consequence of Theorem 4.7 of \cite{ScSch01} is that if $M'$ (a closed, orientable, irreducible 3-manifold) has a JSJ-decomposition with $q$ non-Seifert fibered submanifolds,  then $q \leq g-1$.  Let $\Theta$ be the union of the boundary tori of $W_i$ for $i \leq n$.  As each $W_i$ contains no essential annuli or tori, $\Theta$ is the union of the canonical tori in the JSJ-decomposition of $M'$.  None of the $W_i$ are Seifert fibered, so $q = n + 1$.  Therefore, $q = n+1 \leq g - 1 \leq n$, a contradiction.  We conclude that $S$ is not stabilized. \newline

We have just shown that the manifold $M$ has an infinite genus Heegaard surface $S$ which is not stabilized.  $M$ has infinitely generated fundamental group as the tori $\boundary W_i$ for $i \geq 1$ are all incompressible and non-parallel.  

\begin{remark}
There are many other similar constructions of 3-manifolds with infinitely generated fundamental group which have non-stabilized splittings.  By allowing arbitrary glueing maps between boundary tori of the compact pieces, one can use a theorem of Bachman, Schleimer, and Sedgewick\cite{BaScSe05} to show that the amalgamated splittings are not stabilized.  We do not pursue this route further in this paper.
\end{remark}

We now use the splitting of $M$ to obtain a splitting of the Whitehead manifold $W$.  The manifold $M$ has a single boundary component $\boundary_0 W_1$.  By attaching a solid torus $V$ to $\boundary M$ so that the meridian of the solid torus is equal to the meridian $\mu_0$ of $\boundary_0 W_1$, we obtain the Whitehead manifold $W \supset M$.  As this same process can be achieved by attaching first a 2-handle and then a 3-ball to $\boundary M$, the surface $S$ is still a Heegaard surface for $W$.  As $S$ is not stabilized in $W - V$, every stabilizing ball of $S$ in $W$ must intersect the compact set $V$.  Thus, $S$ is not end-stabilized.  $S$ is, however, stabilized.  To see this, recall that $S$ is formed by amalgamating the splitting $S'_n$ (for any given $n$) to the splittings $S_i$ for $i \geq n$.  Interpreted in $W$, $S'_n$ (for any $n$) is a splitting of a solid torus.  The genus of $S'_n$ is $n + 1$ and so, by the classification of splittings of handlebodies, $S'_n$ can be destabilized $n$ times in $W$.  This means, then, that $S$ can be destabilized infinitely many times in $W$.  The stabilizing balls are not properly embedded in $W$ and so only finitely many destabilizations can occur at once. 

\subsection{Ends of Heegaard Surfaces} The remainder of this section is devoted to showing that the inclusion of a Heegaard surface into $M$ induces a homeomorphism of end spaces.  It may serve to give the reader some feel for the properties of noncompact Heegaard splittings. Before stating the results, we recall the definition of the set of ends of a manifold. (See, for example, \cite{BrTu74}.)

\begin{definition}
A \defn{ray} in a connected manifold $M$ is a proper map\linebreak[4] $r:\R_+ \rightarrow M$.  An \defn{end} of a non-compact manifold $M$ is an equivalence class of rays.  Two rays $r,s:\R_+ \rightarrow M$ are equivalent if for every compact set $C \subset M$ there is a number $t_C \in \R_+$ such that the images of $[t_C,\infty)$ under $r$ and under $s$ are in the same component of $M - C$.  The set of ends is topologized by declaring that for any compact set $C$ and any non-compact component $A$ of the closure of $M - C$ the set of equivalence classes $\{[r] : \exists t \in \R_+ \text{ with } r([t,\infty)) \subset A\}$ is an open set.  These open sets form a basis for the topology on the end space of $M$. The set of ends of $M$ with this topology is 0-dimensional, compact, and Hausdorff \cite{Ra60}.
\end{definition}

The proofs of the following lemma and proposition follow suggestions by Martin Scharlemann.

\begin{lemma}\label{graph exterior}
Let $\Gamma$ be a locally finite graph properly embedded in an open 3-manifold $M$.  Then the inclusion of $M - \interior(\eta(\Gamma))$ into $M$ induces a homeomorphism of ends.
\end{lemma}

\begin{proof}
Let $X = M - \interior(\eta(\Gamma))$.  Let $r$ and $s$ be two rays determining the same end of $X$.  Let $C \subset M$ be a compact set.  $X$ is a closed subset of $M$.  As such, $C \cap X$ is a compact subset of $X$.  Hence, there exists a $t \in [0,\infty)$ such that the images of $[t,\infty)$ under $r$ and $s$ are contained in the same component of $X - C$.  This means that the images of $[t,\infty)$ under $r$ and $s$ are contained in the same component of $M - C$.  Thus, $r$ and $s$ are rays in $M$ and determine the same end of $M$.  Hence, there is a well-defined map on ends induced by the inclusion of $X$ into $M$. \newline

We next prove that the induced map on ends is surjective.  Suppose that $[r]$ is an equivalence class of ends of $M$.  By general position, there is a representative of this equivalence class which is disjoint from $\Gamma$ and, hence, there is a representative $r$ which is contained in $X$.  Under the induced map the equivalence class $[r]$ in the set of ends of $X$ is sent to the equivalence class $[r]$ in the set of ends of $M$.  Thus, the induced map on ends is surjective. \newline

Now suppose that $[r]$ and $[s]$ are equivalence classes in the set of ends of $X$ which have the same image in the set of ends of $M$ under the map induced by the inclusion of $X$ into $M$.  Let $r$ and $s$ be representatives of these equivalence classes in the set of ends of $X$.  Since $r$ and $s$ represent the same equivalence class in the set of ends of $M$, for any compact set $C \subset M$ there is a $t_C \in [0, \infty)$ such that the images of $[t_C,\infty)$ under $r$ and $s$ are contained in the same component of $M - C$.  Let $K \subset X$ be a compact set.  As $X$ is closed in $M$, $K$ is a compact subset of $M$.  The images $r([t_K,\infty))$ and $s([t_K,\infty))$ are contained in the same component of $M - K$.  The components of $M - K$ are also the path components of $M - K$, so there is a path $\gamma$ contained in $M - K$ joining $r([t_K,\infty))$ and $s([t_K,\infty))$.  By general position, we may homotope $\gamma$ so that its image is contained in $M - (K \cup \eta(\Gamma))$. That is, $\gamma$ is a path in $X - K$ joining $r([t_K,\infty))$ and $s([t_K,\infty))$.  
Thus, $r([t_K,\infty))$ and $s([t_K,\infty))$ are contained in the same component of $X - K$.  Since $K$ was an arbitrary compact subset of $X$, $[r] = [s]$ in the set of ends of $X$ and the induced map on ends is injective. \newline

We now prove that the induced map is bicontinuous.  To show continuity, it suffices to show that the preimage of a basis element in the topology of ends of $M$ is open in the ends of $X$.  Let $A'$ be a basis element in the topology of the set of ends of $M$.  By definition, there is a compact set $C \subset M$ and a non-compact component $A$ of $M - C$ such that for each ray $r$ for which $[r] \in A'$ there is $t_r \in [0,\infty)$ such that $r([t_r,\infty))$ is contained in $A$.  By replacing $C$ with $\eta(C)$, we may assume that $C$ and $A$ are submanifolds of $M$.  Since $X$ is closed in $M$, $C \cap X$ is compact and so by choosing representatives $r$ for each $[r] \in A'$ such that $r$ is a ray in $X$, we see that $r([t_r,\infty))$ is contained in $A \cap X$. \newline

We claim that $A \cap X$ is connected and non-compact.  It is easy to see that $A \cap X$ is path-connected: choose two points $x,y \in A \cap X$.  Since $A$ is path-connected, there is a path in $M$ joining them.  By general position we may assume that the path is disjoint from $\Gamma$.  Thus, there is a path in $A$ disjoint from $\eta(\Gamma)$.  Hence, $A \cap X$ is path-connected and therefore connected.  $A \cap X$ is also non-compact since $r$ is a proper map and the image of $[t_r,\infty)$ under $r$ is contained in $A \cap X$.  The preimage of $A'$ is, therefore, contained in the set $A'' = \{[s]: \exists t \in \R_+ \text{ with } s([t,\infty)) \subset (A \cap X)\}$.  Suppose, now, that $s$ is a representative for $[s] \in A''$.  Since $A \cap X \subset A$, $s([t,\infty)) \subset A$.  Thus, the image of $[s]$ under the inclusion map of ends of $X$ into ends of $M$ is contained in $A'$.  Thus, the preimage of $A'$ is $A''$.  $A''$ is, by definition, an open set in the topology of the set of ends of $X$.  Hence, the induced map on ends is continuous.  Since the set of ends of a connected manifold is compact and Hausdorff the induced map also has continuous inverse. Thus, the induced map is a homeomorphism.
\end{proof}

\begin{proposition}\label{end homeomorphism}
Let $M = U \cup_S V$ be an absolute Heegaard splitting of a non-compact manifold with compact boundary.  Then the inclusion of $S$ into $M$ induces a homeomorphism of ends.
\end{proposition}

\begin{proof}
If $\boundary M \neq \nil$ we can attach finitely many 2 and 3-handles to $\boundary M$ to obtain an open 3-manifold $M'$ containing $M$.  An absolute Heegaard splitting for $M$ is also a Heegaard splitting for $M'$, since the 2 and 3-handles were attached to $\boundary_-$ of the compressionbodies.  Since we attached only finitely many 2 and 3-handles, the inclusion of $M$ into $M'$ induces a homeomorphism of ends.  So, without loss of generality, we may assume that $M$ is open. \newline

Choose spines $\Sigma_U$ and $\Sigma_V$ for $U$ and $V$ respectively.  Let $\Gamma = \Sigma_U \cup \Sigma_V$.  $\Gamma$ is a locally finite graph properly embedded in $M$.  Let $X$ be the complement of an open regular neighborhood of $\Gamma$ in $M$.  Since $\Sigma_U$ and $\Sigma_V$ are spines of handlebodies giving a Heegaard splitting of $M$, $X$ is homeomorphic to $S \times I$.  By Lemma \ref{graph exterior}, the inclusion of $X$ into $M$ induces a homeomorphism of ends.  Since $X$ is homeomorphic to $S \times I$ there is a proper deformation retraction of $X$ onto $S \times \{\frac{1}{2}\}$.  Thus the inclusion of $S$ into $X$ is a proper homotopy equivalence and so induces a homeomorphism on ends.  Therefore, the inclusion of $S$ into $M$ induces a homeomorphism of ends.
\end{proof}

\begin{remark}
In \cite{FrMe97}, Frohman and Meeks prove by algebraic means that a Heegaard surface in a 1-ended 3-manifold is 1-ended.
\end{remark}
\section{Slide-Moves}\label{slide-moves}

\subsection{Handle-slides}\label{handle-slides}
Let $H$ be a compressionbody (absolute or relative) with preferred surface $S = \boundary_+ H$.  Suppose that we are given a disc set $\Delta$ for $H$ (with $\boundary \Delta \subset \boundary_+ H$).  We now describe a process which transforms $\Delta$ into a new disc set $\Delta'$. \newline

Let $\alpha \subset \boundary_+ H$ be an oriented arc such that $\alpha \cap \boundary \Delta = \boundary \alpha$.  Suppose that the endpoints of $\alpha$ are on distinct discs of $\Delta$.  Let $D_1$ and $D_2$ be the discs of $\Delta$ containing $\boundary \alpha$ so that $\alpha$ joins $D_1$ to $D_2$.  A regular neighborhood of $D_1 \cup \alpha \cup D_2$ has frontier in $H$ consisting of three discs.  Two of these discs are parallel to $D_1$ and $D_2$, the other has arcs in its boundary which are subarcs of $\eta(\alpha)$.  Let $D_1 \slide{\alpha} D_2$ denote this disc.  Let $\Delta' = (\Delta - D_1) \cup (D_1 \slide{\alpha} D_2)$.

\begin{definition} 
The disc set $\Delta'$ is obtained from $\Delta$ by a \defn{handle-slide} of $\Delta$ along $\alpha$.  If $D_1, D_2$ and $\alpha$ are all disjoint from a closed set $X$ then the handle-slide is said to be done \defn{relative} to $X$ or \defn{(rel $X$)}.
\end{definition}

Suppose that $A \subset H$ is a subcompressionbody with the property that $\fr A \subset \Delta$.  There is a subcompressionbody $A'$ of $H$ with frontier contained in $\Delta'$ which we say is \defn{obtained from $A$ by a handle-slide}.  The definition of $A'$ depends on the location of $D_1$ and $\alpha$:
\begin{itemize}

\item If $D_1$ is not in the frontier of $A$ then $A'$ is equal to $A$.

\item If $D_1$ is in the frontier of $A$ and $\alpha$ is contained in $\boundary_{S} A$ then we remove the interior of a collar neighborhood of $\alpha \cup D_2$ from $A$.  (The neighborhood of $D_2$ should be taken to be just on the side of $D_2$ which $\alpha$ intersects.  This way, if $D_2 \subset \interior A$, the disc $D_2$ itself is not removed.) If $D_2$ wasn't in the frontier of $A$, it is now contained in $\fr A'$.  

\item If $D_1$ is in the frontier of $A$ and $\alpha$ is not contained in $\boundary_{S} A$ then to form $A'$, we add the closure of a regular neighborhood of $\alpha \cup D_2$ to $A$.  (Again, the neighborhood of $D_2$ should be taken to be just on the side of $D_2$ which intersects $\alpha$.) 
\end{itemize}

\begin{remark}
The subcompressionbodies $A$ and $A'$ may not be homeomorphic (if, for example, both $D_1$ and $D_2$ are contained in $\fr A$ and $\alpha$ is not in $\boundary_S A$).  We do have, however, that $\sigma(A;\Delta)$ is homeomorphic to $\sigma(A';\Delta')$. 
\end{remark}

Likewise, if $R$ is a (topologically) closed subsurface of $\boundary_+ H$ with the following three properties:
\begin{enumerate}
\item[$\cdot$] $\boundary D_1$ is either a component of $\boundary R$ or disjoint from $\boundary R$.

\item[$\cdot$] $\boundary D_2$ is either a component of $\boundary R$ or disjoint from $\boundary R$.

\item[$\cdot$] The interior of $\alpha$ is disjoint from $\boundary R$
\end{enumerate}
then we can form a new surface $R'$ which is \defn{obtained from $R$ by a handle-slide}.  If $\boundary D_1 \cap \boundary R = \nil$ then $R'$ is defined to be $R$.  If $\boundary D_1 \subset \boundary R$ and $\alpha \subset R$ then $R'$ is defined to be $\cl(R - \eta(\alpha \cup \boundary D_2))$ where the neighborhood of $\boundary D_2$ is a one-sided neighborhood on the side of $D_2$ which $\alpha$ meets.  This way if $\boundary D_2 \subset \interior(R)$ then $\boundary D_2 \subset \boundary R'$. If $\boundary D_1 \subset \boundary R$ and $\alpha$ is not contained in $R$ then $R'$ is defined to be $R \cup \eta(\alpha \cup \boundary D_2)$.  As before, the neighborhood of $\boundary D_2$ should be taken to be a one-sided neighborhood on the side of $\boundary D_2$ which $\alpha$ meets.\newline

Bonahon developed the use of handle-slides to prove results about compressionbodies.  The following proposition and its corollaries are based on his work in \cite{Bo83}.  For proofs see Appendix B of that paper.  The essence of the proof of Proposition \ref{choosing discs} shows up in Step 6 of the proof of Proposition \ref{good discs exist} of this paper.

\begin{proposition}\label{choosing discs}
If $D$ is a boundary-reducing disc for $H$ then there is a collection of defining discs for $H$ which are disjoint from $D$.
\end{proposition}

\begin{corollary}\label{boundary reducing gives compressionbody}
Boundary-reducing a compressionbody along a finite disc set results in compressionbodies.
\end{corollary}

\begin{corollary} \label{nested reducing discs}
Given any finite disc set for a compressionbody, there is a defining collection of discs for the compressionbody which contains the given disc set.
\end{corollary}

\begin{corollary}
A subcompressionbody with compact frontier is a compressionbody.
\end{corollary}

The following definition will be useful later.  We include it here since Corollary \ref{complementary compressionbodies} follows from Corollary \ref{nested reducing discs}.

\begin{definition}
If $A$ and $B$ are relative compressionbodies with $A \subset B$, we say that $A$ is \defn{correctly embedded} in $B$ if $\boundary_+ A \subset \boundary_+ B$ and if every closed component of $\boundary_- A$ is also a component of $\boundary_- B$.
\end{definition}

Another way of stating the definition is that $A \subset B$ is correctly embedded if each component of $\fr A$ is a component of $\boundary_- A$ which has non-empty boundary and is properly embedded in $B$.

\begin{remark}
The notion of ``correctly embedded" is similar to Canary and McCullough's ``normally imbedded" in \cite[Section 3.4]{CaMc04}.
\end{remark}

\begin{corollary}\label{complementary compressionbodies}
Suppose that $A$ is a compact relative compressionbody correctly embedded in a relative compressionbody $B$.  Then $\cl(B - A)$ is a relative compressionbody. In particular, if each closed component of $\boundary_- B$ is contained in $A$ then $\cl(B - A)$ is a handlebody.
\end{corollary}

\begin{proof}
Choose a defining set of discs $\Delta_A$ for $A$.  Boundary-reduce $B$ along $\Delta_A$ to obtain $B'$.  Corollary \ref{boundary reducing gives compressionbody} implies that each component of $B'$ is a (relative) compressionbody.  Each component of $B'$ was either contained in $A$ or contains a copy of $\fr A \times I$ with $\fr A \times \{1\}$ a subsurface of $\boundary B'$.  Subtracting $\fr A \times I$ from those components of $B'$ is simply removing a collar of a subsurface of $\boundary B'$ from $B'$ and hence leaves us with a compressionbody.  But this is exactly $\cl(B - A)$.  If each closed component of $\boundary_- B$ is contained in $A$ then, if $C$ is a component of $B'$ which is not contained in $A$, $\boundary_- C$ contains no closed components.  By our definition of ``compressionbody", $\boundary C$ is compact and so $C$ is formed by adding one-handles to $F \times I$ where $F$ is a compact surface, no component of which is boundary-less.  Thus, $C$ is obtained by adding one-handles to handlebodies and, so, is a handlebody.  We then form a component of $\cl(B - A)$ by removing a neighborhood of $\fr A \cap C$ from $C$.  The result is homeomorphic to $C$ and so is a handlebody. 
\end{proof}

\begin{remark}
We may not be able to choose $\cl(\boundary_+ B - \boundary_+ A)$ to be the preferred surface of $\cl(B - A)$.  For example, if $B$ is a compact relative compressionbody and we push each non-closed component of $\boundary_- B$ slightly into $B$ and take $A$ to be the closure of the complement of the product regions.
\end{remark}

\subsection{Slide-Moves and Isotopies}\label{slide-moves and isotopies}

For the remainder of this section, let $S$ be an absolute Heegaard surface dividing a 3-manifold $M$ with compact boundary into absolute compressionbodies $U$ and $V$. \newline

If we have disc sets $\ob{\Delta}_1$ for $U$ and $\ob{\Delta}_2$ for $V$ which are disjoint from each other we can perform handle-slides on each disc set individually.  The remainder of this section studies how these handle-slides affect the surface $S$.

\begin{definition}
A \defn{2-sided disc family} $\ob{\Delta}$ for $S$ in $M$ is the union of disc sets $\ob{\Delta}_1$ and $\ob{\Delta}_2$ for $U$ and $V$ with the property that the discs of $\ob{\Delta} = \ob{\Delta}_1 \cup \ob{\Delta}_2$ are pairwise disjoint.
\end{definition}

We can expand the notion of a handle-slide to that of a slide-move on the 2-sided disc family $\ob{\Delta} = \ob{\Delta}_1 \cup \ob{\Delta}_2$:

\begin{definition}
A \defn{slide-move} of $\ob{\Delta}$ is one of the following operations:
\begin{enumerate}
\item[(M1)] Perform a handle-slide (rel $\boundary \ob{\Delta}_2$) of $\ob{\Delta}_1$.
\item[(M2)] Add to $\ob{\Delta}_1$ a boundary-reducing disc for $U$ which is disjoint from $\ob{\Delta}_1 \cup \ob{\Delta}_2$.
\item[(M3)] Perform a handle-slide (rel $\boundary\ob{\Delta}_1$) of $\ob{\Delta}_2$.
\item[(M4)] Add to $\ob{\Delta}_2$ a boundary-reducing disc for $V$ which is disjoint from $\ob{\Delta}_1 \cup \ob{\Delta}_2$.
\end{enumerate}
\end{definition}

Suppose that $A$ is a subcompressionbody of $U$ or $V$ with $\fr A \subset \ob{\Delta}$.  If we perform a slide-move on $\ob{\Delta}$ to obtain a 2-sided disc family $\Delta$ we can obtain from $A$ a subcompressionbody $A'$ with frontier contained in $\Delta$: If slide-move (M2) or (M4) is performed, $A'$ is defined to be equal to $A$.  If $A \subset U$ and slide-move (M1) is performed, $A'$ is defined to be the subcompressionbody obtained from $A$ by the handle-slide (see Section \ref{handle-slides}).  Similarly, if $A \subset V$ and slide-move (M3) is performed, $A'$ is defined to be the subcompressionbody obtained from $A$ by the handle-slide.  If we perform a finite sequence of slide-moves to obtain $\Delta$ from $\ob{\Delta}$ there is a subcompressionbody $A'$ with $\fr A' \subset \Delta$ obtained from $A$ by a finite number of handle-slides.  We say that $\Delta$ is \defn{obtained from $\ob{\Delta}$ by slide-moves} and that $A'$ is \defn{obtained from $A$ by slide-moves}.  \newline

Suppose that $R$ is a properly embedded subsurface of $S$ with $\boundary R \subset \boundary \ob{\Delta}$.  The boundary components of $R$ may bound discs in either $U$ or $V$ (i.e. discs which are in $\ob{\Delta}_1$ or $\ob{\Delta}_2$).  If we perform a finite sequence of slide-moves on $\ob{\Delta}$ to obtain $\Delta$ we may define a subsurface $R'$ of $S$ which is \defn{obtained from $R$ by slide-moves} and has boundary contained in $\boundary \Delta$. The definition is basically the same as the definition when a single handle-slide is performed: If slide-moves (M2) or (M4) are performed, $R$ is left unchanged.  If (M1) or (M3) is performed, so that a disc $D_1$ is slid over a disc $D_2$ via a path $\alpha$, we can define $R'$ as before (see Section \ref{handle-slides}).  The following proposition is an integral part of the proof of Theorem \ref{well-placed}.

\begin{proposition}\label{isotopies of slide-moves}
Suppose that $\ob{\Delta}$ is a 2-sided disc family for $S$ and that $\Delta$ is obtained from $\ob{\Delta}$ by slide-moves.  Then there is a finite collection of disjoint discs $\mc{D}$ with $\boundary \mc{D} \subset \sigma(S;\ob{\Delta})$ and a proper ambient isotopy of $\sigma(S;\Delta)$, fixed outside a compact subset of $M$, with the following properties:
\begin{enumerate}
\item[(i)] The discs $\mc{D} = D_1 \cup \hdots \cup D_p$ have an ordering such that the disc $D_i$ intersects only on its boundary the surface $\sigma(S;\ob{\Delta})$ compressed along $D_1, \hdots, D_{i-1}$. (See the remark below.)
\item[(ii)] The isotopy takes $\sigma(S;\Delta)$ to $\sigma(S;\ob{\Delta})$ compressed along $\mc{D}$.
\item[(iii)]  Let $R$ be a topologically closed subsurface of $S$ such that $\boundary R \subset \boundary \ob{\Delta}$ and $R'$ the subsurface of $S$ obtained from $R$ by that sequence of slide-moves.  The isotopy takes $\sigma(R';\Delta)$ to the surface obtained from $\sigma(R;\ob{\Delta})$ by compressing along whatever discs of $\mc{D}$ have boundary in $R$. 
\end{enumerate}
\end{proposition}

\begin{remark}
The discs $\mc{D}$ may intersect $S$ on their interiors, so part of the conclusion of the theorem is that when we compress $\sigma(S;\ob{\Delta})$ along the discs $D_1, \hdots, D_{i-1}$ we have chosen the regular neighborhoods of $D_1,\hdots, D_{i-1}$ so that $D_i$, although it may intersect $S$, does not intersect $\sigma(S;\ob{\Delta})$ compressed along $D_1, \hdots, D_{i-1}$.  We will abuse notation and write $\sigma(S;\ob{\Delta} \cup \mc{D})$ for the surface obtained from $\sigma(S;\ob{\Delta})$ by compressing along the discs $\mc{D}$ in the order given.  Similarly, if $R$ is a topologically closed subsurface of $S$ with $\boundary R \subset \boundary \ob{\Delta}$ we will use $\sigma(R;\ob{\Delta} \cup \mc{D})$ to indicate the surface obtained from $R$ by compressing along the discs of $\ob{\Delta}$ and then $\mc{D}$ in the given order (rather, compressing along those discs which have boundary on $R$).
\end{remark}

The proof of Proposition \ref{isotopies of slide-moves} will make use of the following lemma:

\begin{lemma}\label{isotopies of handle slides}
If $\ob{\Delta}_i$ is a disc family for $S$ with $\ob{\Delta}_i \subset U$ or $\ob{\Delta}_i \subset V$ and if $\Delta_i$ is obtained from $\ob{\Delta}_i$ by a single handle-slide of the disc $D_1$ over the disc $D_2$ via a path $\alpha$, then there is a proper ambient isotopy of $M$, fixed off a compact set, with the following properties:
\begin{enumerate}
\item[(a)] the isotopy takes $\sigma(S;\Delta_i)$ to $\sigma(S;\ob{\Delta}_i)$.
\item[(b)] if the handle-slide is relative to a closed set $X$ then we can choose the isotopy to be (rel $X$).
\item[(c)] if $R$ is a subsurface of $S$ with all of the following properties:
\begin{enumerate}
\item[$\cdot$]$\boundary D_1$ is either a component of $\boundary R$ or disjoint from $\boundary R$.
\item[$\cdot$]$\boundary D_2$ is either a component of $\boundary R$ or disjoint from $\boundary R$.
\item[$\cdot$]The interior of $\alpha$ is disjoint from $\boundary R$.
\end{enumerate}
then if $R'$ is the subsurface of $S$ obtained from $R$ by the handle-slide then the isotopy takes $\sigma(R';\Delta_i)$ to $\sigma(R;\ob{\Delta}_i)$.  
\end{enumerate}
\end{lemma}

\begin{proof}[Proof of Lemma \ref{isotopies of handle slides}]

\begin{figure}[ht] 
\scalebox{0.5}{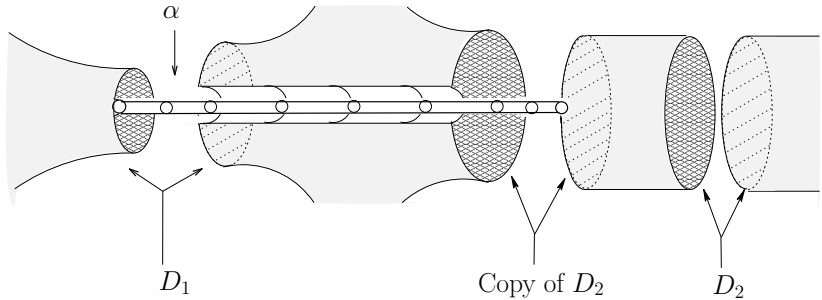}
\caption{Compressing along $\Delta_i$}
\label{Figure 1}
\end{figure}

Recall that $\Delta_i$ is obtained from $\ob{\Delta}_i$ by removing the disc $D_1$ and replacing it with the disc $D_1 \slide{\alpha} D_2$.  Let $S' = \sigma(S;\Delta_i)$.  When we compress along the discs $D_2$ and $D_1 \slide{\alpha} D_2$ we end up with a situation as depicted in Figure \ref{Figure 1}. Note that the figure depicts four discs parallel to $D_2$ since a disc parallel to $D_2$ makes up part of $D_1 \slide{\alpha} D_2$ and both $D_1 \slide{\alpha} D_2$ and $D_2$ are in $\Delta_i$. \newline

After compressing along $D_1 \slide{\alpha} D_2$ we see that there is regular neigborhood $N$ (in the compressionbody containing $\ob{\Delta}_i$) of $\alpha$ homeomorphic to $D^2 \times I$ with $D^2 \times \{0\}$ glued to a copy of $D_1$ and $D^2 \times \{1\}$ glued to a copy of $D_2 \times I$. Take a regular neighborhood in the compressionbody containing $\ob{\Delta}_i$ of $N \cup (D_2 \times I)$ which misses the rest of the surface $S'$. This regular neighborhood is a 3-ball $B$. Choose the regular neighborhoods so that $B \cap S' \subset \boundary B$. The intersection of $B$ with $D_1$ is a disc which is a regular neighborhood (in the compressionbody) of the point $\alpha \cap D_1$.  Slightly enlarge $B$ in $M$ to a ball $B'$ and perform an ambient isotopy supported on $B'$ and which takes $B - D_1$ to $B \cap D_1$.  Next use the regular neighborhood of $\alpha$ to isotope back to $S$ the portion of $S'$ which forms part of the boundary of a regular neighborhood of $\alpha$ (the ``trough").  The result is the same as if we had compressed along $\ob{\Delta}_i$.  This proves statement (a). The isotopy described is the identity off a neighborhood of $D_1 \cup \alpha \cup D_2$ and so is a proper isotopy.  \newline

If the handle-slide was relative to a closed set $X$, then by choosing the neighborhoods of $D_1$, $D_2$, and $\alpha$ to be disjoint from $X$, the isotopy described is relative to $X$.  This proves statement (b). \newline

To prove conclusion (c), we examine the possibilities.  Suppose that $R$ is a subsurface of $S$ as in the statement and suppose that $R'$ is obtained from $R$ by the handle-slide.  Recall that $B$ is the ball which is a regular neighborhood of $\alpha$ and $D_2$.  The important observation is that the isotopy takes $\boundary B - D_1$ into $D_1$.
\begin{itemize}

\item Suppose that $\boundary D_1 \subset \boundary R$ and that $\alpha \subset R$.  In this case, $R'$ equals $R - \boundary B$.  The isotopy fixes $R' - \eta(\boundary D_1 \slide{\alpha} D_2 \cup D_2)$. And so the isotopy takes $\sigma(R';\Delta_i)$ into $\sigma(R;\ob{\Delta}_i)$ \newline

\item Suppose that $\boundary D_1 \subset \boundary R$ and that $\alpha$ is not contained in $R$.  Then $\sigma(R';\Delta_i)$ equals $\sigma(R;\ob{\Delta}_i) \cup \boundary B$.  The isotopy described takes $\boundary B - D_1$ into $\sigma(R;\ob{\Delta}_i)$.  \newline

\item Suppose that $\boundary D_1$ is not contained in $R$.  The previous case shows that $\sigma(\cl(S - R');\Delta_i)$ is taken into $\sigma(\cl(S - R);\ob{\Delta}_i)$ and by part (a) we must have that $\sigma(R';\Delta_i)$ is taken into $\sigma(R;\ob{\Delta}_i)$.\newline

\item Suppose that $\boundary D_1 \subset \interior R$.  In this case, $\boundary B - D_1$ is contained in $\sigma(R';\Delta_i)$ and $(\boundary B - D_1) \cap S$ in $\sigma(R;\ob{\Delta}_i)$.  The isotopy clearly satisfies (c). 
\end{itemize}
\end{proof}

We now turn to the proof of Proposition \ref{isotopies of slide-moves}.

\begin{proof}[Proof of Proposition \ref{isotopies of slide-moves}]

Suppose that the 2-sided disc family $\Delta$ is obtained from the 2-sided disc family $\ob{\Delta}$ by a finite sequence $\{\mu_1,\hdots, \mu_{n}\}$ of slide-moves.  Each $\mu_i$ is a slide-move of type (M1), (M2), (M3), or (M4).  We prove the proposition by induction on the length of the sequence.  If the sequence is of length 0 the result is immediate so suppose that $n \geq 1$ and that the proposition is true for all sequences with $n - 1$ elements.  \newline

Let $\delta$ be the 2-sided disc family obtained from $\ob{\Delta}$ by the sequence $\nu = \{\mu_1, \hdots, \mu_{n-1}\}$.  Using the notation from the statement of the proposition: let $r$ be the subsurface of $S$ obtained from the subsurface $R$ by the sequence $\nu$. \newline

By the induction hypothesis, there is a collection of disjoint discs $\mc{E}$ with boundary on $\sigma(S;\ob{\Delta})$ and there is an ambient isotopy $f$, fixed off a compact set, which takes $\sigma(S;\delta)$ to $\sigma(S;\ob{\Delta} \cup \mc{E})$ and which takes the surface $\sigma(r;\delta)$ into the surface $\sigma(R;\ob{\Delta} \cup \mc{E})$. (Recall that this means $R$ compressed along those discs of $\ob{\Delta} \cup \mc{E}$ with boundary on $R$.)  We assume that $f$ also satisfies conclusions (i), (ii), and (iii).\newline

The 2-sided disc family $\Delta$ is obtained from the 2-sided disc family $\delta$ by a single slide-move $\mu_n$ of type (M1), (M2), (M3), or (M4).  We divide the proof into the case when $\mu_n$ is of type (M2) or (M4) and the case when the slide-move is of type (M1) or (M3).

\subsubsection*{Case: $\mu_n$ is of type (M2) or (M4)} If $\mu_n$ is a slide-move of type (M2) or (M4), $\Delta$ is obtained from $\delta$ by adding a single disc $D'$ to $\delta$.  In this case, $R' = r$.  The ambient isotopy $f$ takes the disc $D'$ to a disc $D$ with boundary on $\sigma(S;\ob{\Delta}\cup \mc{E})$.  By a further isotopy, if necessary, we may arrange that the disc $D$ has boundary disjoint from the remnants of $\mc{E}$ and so has boundary on $\sigma(S;\ob{\Delta})$ and that $D$ is disjoint from the discs of $\mc{E}$, though it may intersect $S$ in a neighborhood of $\mc{E}$.  Let $\mc{D}$ equal $\mc{E} \cup D$.  We need to show that we have satisfied the conclusions of the proposition.  \newline

To prove (i), recall that the discs $\mc{E}$ are numbered.  The disc $D$ should be given the next number.  Since $\interior D'$ is disjoint from $\sigma(S;\delta)$ and the isotopy is an ambient isotopy the disc $D$ has interior disjoint from $\sigma(S;\ob{\Delta} \cup \mc{E})$.  Thus, $D$ intersects $\sigma(S;\ob{\Delta} \cup \mc{E})$ only on $\boundary D$. \newline

Conclusion (ii) is clear, since the isotopy $f$ took $\sigma(S;\Delta' - D')$ to $\sigma(S;\ob{\Delta} \cup \mc{E})$ and also took $D'$ to $D$ which is a disc with boundary on the surface obtained from $\sigma(S;\ob{\Delta})$ by compressing along $\mc{E}$. \newline

To prove (iii), recall that since $\mu_n$ is the slide-move consisting of adding the disc $D'$ to $\delta$, the surface $R'$ equals the surface $r$.  The induction hypothesis says that $f$ takes $\sigma(r;\delta)$ to the surface $\sigma(R;\ob{\Delta} \cup \mc{E})$.  Conclusion (ii) shows that the isotopy $f$ takes the surface $\sigma(R';\Delta = \delta \cup D')$ to $\sigma(R;\ob{\Delta} \cup \mc{D})$. \newline

\subsubsection*{Case: $\mu_n$ is of type (M1) or (M3)}
If the slide-move $\mu_n$ is of type (M1) or (M3) we have obtained $\Delta$ from $\delta$ by a single handle-slide of $\delta_1$ in $U$ or $\delta_2$ in $V$.  Without loss of generality, assume that $\mu_n$ is a slide-move of type (M3), so that $\Delta$ is obtained from $\delta$ by the slide-move $\mu_n$ of $\delta_2$.   By Lemma \ref{isotopies of handle slides}, there is an ambient isotopy $g$ of $M$, fixed off a compact set, which satisfies properties (a), (b), and (c).  In particular, $g$ takes $\sigma(S;\Delta)$ to $\sigma(S;\delta)$ because it takes $\sigma(S;\Delta_2)$ to $\sigma(S;\delta_2)$ and is performed relative to $\boundary \delta_1$ (property (b)).  Let $h$ be the ambient isotopy formed by performing $g$ and then performing $f$.  Let $\mc{D} = \mc{E}$.  We show that $h$ satisfies conclusions (i), (ii), and (iii).  \newline

Conclusions (i) and (ii) follow immediately from the induction hypothesis on $f$ and property (a) of Lemma \ref{isotopies of handle slides}. \newline

To prove conclusion (iii), recall that $r$ denotes the surface obtained from $R$ by the sequence of slide-moves $\{\mu_1, \hdots, \mu_{n-1}\}$.  The surface $R'$ is obtained from $r$ by the handle-slide $\mu_n$.  Property (c) from Lemma \ref{isotopies of handle slides} shows that $g$ takes $\sigma(R';\Delta_2)$ to $\sigma(r;\delta_2)$.  The isotopy $g$ is an ambient isotopy which was performed relative to $\boundary \delta_1$, so $g$ also takes $\sigma(R';\Delta)$ to $\sigma(r;\delta)$.  By induction, the isotopy $f$ takes $\sigma(r;\delta)$ to $\sigma(R;\ob{\Delta}\cup \mc{D})$. And so $h$ satisfies (iii).
\end{proof}

   
\section{Relative Heegaard Splittings}\label{rel HS}

\subsection{The Outer Collar Property}\label{outer collar prop}
Recall that a compact relative compressionbody $H$ is formed by taking a compact surface $F$ (possibly $\nil$), forming $F \times I$, and attaching 1-handles to the interior of $F \times \{1\}$ and finitely many 3-balls.  The surface $F \times \{0\}$ is $\boundary_- H$ and the closure of its complement in $\boundary H$ is the preferred surface of $H$, denoted $\boundary_+ H$.  Recall, also, from Section \ref{handle-slides} that a relative compressionbody $A$ is correctly embedded in a compressionbody $B$ if the frontier of $A$ in $B$ consists only of components of $\boundary_- A$ which have boundary.  Corollary \ref{complementary compressionbodies} states that, in this case, $\cl(B - A)$ is a relative compressionbody with some preferred surface.  However, $\cl(B - A)$ may not be correctly embedded in $B$ as we may not be able to choose $\cl(\boundary_+ B - \boundary_+ A)$ to be the preferred surface of $\cl(B - A)$. \newline

In this section we explore situations in which we can ``come close" to having $\boundary_+ \cl(B - A)$ equal $\cl(\boundary_+ B - \boundary_+ A)$.  These situations will arise when we have exhausting sequences of noncompact absolute compressionbodies.

\begin{definition}
Suppose that $\{K'_i\}$ is an exhausting sequence for a noncompact absolute compressionbody $U$.  If each $K'_i$ is a relative compressionbody correctly embedded in $U$ and each $K'_i$ is correctly embedded in $K'_{i+1}$ then $\{K'_i\}$ is a \defn{correctly embedded exhausting sequence} for $U$.  
\end{definition}

The following definition is somewhat technical, but will be useful for statements of results in Section \ref{Exh. Seq.}.  Recall that a collaring set of discs for a compressionbody $H$ is a set of discs which separates off a copy of $\boundary_- H \times I$.

\begin{definition}
Suppose that $\{K'_i\}$ is a correctly embedded exhausting sequence for $U$. Suppose that for each $i \geq 2$ there is an embedding of $(\fr K'_i \times I, (\boundary \fr K'_i) \times I)$ into $(\cl(K'_i - K'_{i-1}),\boundary_+ U \cap \cl(K'_i - K'_{i-1}))$ so that $\fr K'_i = \fr K'_i \times \{0\}$ and so that $\fr K'_i \times \{1\}$ is a subsurface of $\boundary_+ U$ except at a finite number of open discs.  Then $\{K'_i\}$ is said to have the \defn{outer collar property}.
\end{definition}

\begin{remark}
The open discs of $\fr K'_i \times \{1\}$ which are not contained in $\boundary_+ U$ are the interiors of a set of collaring discs for $K'_i$.
\end{remark}

\begin{definition}
Suppose that $\{K'_i\}$ is a correctly embedded exhausting sequence for $U$.  Suppose that for each $i \geq 2$ there is an embedding of $(\fr K'_{i-1} \times I,(\boundary \fr K'_{i-1}) \times I)$ into $(\cl(K'_i - K'_{i-1}),\boundary_+ U \cap \cl(K'_i - K'_{i-1}))$ so that $\fr K'_{i-1} = \fr K'_{i-1} \times \{0\}$ and so that $\fr K'_{i-1} \times \{1\}$ is a subsurface of $\boundary_+ U$ except at a finite number of open discs.  Then $\{K'_i\}$ is said to have the \defn{inner collar property}.
\end{definition}

The outer collar property plays an important role in this paper, so it may be helpful to give an example of an exhausting sequence of a handlebody with the outer collar property.  Our example is, in fact, an exhausting sequence of a one-ended, infinite genus handlebody which has both the inner and outer collar properties.

\begin{figure}[ht]
\scalebox{1.0}{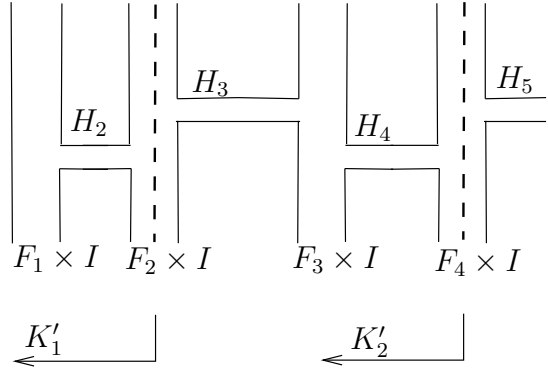}
\caption{The handlebody $H$.}
\label{innerouter}
\end{figure}

\begin{example}
For each natural number $i$, let $F_i$ be a compact, connected surface with non-empty boundary.  Let $P_i = F_i \times I$.  Recall that $P_i$ is a handlebody.  For each $i \geq 2$ join $F_i \times \{0\}$ to $F_{i-1} \times \{1\}$ by a one-handle $H_i$.  Denote the union of all the product regions and all the one-handles by $H$.  See Figure \ref{innerouter} for a schematic depiction of $H$.  Let $D_i$ be a disc which is a cocore of the one-handle $H_i$. Let 

$$K'_i = P_1 \cup H_2 \cup P_2 \cup \hdots \cup P_{2i - 1} \cup H_{2i} \cup (F_{2i} \times [0,\frac{1}{2}]).$$  

The construction makes clear that $\{K'_i\}$ is a correctly embedded exhausting sequence of the handlebody $H$.  The frontier of $K'_i$ is $F_{2i} \times \{\frac{1}{2}\}$ which is an incompressible surface in $H$.  For $i \geq 2$, compressing $K'_i$ along $D_{2i}$ leaves two components, one of which is $F_{2i} \times [0,\frac{1}{2}] = \fr K'_i \times I$.  This component is disjoint from $K'_{i - 1}$.  From the construction, it is clear that $\{K'_i\}$ has the outer collar property.  For $i \geq 2$, boundary-reducing $\cl(K'_i - K'_{i-1})$ along the disc $D_{2i - 1}$ leaves two components, one of which is $F_{2i - 2} \times [\frac{1}{2},1] = \fr K'_{i-1} \times I$.  Again, from the construction it is clear that $\{K'_i\}$ has the inner collar property.         
\end{example}

In this paper, it is the outer collar property which is most used.  The inner collar property makes an appearance in the appendix.  Certainly not every correctly embedded sequence has the outer collar property.  If, for example, for some $i$, $\boundary_- K'_{i-1}$ was not a disc and bounded a product region with $\boundary_- K'_i$, the sequence would not have the outer collar property.   If a sequence has both the inner and outer collar properties we can take $\cl(\boundary_+ K'_{i+1} - \boundary_+ K'_i)$ to be the preferred surface for the relative compressionbody $\cl(K'_{i+1} - K'_i)$.  It is in this sense that a sequence with the outer collar property ``comes close" to having $\cl(\boundary_+ K'_{i+1} - \boundary_+ K'_i)$ the preferred surface for the relative compressionbody $\cl(K'_{i+1} - K'_i)$. \newline

It turns out that sequences with the outer collar property are fairly common:

\begin{lemma}\label{outer collar property}
Suppose that $\{K'_i\}$ is a correctly embedded exhausting sequence for the absolute compressionbody $U$.  Then there is a subsequence with the outer collar property.
\end{lemma}

\begin{proof}
Let $\{L_i\}$ be an exhausting sequence for $U$ by subcompressionbodies.  Recall that the frontier of a subcompressionbody consists of properly embedded discs.  Take subsequences of $\{K'_i\}$ and $\{L_i\}$ so that for all $i$, $L_i \subset K'_i \subset L_{i+1}$.  Each inclusion should be into the interior of the succeeding submanifold. Fix some $i \in \N$.\newline

By Lemma \ref{nested reducing discs}, we may choose a defining collection of discs $\Delta$ for $K'_i$ which includes the discs $\fr L_i$.  Boundary-reducing $K'_i$ along $\Delta$ leaves us with 3-balls and products $\boundary_- K'_i \times I$.  Since $K'_{i-1} \subset L_i$ and $\fr L_i$ separates $U$ we have that the remnants of $K'_{i-1}$ are completely contained in the 3-balls.  Thus the product regions $\boundary_- K'_i \times I$ are contained completely in $\cl(K'_i - K'_{i-1})$. Label $\boundary_- K'_{i}$ with $\boundary_- K'_i \times \{1\}$.  The discs of $\Delta$ which show up on $\boundary_- K'_i \times \{0\}$ can be taken to be our collaring set of discs.  This collaring set of discs and the product region $\boundary_- K'_i \times I$ are contained in $K'_i - K'_{i-1}$ so the sequence $\{K'_i\}$ now has the outer collar property.  
\end{proof}

\subsection{Relative Heegaard Splittings}Suppose that $M = U \cup_S V$ is an absolute Heegaard splitting of a non-compact 3-manifold with compact boundary, containing no $S^2$ components.  If $N \subset M$ is a compact submanifold, the surface $S \cap N$ cannot possibly give an absolute Heegaard splitting of $N$ as $S$ is non-compact and $N$ is compact.  It can, however, give a relative Heegaard splitting of $N$.\newline  

We will eventually look at the relationship between relative Heegaard splittings and absolute Heegaard splittings, but first we show how exhaustions of $M$ by compact submanifolds which inherit relative Heegaard splittings from $S$ give rise to correctly embedded exhausting sequences of $U$.

\begin{definition}
A submanifold $N$ contained in $M$ is \defn{adapted} to $S$ if \linebreak[4]$(U \cap N) \cup_{S \cap N} (V \cap N)$ is a relative Heegaard splitting of $N$ and $(U \cap N)$ is correctly embedded in $U$ and $(V \cap N)$ is correctly embedded in $V$.  An exhausting sequence $\{K_i\}$ is \defn{adapted} to $S$ if each $K_i$ is adapted to $S$.  It is \defn{perfectly adapted} to $S$ if it is adapted to $S$ and, additionally, each $\cl(K_{i+1} - K_i)$ is adapted to $S$.
\end{definition}

\begin{remark}
The requirement that $(U \cap N)$ and $(V \cap N)$ are correctly embedded in $U$ and $V$ respectively means that $\fr N$ can have no closed components which are contained entirely in $U$ or $V$: such a component would have to be a component of $\boundary_- (U \cap N)$ or $\boundary_- (V \cap N)$ as it would not be a subsurface of $S$.  This, however, means that $U \cap N$ or $V \cap N$ is not correctly embedded in $U$. 
\end{remark}

In constructing an exhausting sequence that is adapted to $S$, the requirement that $U \cap N$ and $V \cap N$ are correctly embedded in $U$ and $V$ is a minor one.  To see this, suppose that a compact submanifold $N \subset M$ containing $\boundary M$ has the property that $N \cap U$ and $N \cap V$ are relative compressionbodies with preferred surfaces $S \cap N$.  It is easy to adjust $N$ so that $U \cap N$ and $V \cap N$ are correctly embedded.  If $U \cap N$, say, is not correctly embedded there must be a component $F$ of $\boundary_- (U \cap N) - \boundary M$ which is a closed surface. Since $U$ is an absolute compressionbody, $\operatorname{H}_2(U,\boundary_- U) = 0$.  Thus either $F$, or $F$ and components of $\boundary_- U \cap \boundary M$, bound(s) a compact submanifold $L$ of $U$.  $L$ cannot be interior to $N$ as $N \cap U$ is a relative compressionbody with non-empty preferred surface and $F \cup \boundary_- U$ is contained in $\boundary_- (U \cap N)$.  Thus $L$ is exterior to $N$.  Since $\boundary_- U \subset \boundary M \subset N$, we have that $\boundary L = F$.  In fact, $(N \cap U) \cup L$ must still be a relative compressionbody.  To see this, note that $F$ must be compressible in $L$ as $F$ is incompressible in $N \cap U$. ($\boundary_-$ of a compressionbody is incompressible in the compressionbody).  Every closed incompressible surface in $U$ is parallel to $\boundary_- U$.  Boundary-reducing $L$ is the same as adding a 2-handle to $N \cap U$ along a curve in $F \subset \boundary_- (N \cup L)$.  Adding a 2-handle to $\boundary_-$ of a compressionbody preserves the fact that we have a relative compressionbody (up to the introduction of spherical boundary components).  Eventually, our surface is a collection of spheres, which, since $U$ is irreducible, bound balls in $U$.  Including these balls into $N$ (with the 2-handles attached) also preserves the fact that we have a relative compressionbody.  

\begin{lemma}\label{adapted means correctly embedded}
If $\{K_i\}$ is an exhausting sequence of $M$ adapted to $S$ with $\boundary M \subset K_1$ then $\{K_i \cap U\}$ and $\{K_i \cap V\}$ are correctly embedded exhausting sequences of $U$ and $V$ respectively.
\end{lemma}

\begin{proof}
Let $X$ denote either $U$ or $V$.  Since $\{K_i\}$ is adapted to $S$, by definition each $K_i \cap X$ is correctly embedded in $X$.  Thus, $\boundary_+ (K_i \cap X) \subset \boundary_+ (K_{i+1} \cap X)$.  Furthermore, any closed component of $\boundary_- (K_i \cap X)$ is a component of $\boundary_- X$ which is contained in $\boundary_- (K_{i+1} \cap X)$.  Thus, each closed component of $\boundary_- (K_i \cap X)$ is a component of $\boundary_- (K_{i+1} \cap X)$.  Hence, $K_i \cap X$ is correctly embedded in $K_{i+1} \cap X$.
\end{proof}

\begin{definition}
If $\{K_i\}$ is an exhausting sequence of $M$ adapted to $S$ with $\boundary M \subset K_1$ and such that $\{K_i \cap U\}$ has the outer collar property we say that $\{K_i\}$ has \defn{the outer collar property with respect to $U$}.
\end{definition}

\begin{corollary}\label{outer collar property 2}
If $\{K_i\}$ is an exhausting sequence of $M$ adapted to $S$ with $\boundary M \subset K_1$ then there is a subsequence which has the outer collar property with respect to $U$.
\end{corollary}

\begin{proof}
By Lemma \ref{adapted means correctly embedded}, $\{K_i \cap U\}$ is a correctly embedded exhausting sequence of $U$.  By Lemma \ref{outer collar property}, there is an infinite subset $\mc{N}$ of $\N$ such that $\{K_i \cap U\}_{i \in \mc{N}}$ has the outer collar property.  Hence, $\{K_i\}_{i \in \mc{N}}$ has the outer collar property with respect to $U$.
\end{proof}
\subsection{Balanced Exhausting Sequences}\label{Balanced Exhausting Sequences}
We've shown so far that if $M$ has an exhausting sequence adapted to $S$ we can find one which has the outer collar property.  We've not yet addressed the question of the existence of an exhausting sequence adapted to $S$.  We do that now.  This construction is a variation of the construction given by Frohman and Meeks in \cite{FrMe97}. \newline

Recall that $M = U \cup_S V$ is an absolute Heegaard splitting of a non-compact 3-manifold with compact boundary.  Let $A$ and $B$ be compact subcompressionbodies of $U$ and $V$ respectively with the property that $\boundary_S A \subset \interior(\boundary_S B)$.  Let $C$ be a regular neighborhood of $A \cup B$.

\begin{definition}
A set $C$ constructed in such a manner will be called a \defn{balanced submanifold} of $M$ (with respect to $S$).  An exhausting sequence $\{C_i\}$ of $M$ will be called a \defn{balanced exhausting sequence} for $M$ (with respect to $S$) if each $C_i = \eta(A_i \cup B_i)$ is a balanced submanifold and if, for all $i$, $\boundary_S B_i \subset \interior(\boundary_S A_{i+1})$.
\end{definition}

The next lemma guarantees that balanced submanifolds are adapted to the Heegaard surface.  Consequently, we will say that such a set $C$ is a \defn{balanced submanifold of $M$ (adapted to $S$)}.

\begin{lemma}[{Frohman-Meeks,\cite[Proposition 2.2]{FrMe97}}]\label{balanced are adapted}
If $C$ is a balanced submanifold of $M$ with respect to $S$ then $C$ is adapted to $S$.  
\end{lemma}

\begin{proof}
We must show that $(U \cap C) \cup_{S \cap C} (V \cap C)$ is a relative Heegaard splitting of $C$. In other words, we must show that $U \cap C$ and $V \cap C$ are both relative compressionbodies with preferred surface $S \cap C$. \newline

Assume that $C$ is a regular neighborhood of $A \cup B$ where $A$ and $B$ are compact subcompressionbodies of $U$ and $V$ respectively and $\boundary_S A \subset \boundary_S B$.  We have $C \cap V = \eta(B)$ so $C \cap V$ is a relative compressionbody with preferred surface $S \cap C$.  To obtain $C \cap U$ we take a regular neighborhood of $\cl(\boundary_S B - \boundary_S A)$ in $U$ and glue it to $A$.  An alternative way of performing the construction is as follows. \newline

Let $D$ be the collection of discs which make up the frontier of $A$.  Take a regular neighborhood of $D$ and let $D'$ be the components of the frontier of the neighborhood which are not in $A$.  Let $F = \cl(\boundary_S B - (\boundary_S A \cup \eta(D)))$.  Take a regular neighborhood of $F$ in $U - A$.  Consider $F$ to be $F \times \{1\}$.  Since $D' \subset F$, this regular neighborhood contains $D' \times I$.  See Figure \ref{compressionbody}.  This revised neighborhood is $\boundary_- C \times I$.  We may then add one-handles so that one end of each one-handle is on a disc of $D' \times \{1\}$ and the other end is on the corresponding disc of $D$.  It is clear that $S \cap C$ is the preferred surface of this compressionbody.
\begin{figure}[ht]
\scalebox{0.5}{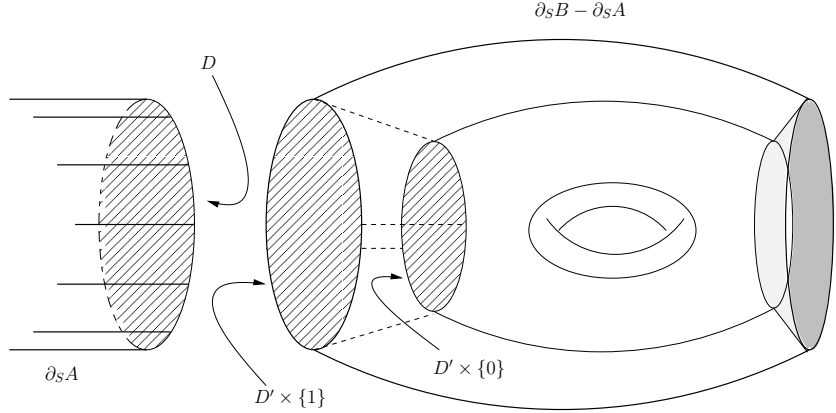}
\caption{Adding a regular neighborhood of $\boundary_S B - \boundary_S A$ to $A$ gives us a relative compressionbody.}
\label{compressionbody}
\end{figure}
\end{proof}

To obtain a balanced exhausting sequence of $M$, start by taking exhausting sequences $\{A_i\}$ and $\{B_i\}$ of $U$ and $V$ by subcompressionbodies.  Since each $A_i$ and each $B_i$ are compact we may take subsequences of $\{A_i\}$ and $\{B_i\}$ so that, for all $i$, $\boundary_S A_i \subset \boundary_S B_i \subset \boundary_S A_{i+1}$.  Each of the inclusions should be into the interior of the succeeding surface. \newline

A component of the frontier of a balanced submanifold $C$ can be thought of as being a compact subsurface of $S$ with discs, each contained entirely in $U$ or $V$, glued onto the boundary components.  In fact, since each component of the frontier of each balanced submanifold intersects $S$, neither $\boundary_-(C \cap U)$ nor $\boundary_- (C \cap V)$ have components which are closed surfaces not contained in $\boundary M$.  Thus, if we have a balanced exhausting sequence $\{C_i\}$ of $M$ adapted to $S$ with $\boundary M \subset C_1$, it is adapted to $S$ in the sense of the definition given at the beginning of this section.  By Corollary \ref{outer collar property 2}, we can take a subsequence of $\{C_i\}$ so that it has the outer collar property. 

\begin{remark}
Even though we have a balanced exhausting sequence $\{C_i\}$ of $M$ which is adapted to $S$ and has the outer collar property, there is no reason to suppose that it is perfectly adapted to $S$.  The difficulty is in the fact that $\cl(C_{i+1} - C_i) \cap U$ may not be a relative compressionbody with preferred surface $S \cap \cl(C_{i+1} - C_i)$.
\end{remark}

Let $\{C_i\}$ be a balanced exhausting sequence for $M$ adapted to $S$.  Each $C_i$ is the neighborhood of $A_i \cup B_i$ where $A_i$ and $B_i$ are subcompressionbodies of $U$ and $V$ respectively.  As such, the collection of discs $\ob{\Delta} = \cup_i(\fr A_i \cup \fr B_i)$ is a 2-sided disc family for $S$ in $M$.  (The notation means the frontier of $A_i$ in $U$ and the frontier of $B_i$ in $V$.)  We can perform a finite sequence of slide-moves (Section \ref{slide-moves and isotopies}) on $\ob{\Delta}$ to obtain a new 2-sided disc family $\Delta$.  This sequence also gives us, for each $i$, subcompressionbodies $A'_i$ and $B'_i$ obtained from $A_i$ and $B_i$ respectively by slide-moves.  The important observation is that since the slide-moves are done relative to $\cup_i(\boundary \fr A_i \cup \boundary \fr B_i)$ we still have, for each $i$, that $\boundary_S A'_i \subset \boundary_S B'_i \subset \boundary_S A'_{i+1}$.  Thus, $C'_i = \eta(A'_i \cup B'_i)$ is a balanced submanifold of $M$ adapted to $S$.  And so, $\{C'_i\}$ is a balanced exhausting sequence of $M$ adapted to $S$.  These observations provide the key to the proof of Theorem \ref{well-placed}.

\begin{definition}
The balanced submanifold $C' = \eta(A' \cup B')$ is obtained from the balanced submanifold $C = \eta(A \cup B)$ \defn{by slide-moves} if there is a finite sequence of slide-moves by which $A'$ is obtained from $A$ and $B'$ is obtained from $B$.
\end{definition}

\subsection{Comparing Absolute and Relative Heegaard Splittings}
In the remainder of this section, we look at the relationship between absolute and relative Heegaard splittings of a compact manifold.  These results will help us to translate facts about absolute Heegaard splittings to relative Heegaard splittings.  Let $N$ denote a 3-manifold, compact or non-compact, with non-empty compact boundary. \newline

Suppose that $N = U \cup_S V$ is a relative Heegaard splitting.  Let $\mc{B}$ denote the boundary components of $N$ which intersect $S$. Define $\hat{U}$ to be $U$ together with a regular neighborhood of $\mc{B}$.  Define $\hat{V}$ to be the closure of the complement of $\hat{U}$ in $N$ and let $\hat{S} = \hat{U} \cap \hat{V}$.

\begin{lemma}
$N = \hat{U} \cup_{\hat{S}} \hat{V}$ is an absolute Heegaard splitting of $N$.
\end{lemma}
\begin{proof}
If $\mc{B} = \nil$ there is nothing to prove, so assume that $\mc{B}$ is non-empty.  $U$ is a relative compressionbody and so is obtained from $F \times I$ by adding 1-handles to $F \times \{1\}$ and countably many 3-balls.  $F$ is a compact surface with boundary.  Let $B$ be a component of $\mc{B}$ and let $B_U = B \cap U$ and $B_V = B \cap V$.  In the process of obtaining $\hat{U}$ we glue $B_V \times I$ to $B_U \times I$ along $\gamma \times I$ where $\gamma = \boundary B_V = \boundary B_U$.  So $\hat{U}$ is $B \times I$ attached by 1-handles to the preferred surface of a compressionbody.  Hence, performing this operation for each boundary component of $N$ which intersects $S$, leaves us with $\hat{U}$, an absolute compressionbody. On the other hand, to form $\hat{V}$ we have removed a collar neighborhood of each component of $\boundary_- V$ which intersected $\boundary_+ V$.  Let $\mc{D}$ be a collaring set of discs for $V$.  The discs $\mc{D}$ are also discs in $\hat{V}$.  Let $\mc{E}$ be the collection of components of $\sigma(V;\mc{D})$ which contain $\boundary_- V \cap \mc{B}$.  Each of these components is a $\text{(surface with boundary)} \times I$.  As such, each component is a handlebody.  Removing a collar neighborhood of $\boundary_- V \cap \mc{B}$ from these components does not change the homeomorphism type.  The space $\hat{V}$ is formed by attaching these handlebodies to the preferred surface of the absolute compressionbody $\sigma(V;\mc{D}) - \mc{E}$ by 1-handles dual to the discs $\mc{D}$.  Thus, $\hat{V}$ is an absolute compressionbody. 
\end{proof}

If we know that $V$ intersects $\boundary N$ in discs, the relationship is stronger.

\begin{lemma}[The Marionette Lemma]\label{equivalence of absolute and relative HS}
Suppose that $U_S \cup_S V_S$ and $U_T \cup_T V_T$ are two relative Heegaard splittings of a 3-manifold $N$.  Suppose also that for each component of $\boundary N$ which intersects $S$, $V_S$ and $V_T$ intersect that component in discs.  If for each such boundary component of $N$, $V_S$ and $V_T$ intersect that boundary component in the same number of discs, then $S$ and $T$ are properly ambient isotopic if and only if $\hat{S}$ and $\hat{T}$ are properly ambient isotopic.
\end{lemma}

We form $\hat{U_S}$ and $\hat{U_T}$ by including a regular neighborhood of $V_S \cap \boundary N$ and $V_T \cap \boundary N$ into $U_S$ and $U_T$.  If we want to undo this operation we can remember the cocores of the discs $V_S \cap \boundary N$ and $V_T \cap \boundary N$.  These give us finite collections of arcs in $\hat{U_S}$ and $\hat{U_T}$ joining $\boundary N$ to $\hat{S}$ and $\hat{T}$ respectively.  To prove the lemma, we need to understand how these arcs can be isotoped within the compressionbodies $\hat{U_S}$ and $\hat{U_T}$.  We will show that if $\hat{S}$ and $\hat{T}$ are isotopic, we can isotope $\hat{S}$ and $\hat{T}$ to coincide and then isotope the arcs to coincide.

\begin{definition}
Let $\psi$ be a finite collection of arcs in an absolute compressionbody $H$ with at least one endpoint of each arc on $\boundary_+ H$.  If $H$ is a 3-ball then $\psi$ is \defn{standard} if it is isotopic to a collection of arcs which lie in $\boundary_+ H = \boundary H$.  If $H = F \times I$ where $F$ is a closed connected surface, then $\psi$ is \defn{standard} if there is an isotopy of $\psi$ so that each spanning arc is vertical in the product structure and each non-spanning arc is contained in $F \times \{1\} = \boundary_+ H$.  For a generic absolute compressionbody, $\psi$ is \defn{standard} if there is a defining collection of discs $\Delta$ for $H$ which is disjoint from $\psi$ and such that $\psi$ is standard in each component of $\sigma(H;\Delta)$.
\end{definition}

We need the following two results which are slightly rephrased from \cite{ScTh93}.  We are allowing our compressionbody to be non-compact, but since the number of arcs is finite the results are still true, as we may restrict attention to a compact subcompressionbody.

\begin{lemma}[{Scharlemann-Thompson, \cite[Lemma 6.4]{ScTh93}}]\label{standard collections}
If $\sigma$ and $\tau$ are standard collections of arcs in an absolute compressionbody $H$, then for any defining collection of discs $\Delta$ for $H$ there is an isotopy of $\sigma$ and an isotopy of $\tau$ so that $\sigma$ and $\tau$ are standard in $\sigma(H;\Delta)$.
\end{lemma}

\begin{lemma}[{Scharlemann-Thompson, \cite[Corollary 6.7]{ScTh93}}]\label{inducing standard collections}
Let $\psi$ be a collection of arcs properly embedded in a compressionbody $H$ such that for every subcollection $\psi' \subset \psi$, the complement of $\psi'$ is a compressionbody.  Then $\psi$ is standard.
\end{lemma}

\begin{proof}[Proof of the Marionette Lemma]
If $S$ and $T$ are ambient isotopic, it is clear that $\hat{S}$ and $\hat{T}$ are.  So suppose that $\hat{S}$ and $\hat{T}$ are ambient isotopic. \newline

As mentioned earlier, we can recover $S$ and $T$ from $\hat{S}$ and $\hat{T}$ by remembering the cocores of the 2-handles that were added to $U_S$ and $U_T$.  Let $\sigma$ be the collection of arcs coming from $V_S \cap \boundary N$ and let $\tau$ be the collection of arcs coming from $V_T \cap \boundary N$.  \newline

Isotope $\hat{S}$ onto $\hat{T}$.  Now we have $\hat{U_S} = \hat{U_T}$.  This isotopy takes $\sigma$ to some collection of arcs which we continue to call $\sigma$.  
If we can show that there is an isotopy of $\sigma$ onto $\tau$ which keeps $\hat{S}$ mapped onto $\hat{T}$ for all time, we will be done.  The isotopy is allowed to move the endpoints of the arcs, but it must keep them on $\boundary N \cup \hat{S}$. \newline

We claim, first, that for each subcollection $\sigma'$ of arcs in $\sigma$ the complement of $\sigma'$ in $\hat{U_S} = \hat{U_T}$ is a compressionbody.  Let $\sigma'$ be a subcollection of arcs from $\sigma$.  Let $s'$ denote the arcs of $\sigma - \sigma'$.  Let $D_{s'}$ be the 2-handles of $\eta(V_T \cap \boundary N)$ which have cocores $s'$.  Consider the relative compressionbody $U_S$.  $U_S$ is formed by taking a surface $F$ with boundary, forming $F \times I$ and adding 1-handles to $F \times \{1\}$.  The surface $F$ has one boundary component for each component of $S \cap \boundary N$.  Let $\gamma$ denote the boundary components of $F \times \{0\}$ which correspond to $s'$.  Adding the 2-handles $D_{\sigma'}$ to $U_S$ is achieved by attaching copies of $D^2 \times I$ to $F$ along $\gamma \times I$.  It's clear that the result is still a compressionbody.  But this is exactly $\cl(\hat{U_S} - \eta(\sigma'))$.  Thus, the complement of every subcollection of $\sigma$ in $\hat{U_S}$ is a compressionbody.  The same result holds for $\tau$. \newline

By Lemma \ref{inducing standard collections}, both $\sigma$ and $\tau$ are standard.  By Lemma \ref{standard collections}, there is a proper isotopy of $\sigma$ and a proper isotopy of $\tau$ so that both $\sigma$ and $\tau$ are disjoint from a defining disc set $\Delta$ for $\hat{U_S} = \hat{U_T}$ and both are standard in $\sigma(U_S;\Delta)$.  Since each arc of $\sigma \cup \tau$ has an endpoint on a component of $\boundary N$, we may assume that the isotopy has made each arc of $\sigma$ and each arc of $\tau$ vertical in the product structure of $(\boundary_N \times I) \cap U_S$.  Since for each component of $\boundary N$ the arcs of $\sigma$ and $\tau$ with an endpoint on that component are in one-to-one correspondence, there is the required isotopy taking $\sigma$ onto $\tau$.  
\end{proof}

The following is a version of Haken's Lemma for relative Heegaard splittings.  It is, perhaps, well-known.  It appears in similar versions as Lemma 5.2 in \cite{BaScSe05} and as a remark following Definition 2.1 in \cite{FrMe97}.

\begin{lemma}[Haken's Lemma]\label{Haken}
Suppose that $U \cup_S V$ is a relative Heegaard splitting of $N$ with the property that each component of $V \cap \boundary N$ is a disc.  Then if $\boundary N$ is compressible in $N$ there is a compressing disc for $\boundary N$ whose intersection with $S$ is a single simple closed curve.  Furthermore, boundary reducing $N$ along this disc leaves us with a relative Heegaard splitting $\cl(U - \eta(D)) \cup_{\cl(S - \eta(D))} \cl(V - \eta(D))$ of the resulting manifold.
\end{lemma}

\begin{proof}
Let $\hat{U} \cup_{\hat S} \hat{V}$ be the absolute Heegaard splitting for $N$ obtained by including $\eta(V \cap B)$ into $U$ for each component $B \subset \boundary N$ which intersects $S$.  Since $\boundary N$ is compressible, by Casson and Gordon's version of Haken's Lemma \cite{CaGo87}, there is a compressing disc $D$ for $\boundary N$ which intersects $\hat{S}$ in a single simple closed curve.  \newline

To obtain $U \cup_S V$ from $\hat{U} \cup_{\hat{S}} \hat{V}$ we include into $\hat{V}$ the neighborhood of a certain collection of arcs $\sigma$.  The arcs $\sigma$ are the cocores of the 2-handles which we added to $U$ in order to obtain $\hat{U}$.  \newline

If $\boundary D$ is on a component of $\boundary N$ contained in $\hat{V}$, then by Lemma \ref{standard collections} we may isotope $\sigma$ to be disjoint from the disc $D \cap \hat{U}$.  Thus, there is a compressing disc for $\boundary N$ which intersects $S$ in a single simple closed curve. \newline

If $\boundary D$ is on a component of $\boundary N$ contained in $\hat{U}$ then $D \cap \hat{U}$ is an annulus.  By performing handle-slides, we may obtain a defining collection of discs $\Delta$ for $\hat{U}$ which are disjoint from that annulus.  We may assume that the annulus $D \cap \hat{U}$ is vertical in the product structure of the component of $\sigma(\hat{U};\Delta)$ containing it.  By Lemma \ref{standard collections}, there is an isotopy of the arcs $\sigma$ so that $\sigma$ is disjoint from $\Delta$ and is vertical in the product structure of the components of $\sigma(\hat{U};\Delta)$ containing it.  It is then easy to isotope $\sigma$ to be disjoint from the annulus $D \cap \hat{U}$.  Hence, when we remove an open regular neighborhood of $\sigma$ from $\hat{U}$ to obtain $U$ we have the disc $D$ intersecting $S$ in a single simple closed curve.  Thus $S$ divides $D$ into a disc and an annulus. \newline

Boundary-reducing $N$ along $D$ leaves us with a 3-manifold $\ob{N} = \sigma(N;D)$.  We have boundary-reduced the relative compressionbody ($U$ or $V$) containing the disc part of $D$ along a disc with boundary in the preferred surface.  Thus, by Lemma \ref{boundary reducing gives compressionbody} it is still a relative compressionbody.  In the other compressionbody $X$ (equal to $V$ or $U$), there is a defining set of discs $\Delta$ disjoint from $D$ and the annulus $D \cap X$ is vertical in the product structure of the component of $\sigma(X;\Delta)$ containing it.  That component is homeomorphic to $F \times I$ where $F$ is a compact surface, possibly with boundary.  Removing the open neighborhood of a vertical annulus in such a component leaves us with a manifold homeomorphic to $G \times I$ where $G$ is a compact surface obtained from $F$ by removing an open annulus.  Thus, $X - \interior(\eta(D \cap X))$ is still a relative compressionbody with preferred surface $S - \interior(\eta(D))$. This implies that $\ob{N} = \cl(U - \eta(D)) \cup_{\cl(S - \eta(D))} \cl(V - \eta(D))$ is a relative Heegaard splitting.
\end{proof}

\section{Heegaard Splittings of Eventually End-Irreducible 3-manifolds}\label{Exh. Seq.}

\subsection{Introduction}\label{Types of Exh. Seq.}
Recall that a non-compact 3-manifold $M$ is \defn{end-irreducible rel $C$} for a compact set $C \subset M$ if there is an exhausting sequence $\{K_i\}_{\N}$ for $M$ such that $C \subset K_1$ and, for all $i$, $\fr K_i$ is incompressible in $M - C$.  Inessential spheres count as incompressible surfaces, so, for example, $\R^3$ is end-irreducible rel $\nil$. Other examples of eventually end-irreducible 3-manifolds are deleted boundary 3-manifolds.  A deleted boundary 3-manifold $M$ contains a compact set $C$ so that $\cl(M - C)$ is homeomorphic to $F \times \R_+$ for some closed surface $F$. \newline

For the remainder of this section, assume that $M$ is an orientable non-compact 3-manifold which is end-irreducible rel $C$ and that $\boundary M \subset C$. Let $M = U \cup_S V$ be an absolute Heegaard splitting of $M$.  \newline

Since we will be dealing with a variety of exhausting sequences for $M$ we collect the following definitions here:

\begin{definition}
Let $\{K_i\}$ be an exhausting sequence for $M$ with $C \subset K_1$.  We say that:

\begin{itemize}
\item $\{K_i\}$ is \defn{frontier-incompressible rel $C$} if, for each $i$, $\fr K_i$ is incompressible in $M - C$. \newline

\item $\{K_i\}$ is \defn{adapted} to $S$ if, for all $i$, $(U \cap K_i) \cup_{(S \cap K_i)} (V \cap K_i)$ is a relative Heegaard splitting of $K_i$ and if $(X \cap K_i)$ is correctly embedded in $X$ for $X = U,V$.  If $\{K_i\}$ is adapted to $S$ there is a subsequence which has the \defn{outer collar property} (Lemma \ref{outer collar property}). \newline

\item $\{K_i\}$ is \defn{perfectly adapted} to $S$ if it is adapted to $S$ and, in addition, each $\cl(K_{i+1} - K_i)$ is adapted to $S$.  That is, each $\cl(K_{i+1} - K_i)$ inherits a relative Heegaard splitting with Heegaard surface $S \cap \cl(K_{i+1} - K_i)$. \newline

\item $\{K_i = \eta(A_i \cup B_i)\}$ is a \defn{balanced exhausting sequence} for $M$ (adapted to $S$) if each $K_i$ is a regular neighborhood of $A_i \cup B_i$ where $A_i$ and $B_i$ are subcompressionbodies of $U$ and $V$ respectively with $\boundary_S A_i \subset \boundary_S B_i \subset \boundary_S A_{i+1}$.  \newline

\item $\{K_i\}$ is \defn{well-placed on $S$ rel $C$} if it is a frontier-incompressible (rel $C$) exhausting sequence for $M$ which is adapted to $S$ and, in addition, has the following properties:

\begin{enumerate}
\item[(WP1)] For each $i$, $V$ intersects each component of $\fr K_i$ in a single disc.
\item[(WP2)] For each $i$, $\fr K_i \cap U$ is incompressible in $U$.
\item[(WP3)] $\{K_i\}$ has the outer collar property with respect to $U$.
\item[(WP4)] For each $i$, no component of $\cl(M - K_i)$ is compact.
\end{enumerate}
\end{itemize}
\end{definition}

The main result of this section is:

\begin{theorem}\label{well-placed}
Suppose that $M$ is a non-compact orientable 3-manifold with compact boundary which is end-irreducible (rel $C$) where $C$ is a compact set containing $\boundary M$.  Suppose also that $U \cup_S V$ is an absolute Heegaard splitting of $M$.  Then there is an exhausting sequence of $M$ which is well-placed on $S$ rel $C$.
\end{theorem}

The most difficult part of the proof is in showing that there is a frontier-incompressible (rel $C$) exhausting sequence which is adapted to $S$.  

\subsection{Balanced Sequences and the Weakly Reducible Theorem}
We begin by showing that there is a balanced exhausting sequence of $M$ adapted to $S$ so that the compressing discs for the frontiers of the exhausting elements are in a ``good position" relative to the Heegaard surface.

\begin{proposition}\label{good discs exist}
There is a balanced exhausting sequence $\{C_i = \eta(A'_i \cup B'_i)\}$ for $M$ adapted to $S$ and a 2-sided disc family $\Psi$ for $S$ which contains $\cup_i (\fr A'_i \cup \fr B'_i)$ such that, for each $i$, $\sigma(\cl(\boundary_S B'_i - \boundary_S A'_i);\Psi)$ is incompressible in $M - C$.  
\end{proposition}

\begin{remark}In \cite{CaGo87} Casson and Gordon prove that if a Heegaard splitting of a compact 3-manifold is weakly reducible then there is a 2-sided disc family for the Heegaard surface such that when the surface is compressed along that family, the result is a collection of incompressible surfaces (possibly inessential spheres)\footnote{This is not how the result is usually stated, but see the proof given in \cite{Sc02}.}.  Since the frontiers of balanced submanifolds consist of surfaces which are obtained from the Heegaard surface by compressions along disjoint discs, it is natural to try to harness the power of the Casson and Gordon theorem.\newline

It is unclear, however, if the Casson and Gordon theorem can be extended to non-compact 3-manifolds in a way that is directly useful in this situation.  Nonetheless, the proof of our theorem is based on the outline of a proof of Casson and Gordon's theorem given in \cite{Sc02}.  We will also need to use Casson and Gordon's version of Haken's Lemma. 
\end{remark}

The proof is rather long so we begin with an outline of the proof:

\begin{enumerate}
\item\label{step 1} Take a balanced exhausting sequence $\{K_i = \eta(A_i \cup B_i)\}$.  For each $K_n$ show how to replace $K_{n-2}, K_{n-1},$ and $K_n$ with ``better" balanced submanifolds $K^L_{k} = \eta(A^L_k \cup B^L_k)$ for $k = n-2, n-1, n$.  Each of these better balanced submanifolds is still contained in $K_{n+1}$ and still contains $K_{n-3}$.  Let $C_n = K^L_n$.  The new manifolds will be obtained from the old ones by a finite sequence of slide-moves $L$.  The process of obtaining $C_n$ will also leave us with a 2-sided disc family $\Delta$ for $S \cap K_{n+1}$. \newline

\item\label{step 2} Suppose that there is a compressing disc $D$ for $\sigma(\cl(\boundary_S B^L_n - \boundary_S A^L_n); \Delta)$.  \newline

\item \label{step 3}Show that we can assume that $D$ is contained in $K_{n+1} - K^L_{n-2}$.  This step is where we use the eventual end-irreducibility of $M$. \newline

\item\label{step 4} Replace $D$ by a compressing disc of $\sigma(S \cap K_{n+1};\Delta)$ which intersects $\sigma(S \cap K_{n+1};\Delta)$ only on $\boundary D$.  We continue calling the disc $D$.\newline

\item\label{step 5} Use Haken's Lemma to replace $D$ by a disc which intersects a certain Heegaard surface exactly once and is contained in $K_{n+1} - K_{n-3}$.  We continue calling the disc $D$. \newline

\item\label{step 6} Follow the arguments of Casson and Gordon's Weakly Reducible theorem to obtain from $C_n$ by slide-moves a balanced submanifold which is even better than $C_n$.  This will contradict the construction of $C_n$. \newline

\item\label{step 7} Use this replacement technique on each element of a subsequence of $\{K_i\}$ to obtain the desired $\{C_i\}$.  Construct the 2-sided disc family $\Psi$ from the 2-sided disc families $\Delta$ which were created in each replacement operation.
\end{enumerate}

\begin{proof}[Proof of Proposition \ref{good discs exist}]

Let $\{K_i = \eta(A_i \cup B_i)\}_{i \geq 0}$ be a balanced exhausting sequence for $M$ adapted to $S$ and let $\{P_i\}$ be a frontier incompressible exhausting sequence (rel $C$).  Choose the exhausting sequences so that $C \subset K_0 \subset P_{i-1} \subset K_i \subset P_i$ for all $i \geq 1$.  Each of the inclusions should be into the interior of the succeeding submanifold. Figure \ref{frontier and balanced sequences} is a schematic of the exhausting sequences.  The frontiers of the submanifolds in $\{P_i\}$ may have a very complicated intersection with the Heegaard surface.  The frontier of each submanifold in the balanced exhausting sequence consists of discs and compact surfaces parallel to subsurfaces of $S$.\newline

\begin{figure}[ht]
\scalebox{0.6}{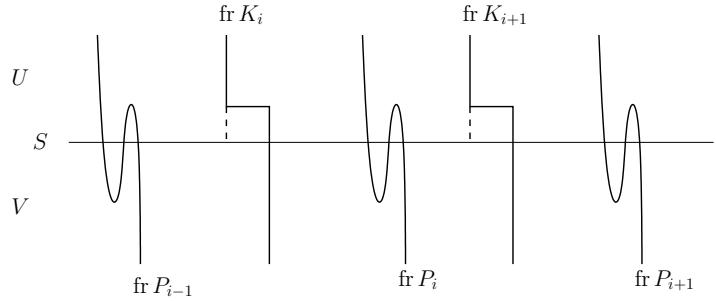}
\caption{A schematic of the exhausting sequences}
\label{frontier and balanced sequences}
\end{figure}

We will show that given a $q \in \N$ and $n \geq q + 3$, $K_n$ can be replaced by a compact submanifold $C_n = \eta(A'_n \cup B'_n)$ with the following properties: \newline

\begin{enumerate}
\item $C_n$ is obtained from $K_n$ by slide moves.\newline

\item There is a 2-sided disc family $\Delta$ for $S$ in $K_{n+1}$ containing $\fr A'_n \cup \fr B'_n$ such that $\sigma(\cl(\boundary_S B'_n - \boundary_S A'_n); \Delta)$ is incompressible in $M - C$. \newline

\item We still have $K_q \subset C_n$ and the discs $\fr A_q \cup \fr B_q$ are contained in $\Delta$. 
\end{enumerate}

Choose some $n \geq q+3$. 
\newline

Let $\ob{\Delta}= \bigcup_{q \leq i \leq n+1} (\fr A_i \cup \fr B_i)$.  Recall from Section \ref{slide-moves and isotopies} that a slide-move of this 2-sided disc family consists of either adding a compressing disc for $S$ to $\ob{\Delta}$ which is disjoint from all other discs of $\ob{\Delta}$ or performing a 2-handle slide of one disc of $\ob{\Delta}$ over another disc of $\ob{\Delta}$.  The arc over which a 2-handle slide is performed must have its interior disjoint from all discs of $\ob{\Delta}$.\newline

Recall from just after Lemma \ref{balanced are adapted} that each slide-move performed on $\ob{\Delta}$ leaves us with new balanced submanifolds obtained from the submanifolds $\{K_i\}_{i \leq n+1}$ by slide-moves.  After performing a slide-move, we still have $K_i \subset K_{i+1}$ for all $i$, since all the slides are performed relative to $\ob{\Delta}$.  \newline

Let $\mc{L}$ denote the set of all finite sequences of slide-moves of $\ob{\Delta}$ subject to the following restrictions: 

\begin{enumerate}
\item Every time a disc is added to $\ob{\Delta}$, the disc has boundary lying on $S \cap K_{n+1}$.\newline

\item No disc of $\fr K_{n+1} \cup \fr K_{q}$ is ever slid over another disc. 
\end{enumerate}

These restrictions mean that performing a sequence of slide-moves in $\mc{L}$ preserves the ordering of submanifolds $K_i$ for $q \leq i \leq n+1$. Furthermore, the manifolds $K_{n+1}$ and $K_q$ are left unchanged.  

\subsubsection*{Step \ref{step 1}} Each sequence $L \in \mc{L}$ leaves us with new balanced submanifolds $K^L_i$ for $q < i < n+1$.  The submanifolds $K_q$ and $K_{n+1}$ are left unchanged.  For ease of notation, let $K^L_q = K_q$ and $K^L_{n+1} = K_{n+1}$.  Let $A^L_i$ be the subcompressionbody of $U$ obtained from $A_i$ by the slide-moves $L$ and let $B^L_i$ be the subcompressionbody of $V$ obtained from $B_i$ by the slide-moves $L$ so that $K^L_i = \eta(A^L_i \cup B^L_i)$. \newline

Recall from \cite{CaGo87} that the complexity of a closed, connected surface $F$ is defined to be $1 - \chi(F)$, unless $F$ is a two-sphere, in which case, it is 0.  The complexity of a disconnected closed surface is the sum of the complexities of the components. \newline

Performing $L$ on $\ob{\Delta}$ leaves us with a disc family $\ob{\Delta}_L$ which contains the discs $\fr A^L_i \cup \fr B^L_i$ for $q \leq i \leq n+1$.  Define the complexity of $\ob{\Delta}_L$ to be the complexity of $\sigma(S \cap K_{n+1};\ob{\Delta}_L)$.  Since complexity is invariant under handle-slides (Lemma \ref{isotopies of handle slides}), the complexity of a 2-sided disc family cannot increase under slide-moves.\newline

Choose an $L \in \mc{L}$ so that $\ob{\Delta}_L$ has minimal complexity.  Let $\Delta = \ob{\Delta}_{L}$ and $C_n = K^L_n$.  Let $\Delta_1$ be those discs of $\Delta$ which lie in $U$ and $\Delta_2$ those discs which lie in $V$.  \newline

Recall that if $R \subset S$ is a compact subsurface of $S$ with $\boundary R \subset \boundary \Delta$, the notation $\sigma(R;\Delta)$ signifies the surface obtained from $R$ by compressing along those discs of $\Delta$ which have boundary on $R$.  Let $R_i = \cl(\boundary_S B_i - \boundary_S A_i)$ and let $R'_i = \cl(\boundary_S B^L_i - \boundary_S A^L_i)$ for $q \leq i \leq n+1$.  Note that $R'_i$ is obtained from $R_i$ by the sequence of slide-moves $L$.  We claim that $\sigma(R'_n; \Delta)$ is incompressible in $M - C$.  The surface $R'_i$ is a subsurface of $S$ which is parallel in $K^L_i$ to $\cl(\fr K^L_i - (\fr A^L_i \cup \fr B^L_i))$.

\begin{figure}[ht]
\scalebox{0.6}{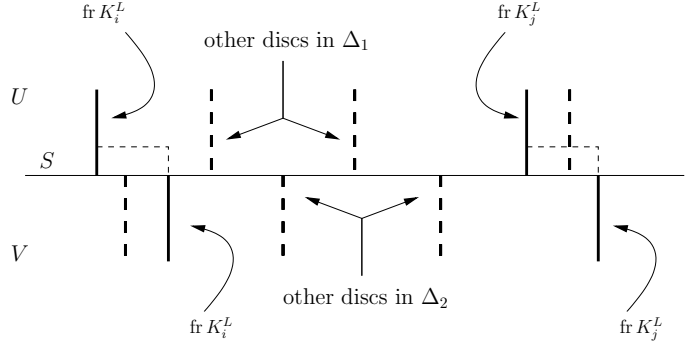}
\caption{A schematic of $\Delta_1$ and $\Delta_2$}
\end{figure}

\subsubsection*{Step \ref{step 2}} Let $S_k = \sigma(S \cap K_{n+1};\Delta_k)$ for $k = 1, 2$.  Let $W_1$ be $U \cap K_{n+1} $ together with the 2-handles coming from $\Delta_2$ minus the 2-handles coming from $\Delta_1$.  Let $W_2$ be the closure of the complement of $W_1$ in $K_{n+1}$.  Let $\ob{S} = \sigma(S \cap K_{n+1}; \Delta)$.  We are trying to show that $\sigma(R'_n;\Delta)$ is incompressible in $M - C$.  We assume the contradiction: suppose that a component $\ob{B}$ of $\sigma(R'_n; \Delta)$ is compressible in $M - C$.


\subsubsection*{Step \ref{step 3}} Our next task is to show that there is a compressing disc for $\ob{B}$ which lies entirely in $K_{n+1} - K^L_{n-2}$. Recall that $\{P_i\}$ is the frontier-incompressible (rel $C$) exhausting sequence for $M$ which is interlaced with $\{K_i\}$.  Let $\Sigma = \fr P_{n-1} \cup \fr P_n$.  The key technique is an application of Proposition \ref{isotopies of slide-moves} and the incompressibility in $M - C$ of $\Sigma$. \newline

By Proposition \ref{isotopies of slide-moves}, there is a proper ambient isotopy $f$ taking $\sigma(S;\Delta)$ to the surface obtained from $\sigma(S;\ob{\Delta})$ by compressing along a certain collection of discs.  In particular, there are disjoint collections of disjoint ordered discs $\mc{E}$ and $\mc{G}$ so that the discs of $\mc{E}$ have boundary on $\sigma(R_n;\ob{\Delta})$ and the discs of $\mc{G}$ have boundary on $\sigma(R_{n-1};\ob{\Delta})$ and the isotopy $f$ takes $\sigma(R'_n;\Delta)$ to $\sigma(R_n;\ob{\Delta} \cup \mc{E})$ and $\sigma(R'_{n-1};\Delta)$ to $\sigma(R_{n-1};\ob{\Delta} \cup \mc{G})$.  The notation $\sigma(R_n;\ob{\Delta}\cup \mc{E})$ means the surface obtained from $\sigma(R_n;\ob{\Delta})$ by compressing along the discs of $\mc{E}$ in the order given.  Similarly, we write $\sigma(R_{n-1};\ob{\Delta} \cup \mc{G})$ for the surface obtained from $\sigma(R_{n-1};\ob{\Delta})$ by compressing along $\mc{G}$. The surfaces $\sigma(R_n;\ob{\Delta} \cup \mc{E})$ and $\sigma(R_{n-1};\ob{\Delta}\cup \mc{G})$ are disjoint. \newline

The discs $\mc{E}$ have boundary on $\sigma(R_n;\ob{\Delta})\subset (P_n - P_{n-1})$.  As $\Sigma$ is incompressible in $M - C$ the intersections of $\mc{E}$ with $\Sigma$ are inessential on $\Sigma$.  Similarly, the discs of $\mc{G}$ have boundary on $\sigma(R_{n-1};\ob{\Delta}) \subset P_{n-1}$ and so $\mc{G}$ intersects $\Sigma$ in loops which are inessential on $\Sigma$.  The surface $\ob{B}$ is taken by the isotopy $f$ to a surface $\ob{K}$ which is a component of $\sigma(R_n;\ob{\Delta}\cup \mc{E})$. Since $\ob{B}$ is compressible in $M - C$ so is $\ob{K}$. \newline

Since the intersections of $\ob{K}$ with the incompressible (in $M - C$) $\Sigma$ come from the intersections of $\mc{E}$ with $\ob{K}$, the loops $\ob{K} \cap \Sigma$ are inessential on both surfaces. There is, therefore, a surface $K' \subset (P_n - P_{n-1})$ which is obtained from $\ob{K}$ by cutting and pasting along the intersections $\ob{K} \cap \Sigma$. (Start with innermost discs of intersection on $\Sigma$ and replace the corresponding discs of $\ob{K}$ with copies of the discs on $\Sigma$ which have been pushed slightly into $(P_n - P_{n-1})$.)  As $\ob{K}$ is compressible in $M - C$, $K'$ is also compressible in $M - C$. Since $\Sigma$ is incompressible in $M - C$ there is a compressing disc $F$ for $K'$ which is contained in $P_n - P_{n-1}$. Our goal is to use $F$ to construct a compressing disc for $\ob{K}$ which is disjoint from $\sigma(R_{n-1};\ob{\Delta} \cup \mc{G})$. \newline

Since $\boundary \mc{E}$ consists of inessential loops on $K'$ we may assume that $\boundary F \cap \boundary \mc{E} = \nil$.  The disc $F$ may intersect the discs $\mc{E}$.  It may also intersect the discs of $\mc{G}$ in simple closed curves.  Since each loop of $F \cap \mc{E}$ is inessential on both $F$ and $\mc{E}$ we may, by cutting and pasting $F$ along the intersections, obtain a compressing disc $F'$ for $\ob{K}$.  Since both $K'$ and $\mc{E}$ were disjoint from $\sigma(R_{n-1};\ob{\Delta} \cup \mc{G})$, any intersections of the disc $F'$ with the surface $\sigma(R_{n-1};\ob{\Delta} \cup \mc{G})$ occur because $F'$ intersects $\mc{G}$ in simple closed curves.  These intersections are inessential on both $F'$ and on $\sigma(R_{n-1};\ob{\Delta}\cup \mc{G})$.  We may cut and paste $F'$ along these intersections to produce a compressing disc $E$ for $\ob{K}$ which is disjoint from $\sigma(R_{n-1};\ob{\Delta}\cup \mc{G})$.  The disc $E$ may intersect $\Sigma$, but that is not of concern. \newline

Reversing the isotopy $f$ takes $E$ to a compressing disc $D$ for $\ob{B}$. $D$ is contained in $K_{n+1}$.  The disc $D$ is disjoint from $\sigma(R'_{n-1};\Delta)$ since $E$ was disjoint from $\sigma(R_{n-1};\ob{\Delta} \cup \mc{G})$. \newline

Recall that we are trying to construct a compressing disc for $\ob{B}$ which is contained in $K_{n+1} - K^L_{n-2}$.  Each disc of $\Delta$ which had boundary on $R'_{n-1}$ was disjoint from $R'_{n-2}$ since no disc of $\Delta$ intersects $S$ except at its boundary and the discs of $\Delta$ are pairwise disjoint.  Thus $K^L_{n-2}$ is contained inside some component of $\sigma(K^L_{n-1};\Delta)$.  But since $D$ is disjoint from $\sigma(R'_{n-1};\Delta)$ which is parallel to $(\fr \sigma(K^L_{n-1};\Delta))$, $D$ can be isotoped so as to not intersect $K^L_{n-2}$.  Hence, there is a compressing disc $D$ for $\ob{B}$ which is contained in $K_{n+1} - K^L_{n-2}$.

\subsubsection*{Step \ref{step 4}} The compressing disc $D$ may intersect the surface $\ob{S} \cap (K_{n+1} - K^L_{n-2})$. By revising the disc $D$ we may assume that no loops of $D \cap \ob{S}$ are inessential on $\ob{S}$.  Replace $D$ by an innermost disc, which we will continue to call $D$, that intersects $\ob{S}$ only on $\boundary D$.  By our construction $D$ is now a compressing disc for $\ob{S}$.  The boundary of $D$ may no longer be on $\ob{B}$.  $D$ lies in either $W_1$ or $W_2$ and is completely contained in $(K_{n+1} - K^L_{n-2})$. Recall that $W_1 = [(U \cap K_{n+1}) - \eta(\Delta_1) \cup \eta(\Delta_2)]$ and that $W_2 = [(V \cap K_{n+1}) - \eta(\Delta_2) \cup \eta(\Delta_1)]$.


\subsubsection*{Step \ref{step 5}}Our goal is to use the disc $D$ to construct a sequence $L' \in \mc{L}$ such that $\ob{\Delta}_{L'}$ has lower complexity than $\Delta = \ob{\Delta}_L$.  This will contradict our original choice of $L$.  As mentioned in the remark preceding this proof, the strategy is to follow the outline of the proof of Casson and Gordon's Weakly Reducible theorem given in \cite{Sc02}.  We will view $S_1$ as a Heegaard surface for $W_1$ or $S_2$ as a Heegaard surface for $W_2$ depending on which side the disc $D$ lies.  In the Casson and Gordon theorem the two cases had identical arguments.  Here, however, the relationship of $W_1$ and $W_2$ to $K_{n+1} - K_q$ is not symmetric due to the asymmetry in the construction of balanced submanifolds.  We will briefly need to consider the two cases separately.  We will eventually be able to combine arguments. 

\begin{remark}
Some care is needed when we consider $S_1$ or $S_2$ has a Heegaard surface, as $\ob{S}$ may contain spheres.  This means that the compressionbodies we are considering may not be irreducible.  This does not really affect the proofs as the only times we would want to use the irreducibility of a compressionbody is when we isotope (in a compressionbody) one disc past another which shares its boundary.  If $\ob{S}$ contains spherical components which get in the way of the isotopy, we may first perform a surgery on the disc we want to isotope so that the two discs with common boundary bound a 3-ball and then perform the isotopy.  We will refer to this process as \defn{revising and isotoping} the disc which, if $\ob{S}$ were irreducible, we would have merely isotoped.
\end{remark}

Suppose, first, that $D$ lies in $W_1$.  By pushing $\ob{S}$ slightly into $W_2$ we can view $S_1$ as a Heegaard surface for the disconnected 3-manifold $W_1$.  $S_1$ divides $W_1$ into (disconnected) absolute compressionbodies $U'$ and $V'$.  Let $V'$ be the absolute compressionbody containing $\ob{S}$.  See Figure \ref{GoodDiscsExist3}.  The disc $D$ is a compressing disc for $\boundary W_1$. \newline   

\begin{figure}[ht]
\scalebox{0.6}{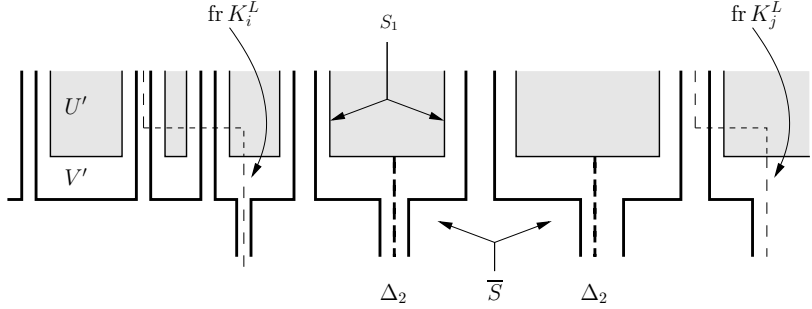}
\caption{$S_1$ as a Heegaard surface for $W_1$.}
\label{GoodDiscsExist3}
\end{figure}

We can apply Haken's Lemma to obtain a compressing disc $D'$, a compressing disc for $\ob{S}$ in $W_1$, which intersects $S_1$ in a single loop and is such that $\boundary D' = \boundary D$.  $W_1 \subset K_{n+1}$ by the definition of $W_1$, so $D'$ does not intersect $\fr K_{n+1}$.  The discs $\fr A^L_{n-2}$ are in $\Delta_1$ and separate $U$.  Thus no component of $W_1$ intersects both $\fr K^L_{n-2}$ and $\fr K^L_{n-3}$. Hence, since $\boundary D$ is in $W_1 \cap (K_{n+1} - K^L_{n-2})$ the disc $D'$ is in $K_{n+1} - K^L_{n - 3}$.  Summarizing: $D'$ is a compressing disc for $\ob{S}$ which intersects $S_1$ in a single loop and is contained in $K_{n+1} - K^L_{n-3}$. \newline  

We now turn to the case when $D \subset W_2$. Push $\ob{S}$ slightly into $W_1$ and view $S_2$ as a Heegaard surface for the 3-manifold $W_2$.  The disc $D$ is a compressing disc for $\boundary W_2$. Let $U'$ and $V'$ be the submanifolds of $W_2$ into which $S_2$ divides $W_2$.  $U'$ is the submanifold which has $\ob{S}$ as its boundary.  See Figure \ref{GoodDiscsExist4}. \newline

\begin{figure}[ht]
\scalebox{0.6}{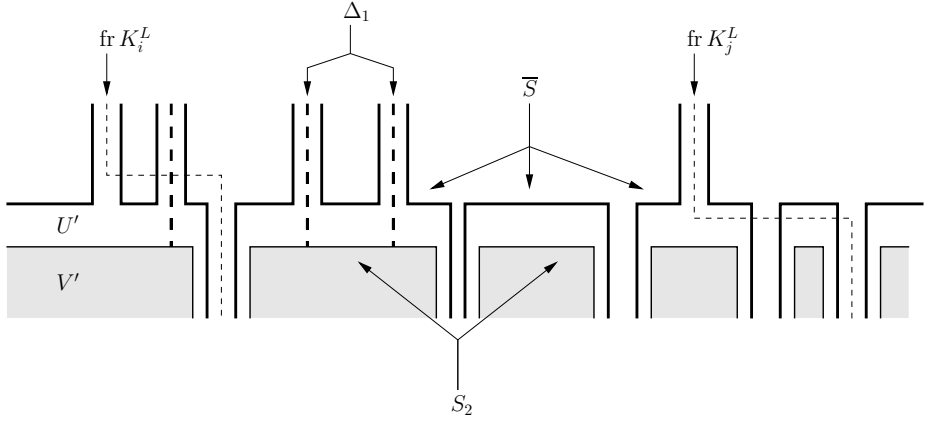}
\caption{$S_2$ as a Heegaard surface for $W_2$.}
\label{GoodDiscsExist4}
\end{figure}

The discs of $\fr B^L_{n-2}$ are contained in $\Delta_2$ and separate $V$.  Thus no component of $W_2$ intersects both $(K_{n+1} - K^L_{n-2})$ and $\interior K^L_{n-2}$.  The disc $D$ is a compressing disc for $\ob{S} \subset \boundary W_2$ which is contained in a component of $W_2$ disjoint from $\interior K^L_{n-2}$.  Applying Haken's Lemma, we can replace $D$ with a disc $D'$ such that $\boundary D' = \boundary D$ and $D'$ intersects $S_2$ in a single loop.  Since $D$ and $D'$ are in the same component of $W_2$, $D' \cap \interior K^L_{n-2} = \nil$.  Summarizing: The disc $D'$ is a compressing disc for $\ob{S}$ which intersects $S_2$ in a single loop and is contained in $K_{n+1} - K^L_{n-2}$.

\subsubsection*{Step \ref{step 6}}Recall that $2 < q \leq (n-3)$.  We may now combine arguments.  In the previous step, we showed that there was a compressing disc for $\ob{S}$ which is located in either $W_1 \cap \cl(K_{n+1} - K_q)$ or $W_2 \cap \cl(K_{n+1} - K_q)$. and intersects $S_1$ or $S_2$ (respectively) in a single loop $\gamma$.  We will now produce a sequence of slide-moves $l$ such that the sequence of slide-moves $L$ followed by $l$ is in $\mc{L}$ and the such that the sequence $L$ followed by $l$ has lower complexity than $L$.  This will contradict our choice of $L$.  The difficult part of this step is nearly identical to Bonahon's proof of Proposition \ref{choosing discs}.  This is Proposition B.1 of \cite{Bo83}.  We include the proof here, because we need to pay careful attention to the type of slide-moves which are required. \newline

Without loss of generality, suppose that $D$ is a compressing disc for $\ob{S}$ which is located in $W_1 \cap \cl(K_{n+1} - K_q)$ and intersects $S_1$ in a single loop $\gamma$. (We were calling this disc $D'$ in the previous step.)  We continue to view $S_1$ as a Heegaard surface for $W_1$.  Recall that $V'$ denotes the compressionbody which is the region between $\ob{S}$ and $S_1$ and that $U'$ is the closure of the complement of $V'$ in $W_1$. See Figure \ref{GoodDiscsExist3}. \newline

We may assume that $D$ is disjoint from the discs of $\Delta_1$; it may, however, intersect  the discs $\Delta_2$ (including the frontiers of some $B^L_i$ (for $q < i < n+1$).  Let $A$ denote the annulus $D \cap V'$ and $D'$ the disc $D \cap U'$.  Consider how $A$ intersects $\Delta_2$.    \newline

By an innermost disc argument we may assume that the annulus $A$ intersects the discs of $\Delta_2$ entirely in arcs with both endpoints on $\gamma$.  Let $a$ be an outermost arc of intersection on $A$.  Let $b$ be the arc of $\gamma$ with endpoints $\boundary a$ which intersects no disc of $\Delta_2$.  Let $G$ be the disc of $\Delta_2$ such that $a \subset G \cap A$.  Let $c$ be an arc of $\boundary G$ which has endpoints $\boundary a$.  The arc $c$, of course, may have other intersections with $\gamma$. \newline

Combining the subdiscs of $A$ and $G$ with boundaries $a \cup b$ and $a \cup c$ respectively and pushing off $\Delta_2$ a little, we obtain a compressing disc for $S_1$ in $V'$ which is disjoint from the complete collection of discs $\Delta_2$ for $V'$.  Thus $b \cup c$  is a loop bounding a disc $Q$ in $\sigma(S_1;\Delta_2)$. (We are not calling this surface $\ob{S}$ since we have pushed $\ob{S}$ into $W_2$.)\newline

We now adapt Bonahon's proof of Proposition \ref{choosing discs} to show that we can perform 2-handle slides of $G$ over the discs of $\Delta_2$ which have boundary in $Q$ and then revise and isotope $D$ to remove all intersections of $D$ with $G$ (see remark in Step 5 about the term ``revise and isotope"). When we compress $S_1$ along $\Delta_2$, the remnants of $\Delta_2$ show up as spots, some of which are in the interior of the disc $Q$.  Each disc of $\Delta_2$ contributes two spots to $\sigma(S_1;\Delta_2)$.  For each spot $F_i$ from $\Delta_2$ which shows up in $Q$, excluding a possible spot coming from $G$, choose oriented arcs $\alpha_i$ contained in $Q$ joining $G$ to the discs of $\Delta_2$ giving rise to those spots.  If a disc of $\Delta_2$ produces two spots contained in $Q$ then we have two oriented arcs joining $G$ to that disc.  Choose the arcs $\alpha_i$ so that $\alpha_i \cap \Delta_2 = \boundary \alpha_i$ and so that the $\{\alpha_i\}$ are pairwise disjoint.  The arcs $\alpha_i$ lie on $S_1$ and for each arc $\alpha_i$ we may perform a handle-slide of $G$ over the the disc to which it is joined by $\alpha_i$.  Continue calling this disc $G$. By performing these slides, we may have increased the number of intersections between $\boundary G$ and $\gamma$.  These handle-slides convert $Q$ into a new disc as the arc $c$ is changed by the handle-slides.  We continue calling the disc $Q$. After these handle-slides $Q$ contains no spots from $\Delta_2$, except perhaps one coming from $G$.  Revise and isotope $D$ (rel $b$) so that $\gamma$ has minimal intersection with $\boundary G$.  Suppose, now, that the disc $Q$ contains a spot arising from $G$.  Let $G_1$ and $G_2$ denote the two spots.  Since they both arise from $G$ we have that $|\gamma \cap \boundary G_1| = |\gamma \cap \boundary G_2|$.  Since any arc of $\gamma$ with both endpoints on $\boundary G_2$ would bound a disc in $S_1$ and could, therefore, be removed by revising and isotoping $D$, each arc of $\gamma$ with an endpoint on $\boundary G_2$ also has an endpoint on $\boundary G_1$.  However $\boundary b \subset \boundary G_1$ and thus $|\gamma \cap \boundary G_1| = |\gamma \cap \boundary G_2| + 2$.  This, however, contradicts the earlier equation and so the spot $G_2$ cannot exist in $Q$.  The disc $Q$, therefore, is now a disc in $S_1$ and we can revise and isotope $D$ to remove the intersection $a$ from $D \cap \Delta_2$.  Since we have previously removed all other intersections, including the intersections introduced earlier, of $c$ with $\boundary \Delta_2$ we have decreased $|D \cap \Delta_2|$ by at least one.  Hence, by induction, we can remove all intersections of $D$ with $\Delta_2$ by revising and isotoping $D$ (rel $\boundary D$) and handle-sliding $\Delta_2$. \newline

This produces a new disc set $\Delta'_2$ which is disjoint from $\Delta_1 \cup \{D'\}$.  At the beginning of the process the curve $\gamma$ does not intersect any disc of $\fr K_{n+1} \cup \fr K_q$.  The set of discs with boundary in $R$ may contain discs that are associated to discs of $\fr K_{n+1} \cup \fr K_q$, but we were able to choose our sliding arcs so that they only intersected the discs of $\fr K_{n+1} \cup \fr K_{q}$ in at most one endpoint.  The only slides we performed were of the disc $G$ over other discs, and since $\gamma$ intersected $G$, $G$ was not a disc of $\fr K_{n+1} \cup \fr K_{q}$.  Furthermore, since the discs of $\Delta_1$ show up as spots on $S_1$ it is easy to arrange these slides to be relative to $\Delta_1$.  Thus, these handle-slides are of the sort allowed in sequences in $\mc{L}$.  Let $l$ denote the sequence of these handle-slides followed by the slide-move (M2) where we add the disc $D'$ to $\Delta_1$.  The sequence of slide-moves consisting of $L$ followed by $l$ does, therefore, give us a sequence of slide-moves in the collection $\mc{L}$.  As $D$ was a compressing disc for $\ob{S}$ this sequence of slide-moves has lower complexity than our original choice from $L$. This, however, contradicts our initial choice $L$ to be such that the complexity of $\sigma(S \cap K_{n+1}; \Delta)$ was minimal. The contradiction arises from our assumption that $\ob{B}$ is compressible: therefore, $\ob{B}$ is incompressible in $M - C$.  
\subsubsection*{Step \ref{step 7}} Recall that $\{K_i\}$ is our balanced exhausting sequence adapted to $S$ which is interspersed with a frontier-incompressible (rel $C$) exhausting sequence.  Let $q_n = 5n$ for $n \geq 2$.  We have shown how to replace $K_{q_n}$ with a balanced submanifold $C_n = \eta(A'_n \cup B'_n)$ which contains $K_{q_n - 3}$.  The sequence $\{C_n\}$ is a balanced exhausting sequence adapted to $S$.  In the construction of each $C_n$ we also constructed a 2-sided disc family $\Delta$ so that $\sigma(\cl(\boundary_S B'_n - \boundary_S A'_n);\Delta)$ is incompressible in $M - C$.  Let $\Delta_n$ denote those discs of $\Delta$ with boundary on $\cl(\boundary_S B'_n - \boundary_S A'_n)$.  Note that $\Delta_n$ is disjoint from $\Delta_i$ for all $i < n$.  Let $\Psi = \cup_n \Delta_n$.  $\Psi$ is a 2-sided disc family for $S$ where each disc of $\Psi$ has boundary on the frontier of some $C_n$.  When we compress $\cup \cl(\boundary_S B'_n - \boundary_S A'_n)$ along $\Psi$ we obtain surfaces which are incompressible (rel $C$).
\end{proof}
\subsection{The Proof of Theorem \ref{well-placed}}
Recall that $M$ is end-irreducible rel $C \supset \boundary M$ and that $U \cup_S V$ is an absolute Heegaard splitting for $M$. 

\begin{proposition}\label{fr incomp and adapted}
$M$ has a frontier-incompressible (rel $C$) exhausting sequence which is adapted to $S$.  Furthermore, $V$ intersects the frontier of each element of the exhausting sequence in discs.
\end{proposition}

\begin{proof}
Let $\{C_i = \eta(A'_i \cup B'_i)\}$ be the balanced exhausting sequence guaranteed by Proposition \ref{good discs exist}.  By the construction of balanced submanifolds, $V$ intersects each $\fr C_i$ in discs.  \newline

Proposition \ref{good discs exist} guarantees the existence of a 2-sided disc family $\Psi$ for $S$ such that $\cup_i (\fr A'_i \cup \fr B'_i) \subset \Psi$ and each $\sigma(\cl(\boundary_S B'_i - \boundary_S A'_i); \Psi)$ is incompressible in $M - C$.  Let $\Psi_1 = \Psi \cap U$ and $\Psi_2= \Psi \cap V$.  We may use the product region between $(\fr C_i - (\fr A_i \cup \fr B_i))$ and $\boundary_S B'_i - \boundary_S A'_i$ to extend the discs in $\Psi_2$ with boundary on $(\boundary_S B'_i - \boundary_S A'_i)$ to have boundary on $\fr C_i$.  \newline

If we boundary-reduce $C_i$ along $\Psi_2$ and add the 2-handles $\eta(\Psi_1)$ to $\fr C_i$ we end up with a new submanifold $\ob{C}_i$ of $M$.  By construction, the discs of $\Psi$ with boundary on $\fr C_i$ are disjoint from $\fr C_{i-1} \cup \fr C_{i+1}$.  Hence, $C_i$ is contained in a single component of $\ob{C}_{i+1}$ and so $M = \cup \ob{C}_i$.  Let $K_1$ be the component of $\ob{C}_2$ containing $C$ and, for each $n > 1$, let $K_n$ be the component of $\ob{C}_{n+1}$ containing $C_n$.  Since $C_n \subset K_n$ the sequence $\{K_n\}$ is an exhausting sequence for $M$.  Since the frontier of each $K_i$ is incompressible in $M - C$, the sequence $\{K_i\}$ is frontier-incompressible rel $C$. \newline

When we boundary-reduce $C_i$ along $\Psi_2$ we are boundary-reducing $C_i$ along disjoint discs which each intersect the relative Heegaard surface $S \cap C_i$ in a single simple closed curve.  By Haken's Lemma (Lemma \ref{Haken}), the resulting submanifold still has its intersection with $S$ a relative Heegaard surface.  When we add the 2-handles $\Psi_1$ to $\fr C_i$ we are adding 2-handles to $\boundary_- (U \cap C_i)$.  Hence, the resulting submanifold still has a relative Heegaard splitting coming from its intersection with $S$, apart from the introduction of 2-sphere components to $\boundary_-(U \cap C_i)$.  If there are any, we may add to $U \cap K_i$ the 3-balls bounded by those 2-spheres in $U$.  After we have added these 3-balls, $\{K_i\}$ is a correctly embedded exhausting sequence.  Therefore, $\{K_i\}$ is adapted to $S$.  Since the sequence is also frontier-incompressible (rel $C$) the proposition is proved.
\end{proof}

We now embark on proving that there is a frontier-incompressible (rel $C$) exhausting sequence for $M$ which is adapted to $S$ and has properties (WP1), (WP2), (WP3), and (WP4) in the definition of ``well placed exhausting sequence".  

\begin{lemma}\label{single disc}
Let $\{K_i\}$ be a frontier-incompressible (rel $C$) exhausting sequence for $M$ which is adapted to $S$.  Suppose that $V$ intersects $\fr K_i$ in discs for each $i$.  Then after taking a subsequence of $\{K_i\}$ and performing a proper ambient isotopy of $\cup_i \fr K_i$ we may arrange that $V$ intersects each component of each $\fr K_i$ in a single disc.  Additionally, $\{K_i\}$ has the outer collar property.  
\end{lemma}

\begin{proof}
Begin by taking a subsequence of $\{K_i\}$ such that $\{K_i\}$ has the outer collar property.  Let $K = K_j$ (for $j \geq 2$) be an element of this revised exhausting sequence.  Suppose that $B$ is a component of $\fr K$ such that $|V \cap B| \geq 2$.  We will describe an ambient isotopy of $\fr K$ which is the identity outside of $\cl(K_{j+1} - K_{j-1})$ to reduce the number of components of $|B \cap V|$ by one.  We may then perform this ambient isotopy on each element of $\{K_{2i}\}$ as needed in order to arrange that $V$ intersects each component of $\fr K_{2i}$ in a single disc.  The union of these isotopies is a proper ambient isotopy of $\{K_{2i}\}$.  After performing this isotopy, it will be clear that $\{K_{2i}\}$ still has the outer collar property.   \newline  

Let $B' = U \cap B$.  Since $V \cap B$ consists of discs, $B'$ is connected and has at least two boundary components.  $B'$ makes up part of the frontier of the relative compressionbody $K \cap U$.  $B'$ is a component of $\boundary_- (K \cap U)$ since $\{K_i\}$ is a correctly embedded exhausting sequence.  Since $\{K_i\}$ has the outer collar property, there is a product region $P = B' \times I$ which is embedded in $\cl((K - (K_{j-1}) \cap U))$ such that $B' = B' \times \{0\}$ and $B' \times \{1\}$ is a subsurface of $S \cap K$ except at a finite number of open discs $\delta$.  Choose an arc $\alpha \subset B' \times \{1\}$ so that $\alpha \cap \boundary (B' \times \{1\}) = \boundary \alpha$, $\alpha$ joins different components of $\boundary (B' \times \{1\})$, and $\alpha$ is disjoint from the discs $\delta$.  Let $D = \alpha \times I \subset P$ so that $\alpha = \alpha \times \{1\}$.  $D$ is an embedded disc in $P$ such that $\boundary D$ is composed of two arcs, one on $B'$ and one on $S \cap K$.  Isotope $B \cap \eta(D)$ across the disc $D$.  After this isotopy, the number of intersections $B \cap S$ has been reduced by one.  \newline

We now inspect the effect of this isotopy on $V \cap K$ and $U \cap K$.  In $V \cap K$ we have changed $\boundary_- V$ by banding together two discs.  Since $V \cap K$ was a relative compressionbody with $\boundary_- (V \cap K)$ consisting of discs, we have not changed the homeomorphism type of $V \cap K$, we have changed only the preferred surface. \newline

The effect of the isotopy on $U \cap K$ is to replace $B' \times I$ with $C' \times I$ where $C'$ is the surface obtained from $B'$ by removing a neighborhood of an arc joining two components of $\boundary B'$.  Clearly, $U \cap K$ is still a relative compressionbody with preferred surface $S \cap K$.  Furthermore, the presence of the product region $C' \times I$ shows that the sequence $\{K_i\}$ still has the outer collar property.  The isotopy we have described is the identity outside of $K_{j+1} - K_{j-1}$.  
\end{proof}

\begin{proof}[Proof of Theorem \ref{well-placed}]
Take the exhausting sequence $\{K_i\}$ given by Lemma \ref{single disc}.  The only properties we have left to achieve are (WP2) and (WP4).  We now prove that we have, in fact, already achieved (WP2) and that we can achieve (WP4) without ruining the others. \newline

Suppose that $B$ is some component of $\fr K_i$ such that $B \cap U$ has a compressing disc $D$ which is contained in $U$.  Since $K_i \cap U$ is a relative compressionbody and $(B \cap U) \subset \boundary_- (K_i \cap U)$, the compressing disc $D$ must be on the outside of $K_i$.  The curve $\boundary D$ bounds a disc $E \subset B$ since $B$ is incompressible in $M - C$ and $C \subset K$.  Since $D$ is a compressing disc for $B \cap U$, the disc $E$ is not contained in $B \cap U$.  Thus $(V \cap B) \subset E$.  Forming $K'_i$ by adding $\eta(D)$ to $K_i$ cuts $B$ into two surfaces: $B'$ which is homeomorphic to $B$ and $B''$ which is a 2-sphere.  Note that both $B'$ and $B''$ are components of $\boundary K'_i$.  The surface $B'$ is contained in $U$ and the sphere $B''$ intersects $V$ in a single disc.\newline

Since $B$ was incompressible in $M - C$ and $B'$ was obtained from $B$ by cutting off a 2-sphere, $B'$ is also incompressible in $M - C$.  The surface $B' \subset U$ is closed and incompressible in $U$.  Hence, $B'$ is parallel to a component of $\boundary_- U \subset \boundary M$.  This product region has boundary consisting of two components both of which are components of $\boundary K'_i$.  Thus the product region is actually $K'_i$.  But $B''$ is also a component of $\boundary K'_i$, so this is a contradiction.  Hence, $B \cap U$ is incompressible in $U$. Thus $\{K_i\}$ satisfies (WP2). \newline

Finally, we need to achieve (WP4).  Suppose that $\cl(M - K_1)$ has a compact component $L$.  There is some $K_n$ so that every compact component of $\cl(M - K_1)$ is contained in $K_n$.  By Corollary \ref{complementary compressionbodies}, $U \cap L$ and $V \cap L$ are relative compressionbodies.  Since there are no closed components of $\boundary_- (U \cap L)$ or $\boundary_-(V \cap L)$, both are also handlebodies.  Let $Q = L \cap K_1$.  $Q \cap U$ is an incompressible surface in $U$ which makes up part of $\boundary_- (U \cap K_1)$.  Choose a collaring set of discs $\delta$ for $U \cap K_1$.  Boundary-reducing $K_1 \cap U$ along $\delta$ leaves us with components homeomorphic to $(Q \cap U) \times I$.  Let $L' = (L \cap U) \cup ((Q \cap U) \times I)$.  This does not change the homeomorphism type of $L \cap U$, so $L'$ is a handlebody.  We may now reassemble $K_1 \cap U$ by attaching 1-handles corresponding to the discs $\delta$.  When we do this, we are attaching the handlebody $L'$ to the $\boundary_+$ of a relative compressionbody and so the result is a relative compressionbody with preferred surface $S \cap ((K_1 \cap U) \cup L')$.  Since $V$ intersected each component of $B$ in a single disc, $V \cap L$ is a handlebody and so $V \cap (K_1 \cup L)$ is also a relative compressionbody with preferred surface $S \cap (K_1 \cup L)$.\newline

Thus, if we include each compact component of $\cl(M - K_1)$ into $K_1$ to form $K'_1$ we still have a relative Heegaard splitting $K'_1 = (U \cap K'_1) \cup_{S \cap K'_1} (V \cap K'_1)$.  Assume that we have defined $K'_j$ for $j \geq 1$.  There exists an $n_j$ so that $K'_j \subset K_{n_j}$.  Let $K'_{j+1}$ be the union of $K_{n_j}$ and all of the compact components of $\cl(M - K_{n_j})$.  By the previous argument, $S$ gives a relative Heegaard splitting of $K_{j+1}$.  In such a way we obtain an exhausting sequence $\{K'_n\}$ for $M$ with property (WP4).  It is clear from the construction that $\{K'_n\}$ is, in fact, an exhausting sequence well-placed on $S$.
\end{proof}

\begin{remark}
Theorem \ref{well-placed} tells us that there is a frontier-incompressible (rel $C$) exhausting sequence $\{K_i\}$ for $M$ such that each $K_i$ inherits a relative Heegaard splitting from $U \cup_S V$.  If for some $j \geq 2$ each component of $\cl(K_{j+1} - K_{j})$ intersects $K_j$ in a connected surface, then examining the structure of the absolute Heegaard splitting of $K_{j+1}$ induced by the relative Heegaard splitting coming from $S$, shows that this absolute Heegaard splitting is obtained by amalgamating Heegaard splittings of $K_j$ and each component of $\cl(K_{j+1} - K_j)$.  
\end{remark}

\section{Heegaard Splittings of Deleted Boundary 3-Manifolds}\label{Deleted Boundary}

\subsection{Introduction}
\begin{definition}
A 3-manifold $M$ is \defn{almost compact} if there is a compact 3-manifold $\ob{M}$ with non-empty boundary and a non-empty closed set $J \subset \boundary \ob{M}$ such that $M$ is homeomorphic to $\ob{M} - J$.  If $J$ is the union of components of $\boundary \ob{M}$ then $M$ is a \defn{deleted boundary manifold}.
\end{definition}

Let $M$ be a deleted boundary manifold obtained from the compact manifold $\ob{M}$ by removing the union $J$ of boundary components.  By removing an open collar neighborhood of $J$ from $\ob{M}$ we obtain a compact manifold $C$ which  resides in $M$.  The closure of $M - C$ is homeomorphic to $J \times \R_+$.  Since $J$ is the union of components of $\boundary \ob{M}$, $J$ is a closed, possibly disconnected, surface.  $M$ is obviously end-irreducible (rel $C$) and $\boundary M \subset C$.  We will also assume that $\boundary M$ contains no spherical components, but, except where noted, $J$ may have spherical components. If $|J| \geq 2$ and if at least one component is a sphere, $M$ has Heegaard splittings which have infinitely many properly embedded reducing balls but not end-stabilized.  The following definitions (which make sense even when $M$ is not a deleted boundary 3-manifold) assist the classification in this case.

\begin{definition}
Let $e$ be an end of $M$ represented by non-compact submanifolds $\{W_i\}$ such that $\cl(W_i)$ is non-compact, $W_{i+1} \subset W_{i}$ for all $i$, and $M = \cup (M - W_i)$.  A Heegaard splitting $M = U \cup_S V$ is \defn{$e$-stabilized} if for each $i$ there is a reducing ball for $S$ contained in $W_i$.  Recall that $M$ is \defn{infinitely-stabilized} if it is $e$-stabilized for some end $e$ and $\defn{end-stabilized}$ if it is $e$-stabilized for every end $e$.
\end{definition}

The notion of being $e$-stabilized is a proper ambient isotopy invariant, as the next lemma shows.

\begin{lemma}\label{proper ambient isotopy invariant}
Suppose that $S$ and $T$ are Heegaard surfaces for $M$.  If there is an end $e$ of $M$ such that $S$ is $e$-stabilized but $T$ is not then $S$ and $T$ are not properly ambient isotopic. 
\end{lemma}

\begin{proof}
This follows directly from the fact that including a Heegaard surface into $M$ induces a homeomorphism on ends (Proposition \ref{end homeomorphism}) and that proper ambient isotopies fix each end of a manifold.
\end{proof}

\begin{definition}
Suppose that $U_S \cup_S V_S$ and $U_T \cup_T V_T$ are two absolute Heegaard splittings of $M$.  Then they are \defn{approximately isotopic} if for any compact set $C$ there are proper ambient isotopies of $S$ and $T$ so that $S \cap C = T \cap C$.
\end{definition}

The goal of this section is to completely classify Heegaard splittings of $M$ up to proper ambient isotopy and up to approximate isotopy. In particular, if $J$ contains no spherical components, $M$ has, up to proper ambient isotopy, exactly one Heegaard splitting and that splitting is end-stabilized.\newline

The following three theorems provide key ingredients in the classification.  

\begin{theorem}[Reidemeister-Singer]\label{Reidemeister-Singer}
After finitely many stabilizations, any two absolute Heegaard splittings of a compact 3-manifold which have the same partition of boundary are ambient isotopic.
\end{theorem}

The next is a version of Theorem 2.1 of \cite{FrMe97}.  A proof is provided in the Appendix. 

\begin{theorem}[Frohman-Meeks]\label{Frohman-Meeks}
Any two end-stabilized absolute Heegaard splittings with the same partition $\boundary M$ are properly ambient isotopic.  Any two infinitely-stabilized Heegaard splittings with the same partition of $\boundary M$ are approximately isotopic.
\end{theorem}

The following is the most involved result of this section. Its proof uses Scharlemann and Thompson's classification of splittings of $\text{(closed surface)} \times I$. \newline

Let $W_1, \hdots, W_n$ denote the components of $\cl(M - C)$ and let $X_1, \hdots, X_n$ denote the components of $J$ so that $W_i$ is homeomorphic to $X_i \times \R_+$.  Let $e_1, \hdots, e_n$ denote the ends of $M$ corresponding to $W_1, \hdots, W_n$ respectively. 

\begin{theorem}\label{end-stabilized}
Let $S$ be any Heegaard surface for $M$.  If $S \cap W_i$ is of infinite genus then $S$ is $e_i$-stabilized.  Furthermore, if $X_i$ is not a sphere $S \cap W_i$ is of infinite genus and, therefore, $S$ is $e_i$-stabilized. 
\end{theorem}

The promised classification is contained in the following propositions.  The proofs of these propositions use Theorem \ref{end-stabilized} to give information about stabilizations and then appeal to Frohman and Meeks' theorem for the existence of the desired isotopies. \newline

In Section \ref{examples}, it was explained how to obtain finite genus splittings of non-compact 3-manifolds: remove some finite number of closed balls from a compact 3-manifold.  All such 3-manifolds are deleted boundary 3-manifolds.  One consequence of Theorem \ref{end-stabilized} is that these are the only deleted boundary 3-manifolds with finite genus Heegaard splittings.   All others have only infinite genus splittings and we can classify them up to approximate isotopy and up to proper ambient isotopy. \newline

The following propositions provide the classification.  Recall that $M = \ob{M} - J$ is a deleted boundary 3-manifold:

\begin{proposition}\label{all spheres}
Suppose that $J$ consists of 2-spheres and that $M'$ is obtained from $\ob{M}$ by attaching 3-balls to $J$.  Then, up to proper ambient isotopy of $M$, any finite genus Heegaard surface in $M$ is the intersection of a Heegaard surface for $M'$ with $M$.  The Heegaard surface in $M'$ intersects each attached 3-ball in a properly embedded disc.   If two such splittings of $M'$ were isotopic then the resulting splittings of $M$ are properly ambient isotopic.
\end{proposition}

\begin{proposition}\label{approx isotopic}
Suppose that $S$ and $T$ are infinite genus Heegaard surfaces for $M$ whose splittings have the same partition of $\boundary M$.  Then $S$ and $T$ are approximately isotopic.
\end{proposition}

\begin{proposition}\label{classification}
Suppose that $S$ and $T$ are infinite genus Heegaard surfaces for $M$ with the same partition of $\boundary M$.  Consider the following condition:
\begin{enumerate}
\item[(*)] For each $i$, $S \cap W_i$ has infinite genus if and only if $T \cap W_i$ is of infinite genus.
\end{enumerate}
Then (*) holds if and only if $S$ and $T$ are properly ambient isotopic.  
\end{proposition}

\begin{proposition}\label{no spheres}
If no $X_i$ is a 2-sphere then any two Heegaard splittings of $M$ with the same partition of $\boundary M$ are equivalent up to proper ambient isotopy.
\end{proposition}

Before we prove the theorem and the classifications, we review a technique developed by Scharlemann and Thompson \cite{ScTh94a} which was inspired by work of Otal.  We also need to review the classification of Heegaard splittings of $G \times I$ where $G$ is a closed surface.

\subsection{Edge-Slides of Reduced Spines}\label{edge-slides of reduced spines}

\begin{definition}
Suppose that $Q$ is a compact 3-manifold and that $\Sigma$ is a finite graph in $Q$  such that $\Sigma$ intersects $\boundary Q$ in valence one vertices.  Let $B$ denote the components of $\boundary Q$ which intersect $\Sigma$.  If $\cl(Q - \eta(B \cup \Sigma))$ is a compressionbody then $\Sigma$ is a \defn{reduced spine}. 
\end{definition}

Choose an edge $e \subset \Sigma$ and a path $\gamma \subset \boundary Q \cup \Sigma$ with $\gamma$ beginning at an endpoint of $e$ but otherwise disjoint from $e$.  An edge-slide of $e$ over $\gamma$ replaces $e$ with the union of $e$ and a copy of $\interior(\gamma)$ pushed slightly away from $\Sigma \cup B$.  See \cite{SaScSch01,ScTh93,ScTh94a} for more detail.  Edge slides give isotopies of the surface $S = (B - \interior(\eta(\Sigma))) \cup \boundary \eta(\Sigma)$.  Conversely, an isotopy of a Heegaard surface can be converted into a sequence of edge-slides and isotopies of a reduced spine for one of the compressionbodies.  The correspondence between edge-slides of reduced spines and isotopies of the Heegaard surface will be useful for the proof of Theorem \ref{end-stabilized}.  The reason that this viewpoint is helpful is that if $Q$ is a compact submanifold of a non-compact manifold and if $(\boundary Q - \interior(\eta(\Sigma))) \cup \boundary \eta(\Sigma)$ is part of a Heegaard surface $S$ for $M$ then the isotopies described by edge-slides in $Q$ of $\Sigma$ are fixed off a regular neighborhood of $Q$ and so describe a proper isotopy of $S$. \newline

To increase the genus of the Heegaard surface obtained from the reduced spine, we may stabilize a reduced spine by choosing an edge $e \subset \Sigma$.  The edge $e$ is homeomorphic to $[0,1]$ and, choosing some homeomorphism, let $e'$ denote the subarc $[\frac{1}{4},\frac{3}{4}]$.  Introduce new vertices on $e$ at $\frac{1}{4}$ and $\frac{3}{4}$ and push the interior of $e'$ slightly off of $e$ to form a new edge $e''$ with endpoints on $e$ at the vertices $\frac{1}{4}$ and $\frac{3}{4}$.  The new edges $e''$ and $e'$ of $\Sigma$ bound a disc $D$ whose interior is disjoint from $\Sigma$.  The induced Heegaard splitting is stabilized in the usual sense as the boundary of the disc $D$ intersects a meridian disc of $\eta(\Sigma)$ exactly once. \newline

The final lemma of this section produces a reduced spine for $\text{(surface)} \times I$ with particular properties.  The spine gives rise to a relative version of a standard splitting of $\text{(surface)} \times I$.

\begin{lemma}\label{model surface}
Let $G$ be a closed surface of positive genus. Let $G'$ and $G''$ be the surfaces $G \times \{\frac{1}{4}\}$ and $G \times \{\frac{3}{4}\}$ in $G \times I$.  Let $n$ be an fixed integer bigger than or equal to twice the genus of $G$.  Let $P_0 = G \times [0,\frac{1}{4}]$. Then there is a connected reduced spine $\Sigma = \Sigma(G,n)$ in $G \times I$ such that $\Sigma$ intersects both boundary components of $G \times I$, $\Sigma$ intersects $P_0$ in a vertical arc, the rank of $H_1(\Sigma) = n$, and $\boundary \eta(\Sigma)$ is a relative Heegaard surface for $G \times [\frac{1}{4},1]$.
\end{lemma}

\begin{proof}
Consider $Q' = (G \times [\frac{7}{16},\frac{9}{16}]) - (\eta(* \times [\frac{7}{16},\frac{9}{16}]))$ where $*$ is a point on $G$.  Then $Q'$ is a handlebody of genus twice the genus of $G$.  Choose $\text{genus}(G)$ loops $L$ based at a point $b \in \interior Q'$ which represent generators of $\pi_1(Q',b)$.  Let $a$ be the arc $b \times I$ in $G \times I$ and assume, by general position, that the interior of each loop of $L$ is disjoint from $a$.  Since $\boundary Q'$ is a Heegaard surface for $G \times I$, $\boundary (Q' \cup \eta(a))$ is a relative Heegaard surface for $G \times I$.  Stabilize the reduced spine $\Sigma$ enough times so that the rank of its first homology is $n$.  Be sure that the stabilizations take place in the interval $[\frac{1}{4},1]$.  Then $a \cup L$ is a reduced spine for $G \times I$ satisfying the desired properties.
\end{proof}

\subsection{Heegaard Splittings of $\text{(closed surface)} \mathbf{\times I}$.}

Scharlemann and Thompson classified Heegaard splittings of $G \times I$, where $G$ is a closed connected surface.  In Theorem 6.1 of \cite{ScTh93} they give a way of interpreting their classification in terms of edge slides of spines (reduced or non-reduced).  The following are the versions of their results which we will need.

\begin{theorem}[Scharlemann-Thompson \cite{ScTh93}]\label{STclass1}
Suppose that $\Sigma$ and $\Psi$ are connected reduced spines for $G \times I$ which intersect both boundary components of $G \times I$ and whose first homology groups have the same rank.  Then there is a finite sequence of edge-slides and isotopies taking $\Sigma$ to $\Psi$.
\end{theorem}

\begin{theorem}[Scharlemann-Thompson \cite{ScTh93}]\label{STclass2}
If a Heegaard splitting of $G \times I$ has both boundary components of $G \times I$ contained in the same compressionbody and if the splitting surfaces has genus greater than twice the genus of $G$ then the splitting is stabilized.
\end{theorem}

\subsection{The proofs}
Before beginning each proof, the theorem or proposition has been repeated for the convenience of the reader.

\begin{thm6.4} 
If $S \cap W_i$ is of infinite genus then $S$ is $e_i$-stabilized.  Furthermore, if $X_i$ is not a sphere $S \cap W_i$ is of infinite genus and, therefore, $S$ is $e_i$-stabilized. 
\end{thm6.4}

\subsubsection*{Proof of Theorem \ref{end-stabilized}}

Since $M$ is end-irreducible (rel $C$) and $\boundary M \subset C$, Theorem \ref{well-placed} guarantees that there is an exhausting sequence $\{K_n\}$ which is well-placed on $S$.  In particular, $\fr K_n$ is incompressible in $M - C$ and no component of $\cl(M - K_n)$ is compact.  Recall that $W_i$ is a component of $\cl(M - C)$ and is homeomorphic to $X_i \times \R_+$ where $X_i$ is a closed connected surface.  For each $n$, the surface $\fr K_n \cap W_i$ is an incompressible surface in $W_i$.  Furthermore, as $H_2(W_i,\boundary W_i) = 0$ and $\cl(M - K_n)$ has no compact components, $\fr K_n \cap W_i$ is connected and is not a 2-sphere which is inessential in $W_i$.  

\begin{lemma}\label{h-cobord}
For each $i$ and for each $n$ the submanifold $\cl(K_{n+1} - K_n) \cap W_i$ is homeomorphic to $X_i \times I$.
\end{lemma}
\begin{proof}
The proof is well-known, but we include it for completeness.  Let $F = \fr K_{n+1} \cap W_i$.  $F$ is incompressible in $W_i$.  Let $N_n = \cl(K_{n+1} - K_n) \cap W_i$.  Suppose first that $X_i = S^2$.  In this case, $F$ is also homeomorphic to $S^2$.  As $F$ is essential it does not bound a ball in $W_i$.  By \cite[Theorem 3.1]{Br66}, $N_n$ is homeomorphic to $S^2 \times I$. \newline 

Now suppose that $X_i \neq S^2$.  As $W_i$ is irreducible, $F \neq S^2$.  The inclusion map of $F$ into $N_n$ induces an injective map on fundamental groups.  Since $W_i$ is homeomorphic to $X_i \times \R_+$, each loop in $N_n$ with basepoint on $F$ is homotopic (rel basepoint) to a loop outside of $N_n$.  Hence, each loop is homotopic into $F$.  Thus the inclusion of $F$ into $N_n$ induces an isomorphism of fundamental groups and, so by the h-cobordism theorem \cite[Theorem 10.2]{He04}, $N_n$ is homeomorphic to $F \times I$.  A similar argument shows that the submanifold bounded by $X_i$ and $F$ is homeomorphic to $F \times I$ and so $F$ is homeomorphic to $X_i$.
\end{proof}

Fix some $i$.  Let $W = \cl(W_i - K_2)$.  We will show that there is a subsequence of $\{K_n\}$ and a proper ambient isotopy of $S$ which is fixed off $\cl(W_i - K_1)$ so that either $W \cap \cl(K_{n+1} - K_n)$ is homeomorphic to $S^2 \times I$ and $S \cap W \cap \cl(K_{n+1} - K_n)$ is a genus 0 relative Heegaard surface or $S \cap W \cap \cl(K_{n+1} - K_n)$ is a stabilized relative Heegaard surface of $W \cap \cl(K_{n+1} - K_n)$. \newline

We deal first with the case when $X_i = S^2$.  Let $N_n = W \cap \cl(K_{n+1} - K_n)$ for each $n \geq 2$. 

\begin{lemma}\label{spheres inherit}
If $X_i = S^2$ then $S \cap N_n$ is a relative Heegaard surface for $N_n$.
\end{lemma}

\begin{proof}
Recall that for each $n$, $\fr K_n \cap W$ is an essential 2-sphere and, by property (WP1) of well-placed exhausting sequences, $V \cap (\fr K_n \cap W)$ is a single disc.  This implies that $U \cap (\fr K_n \cap W)$ is a single disc. Thus, for each $n \geq 2$, $U \cap N_n$ is a relative compressionbody with preferred surface $S \cap N_n$.  Similarly, for each $n \geq 2$, $V \cap N_n$ is a relative compressionbody with preferred surface $S \cap N_n$.  Thus $S \cap N_n$ is a relative Heegaard surface for $N_n$.
\end{proof}

By Lemma \ref{h-cobord}, $N_n$ is homeomorphic to $S^2 \times I$.  By the classification of Heegaard splittings of $S^2 \times I$, if $S \cap N_n$ has positive genus, there is a reducing ball for $S \cap N_n$ which is contained in $N_n$.  If $S \cap W_i$ is of infinite genus, there are infinitely many $n$ so that $S \cap N_n$ is of positive genus, and hence $S$ is $e_i$-stabilized. If $S \cap W_i$ is of finite genus, we can take a subsequence of $\{K_i\}$ so that $S \cap N_n$ has genus 0.  This concludes the case when $X_i = S^2$. \newline

Suppose, for the remainder, that $X_i$ is a closed orientable surface of positive genus $g$.  We do not begin by supposing that $S \cap W_i$ is of infinite genus but, rather, draw that as our first conclusion. \newline

Recall that since $\{K_n\}$ is well-placed on $S$, $V$ intersects each $\fr K_n \cap W$ is a single disc.  Let $N_n = W \cap \cl(K_{n+1} - K_n)$ for each $n \geq 1$.  Since $\{K_n\}$ is well-placed on $S$ the sequence $\{K_n\}$ has the outer collar property with respect to $U$.  This means that in each $U \cap N_n$ there is a collection of discs $\delta_n$ with boundary on $S \cap N_n$ so that $\sigma(U \cap N_n ; \delta_n)$ has a component which is $(\fr K_{n+1} \cap U \cap N_n) \times I$.  The frontier of $K_{n+1} \cap U \cap N_n$ is $(\fr K_{n+1} \cap U \cap N_n) \times \{0\}$.  On the other hand, $(\fr K_{n+1} \cap U \cap N_n) \times \{1\}$ is a subsurface of $S$ except at the remnants of $\delta_n$.  Since $V \cap N_n \cap \fr K_{n+1}$ is a single disc and since $\fr K_{n+1} \cap N_n$ is homeomorphic to $X_i$, the surface $\fr K_{n+1} \cap N_n \cap U$ is homeomorphic to $X_i$ with a single puncture.  As $X_i$ has positive genus $g$, the surface $\sigma(S \cap N_n ; \delta_n)$ has positive genus, and, therefore, $S \cap N_n$ has positive genus for all $n \geq 1$.  This implies that $S \cap W$ has infinite genus. \newline

Take a subsequence of $\{K_n\}$ so that the first two terms of the new exhausting sequence are still $K_1$ and $K_2$ but so that the genus of $S \cap \cl(K_{n+1} - K_n) \cap W$ is at least $3g$ for $n \geq 1$.  We continue referring to $\cl(K_{n+1} - K_n) \cap W$ as $N_n$. \newline

Fix some $n \geq 2$ and let $N = N_n$.  By Lemma \ref{h-cobord}, $N$ is homeomorphic to $X_i \times I$.  Let $F_0 = \fr K_n \cap N$ and $F_1 = \fr K_{n+1} \cap N$.  $V$ intersects $F_i$ in a single disc $D_i$ for $i \in \{0,1\}$. Since $\{K_i\}$ has the outer collar property with respect to $U$, there is a collection of boundary-reducing discs $\delta_0$ for $U \cap K_n \cap W$ with boundary on $S$ and such that $\sigma(U \cap K_n \cap W;\delta_0)$ contains a component $P^U_0$ with boundary containing $F_0 \cap U$ and which is homeomorphic to $(F_0 \cap U) \times I$.  Since $S$ is the preferred surface of $U \cap K_n$, there is a copy of $D^2 \times I$ embedded in $V$ so that $D^2 \times \{0\} = V \cap F_0$ and $\boundary D^2 \times I = S \cap P^U_0$.  Let $P_0$ be the union of $P^U_0$ and this $D^2 \times I$.  Note that $P_0$ is homemorphic to $F_0 \times I$, has $F_0$ as a boundary component, and has $V$ running through $P_0$ as the neighborhood of an arc which is vertical in the product structure.  Let $F'_0 = \boundary P_0 - F_0$. \newline

\begin{figure}[ht]
\scalebox{0.45}{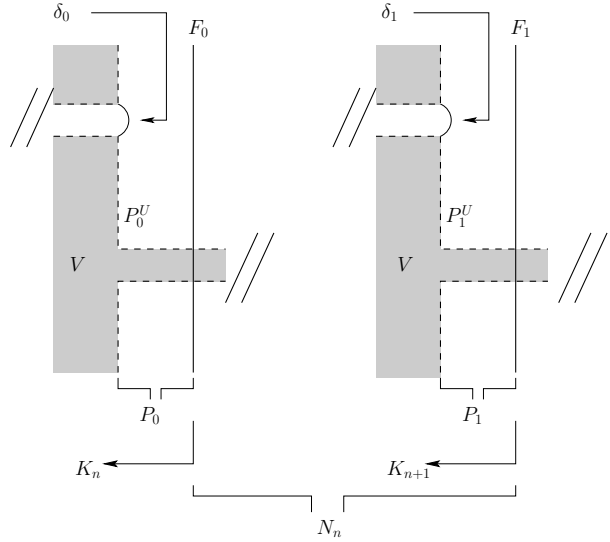}
\caption{A schematic representing $N$.}
\label{spine}
\end{figure}

We can perform a similar construction on $K_{n+1}$ to obtain, embedded in $N_n$, a submanifold $P_1$ homeomorphic to $F_1 \times I$, with $\boundary P_1 = F_1 \cup F'_1$ and $V \cap P_1$ a neighborhood of a vertical arc. Let $N' = N \cup P_0$ and $N'' = \cl(N' - P_1)$.  Note that $N'$ and $N''$ are homeomorphic to $X_i \times I$, since $F_0,F_1,F'_0,$ and $F'_1$ are all homeomorphic to $X_i$.  See Figure \ref{spine}.\newline

Let $\Sigma_V$ be a spine for $V$ in $M$ which intersects each surface $F'_0,F_0,F'_1,F_1$ exactly once and which is a vertical arc in $P_0$ and $P_1$.  Let $\Sigma_S = \Sigma_V \cap N'$.  Note that $\Sigma_S$ is a reduced spine for a Heegaard splitting of $N''$.  To see this, recall that $U \cap (K_{n+1} - K_n)$ is a handlebody (Corollary \ref{complementary compressionbodies}) and notice that $N'' - \eta(\Sigma_S \cup \boundary N')$ is homeomorphic to $U \cap (K_{n+1} - K_n)$. We wish to show that after a proper ambient isotopy of $S$ which is the identity off of $\eta(N')$, $S\cap N$ is a Heegaard surface for $N$.  \newline

Choose a connected reduced spine $\Sigma_T$ for a Heegaard splitting of $N''$ such that $\Sigma_T$ intersects $P_0$ in a vertical arc, the rank of $H_1(\Sigma_T)$ is the same as the rank of $H_1(\Sigma_S)$, $\Sigma_T \cap F'_1 \neq \nil$, $\Sigma_T \cap F'_0 \neq \nil$, and $\boundary \eta(\Sigma_T)$ is a hollow Heegaard surface for $N$.  Such a spine exists by Lemma \ref{model surface}. We call $\Sigma_T$ the \defn{model spine}. \newline

By the Scharlemann-Thompson classification of Heegaard splittings of $\text{(surface)} \times I$ (Theorem \ref{STclass1}) since $\Sigma_S$ and $\Sigma_T$ are both reduced spines with first homologies of the same rank and since they have the same partition of $\boundary N''$ there is a sequence of edge-slides and isotopies which takes $\Sigma_S$ to $\Sigma_T$.  It is easy to arrange these slides to be away from $\delta_0 \cup \delta_1$.  The sequence of edge slides thus describes an isotopy of the surface $S \cap \eta(N'')$.  By the choice $\Sigma_T$, we have that after the isotopy, $S \cap N$ is a relative Heegaard surface of genus at least $3g$ for $N'$. \newline

The next corollary follows from our work so far; it is a technical result which will be useful for the classifications.

\begin{corollary}\label{induced rel HS}
If $X_i$ is a closed surface of positive genus then after a proper ambient isotopy of $S$ which is supported on a neighborhood of $W'_i = \cl(M - K_1) \cap W_i$ we have that $S \cap K_1$ is the same before and after the isotopy and afterwards $S \cap W'_i$ is a relative Heegaard surface for $W'_i$.
\end{corollary}

\begin{proof}
Perform the isotopy just described so that $N_2 = \cl(K_2 - K_1)$ inherits a relative Heegaard splitting from $S$.  This isotopy is fixed off a neighborhood of $W'_i= \cl(M - K_1) \cap W_i$ and $S \cap K_1$ is the same before and after the isotopy.  Since $N_2 \subset W'_i$, there are now discs in $U \cap W'_i$ with boundary on $S$ so that boundary reducing $U \cap W'_i$ along those discs leaves a component homeomorphic to $(\fr W'_i \cap U) \times I$.  Since $V \cap W'_i$ is a disc we have that $U \cap W'_i$ and $V \cap W'_i$ are relative compressionbodies with preferred surface $S \cap W'_i$.  Thus, $S \cap W'_i$ is a relative Heegaard surface for $W'_i$. 
\end{proof}

We now continue the proof of Theorem \ref{end-stabilized}.  For each even $n$, perform this ambient isotopy on $N_n$.  By construction, the union of these ambient isotopies is a proper ambient isotopy of $S \cap W_i$.  After the isotopy, for each even $n$, $S \cap N_n$ is a relative Heegaard surface of genus at least $3g$ for a space homeomorphic to $X_i \times I$ where the genus of $X_i$ is $g$.  By the Scharlemann-Thompson classification of splittings of $\text{(surface)} \times I$ (Theorem \ref{STclass2}), there is a reducing ball for $S$ in each $N_n$ for $n$ even.  Hence $S$ is $e_i$-stabilized.  This concludes the proof of Theorem \ref{end-stabilized}
\qed

\subsubsection*{Proof of Classification}

\begin{prop6.5}
Suppose that $J$ consists of 2-spheres and that $M'$ is obtained from $\ob{M}$ by attaching 3-balls to $J$.  Then, up to proper ambient isotopy of $M$, any finite genus Heegaard surface in $M$ is the intersection of a Heegaard surface for $M'$ with $M$.  The Heegaard surface in $M'$ intersects each attached 3-ball in a properly embedded disc.   If two such splittings of $M'$ were isotopic then the resulting splittings of $M$ are properly ambient isotopic.
\end{prop6.5}

\begin{proof}[Proof of Proposition \ref{all spheres}]
Suppose that $U_S \cup_S V_S$ and $U_T \cup_T V_T$ are both finite genus Heegaard splittings of $M$.  Let $M'$ be the compact 3-manifold obtained from $\ob{M}$ by removing only the interiors of the 3-balls whose removal created $M$. \newline

There is an exhausting sequence $\{K_i\}$ for $M$ which is well-placed on $S$ and an exhausting sequence $\{L_i\}$ which is well-placed on $T$.  We may assume that $K_1$ and $L_1$ are homeomorphic to $M'$ and that $S \cap (M - K_1)$ and $T \cap (M - L_1)$ have genus zero.  The frontiers of the exhausting elements are essential spheres in $S^2 \times \R_+$ so, after taking a subsequence of each, there is a proper ambient isotopy of $M$ which takes $\fr L_i$ to $\fr K_i$ for each $i$ and so that $V_S \cap \fr K_i$ equals $V_S \cap \fr L_i$. \newline

Let $N_n = \cl(K_{n+1} - K_n)$.  By Lemma \ref{spheres inherit}, each component of $N_n$ inherits a genus zero relative Heegaard splitting from $S$ and also from $T$.  By Waldhausen's classification of splittings of $S^2 \times I$, there is a proper ambient isotopy taking $S \cap N_n$ to $T \cap N_n$ which is fixed on $\fr N_n$.  The union of these isotopies over all $n$ is a proper ambient isotopy of $M$ taking $S \cap \cl(M - K_1)$ to $T \cap \cl(M - K_1)$.  In particular, we may assume that $S \cap \cl(M - K_1)$ and $T \cap \cl(M - K_1)$ are vertical annuli in $S^2 \times \R_+$. \newline

When we compactify $M$ to $\ob{M}$, $S \cap \cl(M - K_1)$ and $T \cap \cl(M - K_1)$ compactify to compact annuli.  $V_S \cap \cl(M - K_1) = V_T \cap \cl(M - K_1)$ compactifies to $V' = D^2 \times I$.  Let $V'_S$ and $V'_T$ be the compactified versions of $V_S$ and $V_T$ respectively.  Attach the 3-balls to $\ob{M}$ to create $M'$ and let $U'_S = \cl(M' - V'_S)$ and $U'_T = \cl(M' - V'_T)$. It is clear from the construction that $U'_S$, $U'_T$, $V'_S$ and $V'_T$ are absolute compressionbodies and that the splittings $U_S \cup_S V_S$ and $U_T \cup_T V_T$ are obtained from splittings of $M'$ in the correct fashion. \newline

Furthermore, if we remove the open 3-balls from $\ob{M}$ to create $M'$ we can extend the splittings $U'_S \cup_{S} V'_S$ and $U'_T \cup_T V'_T$ of $\ob{M}$ to be relative Heegaard splittings of $M'$.  By the Marionette Lemma the relative Heegaard splittings of $M'$ are isotopic if and only if the absolute splittings of $\ob{M}$ are isotopic.  If the splittings of $M'$ are isotopic then since $K_1$ is homeomorphic to $M'$, the surfaces $S \cap K_1$ and $T \cap K_1$ are isotopic.  Thus, since we already have $S \cap \cl(M - K_1) = T \cap \cl(M - K_1)$ we can arrange by a proper ambient isotopy for $S$ to be equal to $T$.
\end{proof}

\begin{prop6.6}
Suppose that $S$ and $T$ are infinite genus Heegaard surfaces for $M$ whose splittings have the same partition of $\boundary M$.  Then $S$ and $T$ are approximately isotopic.
\end{prop6.6}

\begin{proof}[Proof of Proposition \ref{approx isotopic}]
If $S$ and $T$ have infinite genus then $S \cap W_i$ and $T \cap W_j$ have infinite genus for some $i,j$.  Since each $W_k$ is homeomorphic to $X_k \times I$ where $X_k$ is a closed surface, Theorem \ref{end-stabilized} shows that $S$ must be $e_i$-stabilized and $T$ must be $e_j$-stabilized.  Theorem \ref{Frohman-Meeks} then shows that $S$ and $T$ are approximately isotopic.
\end{proof}

\begin{prop6.7}
Suppose that $S$ and $T$ are infinite genus Heegaard surfaces for $M$ with the same partition of $\boundary M$.  Consider the following condition:
\begin{enumerate}
\item[(*)] For each $i$, $S \cap W_i$ has infinite genus if and only if $T \cap W_i$ is of infinite genus.
\end{enumerate}
Then (*) holds if and only if $S$ and $T$ are properly ambient isotopic.
\end{prop6.7}

\begin{proof}[Proof of Proposition \ref{classification}]
The proof of Lemma \ref{proper ambient isotopy invariant} can be adapted to show that if $S$ and $T$ are properly ambient isotopic then (*) holds. \newline

Suppose, then, that $S$ and $T$ satisfy (*).  We desire to show that $S$ and $T$ are properly ambient isotopic. Using Proposition \ref{approx isotopic}, we will be able to enlarge $C$ to a compact set $C'$ such that (after performing proper ambient isotopies of $S$ and $T$) $C'$ has the following properties:
\begin{enumerate}
\item $\cl(M - C')$ is homeomorphic to $\cup X_i \times \R_+$ 
\item $S \cap C' = T \cap C'$
\item $V_S \cap \fr C' = V_T \cap \fr C'$ and each of these consists of a single disc on each component of $\fr C'$.
\item For each $W'_i = \cl(M - C) \cap W_i$ where $S$ and $T$ are of infinite genus, the surfaces $S\cap W'_i$ and $T\cap W'_i$ are relative Heegaard surfaces for $W'_i$.
\item For each $W'_i$ where $S$ and $T$ are not of infinite genus, the surfaces $S \cap W'_i$ and $T \cap W'_i$ are of genus zero.
\end{enumerate}

The way to achieve this is to take an exhausting sequence $\{K_i\}$ for $M$ which is well-placed on $S$ such that in each component of $\cl(M - K_1)$ $S$ and $T$ are both of either infinite genus or of genus zero.  Then use the fact that $S$ and $T$ are approximately isotopic to isotope them so that $S \cap K_1 = T \cap K_1$.  Let $C' = K_1$.  If a certain $X_i$ is not a 2-sphere, Corollary \ref{induced rel HS} guarantees that a further proper ambient isotopy of $S$ and $T$ can be performed which is supported on a neighborhood of $W'_i = \cl(M - C') \cap W_i$ so that after the isotopy $S \cap C'$ still equals $T \cap C'$ but we now have property (4) in addition to property (3) for that $W'_i$. In the case when $X_i = S^2$, $S$ and $T$ automatically give relative Heegaard splittings of $W'_i$ as $\boundary_-(U \cap W'_i)$ and $\boundary_-(V \cap W'_i)$ can be taken to be the discs $U \cap X_i$ and $V \cap X_i$ respectively. \newline

For each $W'_i$ in which $S$ and $T$ are of infinite genus, Theorem \ref{end-stabilized} guarantees $S \cap W'_i$ and $T \cap W'_i$ are infinitely stabilized.  Since $W'_i$ is 1-ended, Theorem \ref{Frohman-Meeks} guarantees that the absolute Heegaard splittings of $W'_i$ induced by $S \cap W'_i$ and $T \cap W'_i$ are equivalent by a proper ambient isotopy in $W'_i$.  By the Marionette Lemma, $S \cap W'_i$ and $T \cap W'_i$ are properly ambient isotopic within $W'_i$.  For each $W'_i$ where $S$ and $T$ are of genus zero, the fact that $S$ and $T$ are properly ambient isotopic in $W'_i$ follows from Proposition \ref{all spheres}. \newline

Since in each component of $\cl(M - C')$ there is a proper ambient isotopy of $S$ and $T$ in that component so that they coincide, and since $S$ and $T$ already coincide in $C'$ there is a proper ambient isotopy of $M$ taking $T$ to  $S$.
\end{proof}

\begin{prop6.8}
If no $X_i$ is a 2-sphere then any two Heegaard splittings of $M$ with the same partition of $\boundary M$ are equivalent up to proper ambient isotopy.
\end{prop6.8}

\begin{proof}[Proof of Proposition \ref{no spheres}]
By Theorem \ref{end-stabilized}, $S$ and $T$ are end-stabilized.  Theorem \ref{Frohman-Meeks} then implies that they are properly ambient isotopic.
\end{proof}
\appendix \section{Infinitely Stabilized Heegaard Splittings}
The goal of this section is to give a detailed proof the following theorem which is due, essentially, to Frohman and Meeks.  Our methods are the same but we elaborate in order to fix the error mentioned previously.  We refer the reader to earlier sections for the definitions of the terms used here.  

\begin{theorem}\label{Frohman-Meeks 2}
Let $M$ be a non-compact orientable 3-manifold with compact boundary not containing any 2-sphere components.  Suppose that $M = U_S \cup_S V_S$ and $M = U_T \cup_T V_T$ are two Heegaard splittings of $M$ with the same partition of $\boundary M$.  If both $S$ and $T$ are infinitely stablized then they are approximately isotopic.  If both $S$ and $T$ are end-stabilized then they are properly ambient isotopic.
\end{theorem}

In \cite{FrMe97}, Frohman and Meeks introduce a technique which they call ``stealing handles from infinity".  This method provides a proper isotopy of an infinitely stabilized splitting so that for any compact submanifold $K$, $S \cap K$ is stabilized an arbitrary number of times.

\begin{proposition}[{Frohman-Meeks \cite[Proposition 2.1]{FrMe97}}] \label{stealing handles}
Suppose that $M = U \cup_S V$ is an infinitely stabilized Heegaard splitting of $M$.  Let $C$ be a submanifold of $M$ which is adapted to $S$.  Then for any given $n \in \N$ there is a proper ambient isotopy of $S$ so that $S \cap C$ has been stabilized at least $n$ times.
\end{proposition}

\begin{proof}[Sketch of Proof] Since $S$ is infinitely stabilized, we can find $n$ disjoint reducing balls for $S$ in the complement of $C$.  We may then use paths in the surface $S$ to isotope these balls along $S$ into $C$.
\end{proof}

\begin{definition}
An exhausting sequence $\{K_i\}$ is \defn{perfectly adapted} to $S$ if it is adapted to $S$ and, additionally, each $\cl(K_{i+1} - K_i)$ is adapted to $S$.  (See Section \ref{Types of Exh. Seq.}.)  Note that a subsequence of a perfectly adapted sequence is perfectly adapted.
\end{definition}

A useful corollary of Proposition \ref{stealing handles} is:

\begin{corollary}\label{key cor}
Suppose that $U_S \cup_S V_S$ and $U_T \cup_T V_T$ are two end-stabilized splittings of $M$ with the same partition of $\boundary M$.  If there is an exhausting sequence $\{K_i\}$ for $M$ with the following properties:
\begin{enumerate}
\item[(i)] $\boundary M \subset K_1$
\item[(ii)] $V_S \cap \fr K_i$ and $V_T \cap \fr K_i$ consist of discs for all $i$.
\item[(iii)]$V_S \cap \fr K_i = V_T \cap \fr K_i$ for all $i$.
\item[(iv)] $\{K_i\}$ is perfectly adapted to both $S$ and $T$.
\end{enumerate}
then $S$ and $T$ are equivalent up to proper ambient isotopy.
\end{corollary}

\begin{proof}
By the Reidemeister-Singer theorem and the Marionette Lemma, after finitely many stabilizations of $S \cap K_1$ and $T \cap K_1$ there is an ambient isotopy of $K_1$ so that $S \cap K_1 = T \cap K_1$.  Since both $S$ and $T$ are end-stabilized, these stabilizations can be achieved by stealing handles from infinity.  Thus, we may assume that $S \cap K_1 = T \cap K_1$.  By the assumption that $\{K_i\}$ is perfectly adapted to both $S$ and $T$, the intersections of $U_S \cup_S V_S$ and $U_T \cup_T V_T$ with any compact component $L$ of $\cl(M - K_1)$ give a relative Heegaard splittings of $L$.  By stealing more handles from infinity and passing them through $K_1$ we may stabilize $S \cap L$ and $T \cap L$ enough times so that after performing an ambient isotopy of $L$, $S$ and $T$ coincide in $K_1 \cup L$.  We may do this for each compact component of $\cl(M - K_1)$.  Since there are only finitely many such components, we have constructed proper ambient isotopies of $S$ and $T$ so that they coincide on $K_1$ and each compact component of $\cl(M - K_1)$.  We proceed by induction. \newline

Suppose that we have performed proper ambient isotopies of $M$ so that $S \cap K_{n-1} = T \cap K_{n-1}$ and $S$ and $T$ coincide on each compact component of $\cl(M - K_{n-1})$.  We will show that there are proper ambient isotopies of $S$ and $T$ which are fixed on $K_{n-1}$ so that after the isotopies $S$ and $T$ coincide on $K_n$ and each compact component of $\cl(M - K_n)$.  This will show that the composition of the isotopies of $S$ converges to a proper ambient isotopy of $S$ and the composition of the isotopies of $T$ converges to a proper ambient isotopy of $T$.  Thus, we will have shown that there are proper ambient isotopies of $S$ and $T$ which make them coincide with a third Heegaard surface for $M$.  Hence, $S$ and $T$ are properly ambient isotopic. \newline

Let $L$ be a component of $\cl(K_n - K_{n-1})$.  By hypothesis, both $S$ and $T$ are relative Heegaard surfaces for $L$.  If every non-compact component of $\cl(M - L)$ contains $K_{n-1}$ then $L$ is contained in a compact component of $\cl(M - K_{n-1})$ and so $S \cap L = T \cap L$.  \newline

We may, thus, suppose that there is a non-compact component of $\cl(M - L)$ which does not contain $K_{n-1}$. The surfaces $S$ and $T$ are both end-stabilized and so we may steal handles from that non-compact component of $\cl(M - L)$ in order to stabilize $S \cap L$ and $T \cap L$ enough times so that they are ambient isotopic in $L$.  Since, $S$ and $T$ already coincide on $\fr K_{n-1}$ we may take the ambient isotopy to be the identity on $\fr K_{n-1} \cap L$.  Thus, there is a proper ambient isotopy of $S$ and a proper ambient isotopy of $T$, each fixed on $K_{n-1}$ so that after the isotopies $S \cap K_n = T \cap K_n$. \newline

Now suppose that $L'$ is a compact component of $\cl(M - K_n)$.  As before, $S$ and $T$ both give relative Heegaard splittings of $L'$.  If $S \cap L' \neq T \cap L'$ then $L'$ is not contained in a compact component of $\cl(M - K_{n-1})$.  As in each component of $\cl(K_n - K_{n-1})$ $S$ and $T$ are connected surfaces, this implies that there are paths in $S$ and $T$ from a non-compact component of $\cl(M - K_n)$ to $L'$ which do not intersect $K_{n-1}$.  Thus, we may stabilize $S \cap L'$ and $T \cap L'$ as much as we wish by stealing handles from infinity via paths that do not intersect $K_{n-1}$.  Now isotope in $L'$ so that the splittings coincide.  We have, therefore, constructed proper ambient isotopies of $S$ and $T$ which are fixed on $K_{n-1}$ such that after performing the isotopies $S \cap K_n$ equals $T \cap K_n$ and $S$ and $T$ also coincide on each compact component of $\cl(M - K_n)$. Thus, $S$ and $T$ are properly ambient isotopic in $M$.  
\end{proof}

To show that two end-stabilized splittings of $M$ with the same partition of $\boundary M$ are properly ambient isotopic, we will show that there is an exhausting sequence for $M$ satisfying the requirements of Corollary \ref{key cor}.  The first task is to show that if $S$ and $T$ have perfectly adapted exhausting sequences then there is a perfectly adapted sequence of $M$ adapted to both $S$ and $T$ simultaneously.

\begin{lemma}[{Frohman-Meeks \cite[Proposition 2.3]{FrMe97}}]\label{mutually adapted}
Suppose that $K_1$ and $K_2$ are two submanifolds of $M$ such that $K_1, K_2,$ and $\cl(K_2 - K_1)$ are adapted to $S$.  Suppose that $L_1$ and $L_2$ are two submanifolds of $M$ such that $L_1, L_2$ and $\cl(L_2 - L_1)$ are adapted to $T$.  Assume also that $K_1 \subset L_1 \subset K_2 \subset L_2$ where each inclusion is into the interior of the succeeding submanifold. \newline

Then after stabilizing and isotoping $S$ in $\cl(K_2 - K_1)$ and stabilizing and isotoping $T$ in $\cl(L_2 - L_1)$ there is a submanifold $J_1$ of $M$ adapted to both $S$ and $T$ so that $V_S \cap \fr J_1$ equals $V_T \cap \fr J_1$ and these intersections consist of discs.  
\end{lemma}

\begin{proof}
Push the frontier of $K_2$ slightly into $K_2$ to form a surface $F \subset K_2$.  Let $M_1$ be the submanifold bounded by $\fr K_2$ and $F$.  ($M_1$ is, of course, homeomorphic to $\fr K_2 \times I$.)  Let $M_2$ be the submanifold bounded by $F$ and $\fr K_1$.  Let $N_1$ be the submanifold with boundary $\fr L_2 \cup F$ and let $N_2$ be the submanifold with boundary $F \cup \fr L_1$.  Let $J_1 = K_1 \cup M_2$.  Take Heegaard splittings of $M_1, M_2, N_1$ and $N_2$ with Heegaard surfaces $S_1, S_2, T_1$ and $T_2$ respectively. We should choose these splittings so that all the boundary components of each submanifold are contained in the same compressionbody of the splitting.\newline

We can use the Heegaard surfaces $S_1$ and $S_2$ to form a Heegaard surface $\ob{S}$ for $\cl(K_2 - K_1)$.  To do this, note that there are surfaces $S'_1$ and $S'_2$ in $M_1$ and $M_2$ which are subsurfaces of $S_1$ and $S_2$ except at a finite number of open discs which are parallel to $F = M_1 \cap M_2$.  The surfaces $S'_1$ and $S'_2$ cobound a product region $S'_2 \times I$.  The surface $F$ may be assumed to be $S'_2 \times \{\frac{1}{2}\}$.  Take a disc $D \subset S'_2 \cap S_2$ so that in the product region $S'_2 \times I$ the tube $D \times I$ is disjoint from $\cl(S'_1 - S_1)$.  The Heegard surface $\ob{S}$ for $\cl(K_2 - K_1)$ is formed by taking $(S_1 \cup S_2 \cup D \times I) - \interior (D \times I)$.  We say that $\ob{S}$ is formed by \defn{tubing together} $S_1$ and $S_2$.  This process is different from the amalgamation of Heegaard splittings.  Similarly, we may form a Heegaard surface $\ob{T}$ for $\cl(L_2 - L_1)$ by tubing together $T_1$ and $T_2$.  Since in both constructions the tube intersects $F$ in a single disc, we may arrange that $\ob{S} \cap F = \ob{T} \cap F$ and that these intersections are a single inessential loop on $F$.  Finally, using the product region in the compressionbodies containing $\fr(K_2 - K_1)$ we may use vertical tubes to extend $\ob{S}$ to be a relative Heegaard splitting for $\cl(K_2 - K_1)$ which coincides with $S$ on $\fr(K_2 - K_1)$.  Similarly, extend $\ob{T}$ to be a relative Heegaard splitting for $\cl(L_2 - L_1)$ which coincides with $T$ on $\fr(L_2 - L_1)$.  We call the Heegaard splittings given by $\ob{S}$ and $\ob{T}$ the \defn{model splittings}. \newline

The Reidemeister-Singer theorem and the Marionette Lemma imply that by stabilizing $S$ and $\ob{S}$ enough in $\cl(K_2 - K_1)$ we may perform an ambient isotopy of $\cl(K_2 - K_1)$ which brings $S \cap \cl(K_2 - K_1)$ to $\ob{S}$.  Similarly, we may stabilize $T \cap \cl(L_2 - L_1)$ and $\ob{T}$ enough times so that there is an ambient isotopy of $\cl(L_2 - L_1)$ which brings $T \cap \cl(L_2 - L-1)$ to $\ob{T}$.  Since $\ob{S}$ and $\ob{T}$ coincide on $\fr J_1 = F$ we have now arranged that $J_1$ is a submanifold adapted to both $S$ and $T$ and that $S \cap \fr J_1 = T \cap \fr J_1$ and these intersections consists of a single inessential loop on each component of $\fr J_1$.
\end{proof}

\begin{corollary}\label{mutually adapted exh seq}
Suppose that $\{K_i\}$ is an exhausting sequence perfectly adapted to $S$ and that $\{L_i\}$ is an exhausting sequence perfectly adapted to $T$.  Assume that, for all $i$, $K_i \subset L_i \subset K_{i+1}$.  Then after stabilizing $S$ and $T$ in each component of $\cl(K_{i+1} - K_i)$ and $\cl(L_{i+1} - L_i)$ respectively we may properly isotope $S$ and $T$ so that there is an exhausting sequence $\{J_i\}$ which is perfectly adapted to both $S$ and $T$ and is such that $S \cap \fr J_i = T \cap \fr J_i$ and the intersection consists of a single inessential loop on each component of $\fr J_i$.
\end{corollary}

\begin{proof}
Construct $J_1$ as in the proposition.  Assuming that we have constructed $J_{n-1}$ we will demonstrate how to construct $J_n$.  Build $J_n$ as in the proposition, letting $K_{n+1}$, $K_n$, $L_{n+1}$, $L_n$ play the roles of $K_2, K_1, L_2$ and $L_1$.  Choose model splittings for each component of $\cl(K_{n+1} - K_{n})$ and $\cl(L_{n+1} - L_n)$ which coincide with the model splittings of $\cl(K_n - K_{n-1})$ and $\cl(L_n - L_{n-1})$ on $\fr K_n$ and $\fr L_n$ respectively.  Stabilize the model splittings enough times so that after stabilizing $S \cap \cl(K_{n+1} - K_n)$ and $T \cap \cl(L_{n+1} - L_n)$ we may perform ambient isotopies of $S \cap \cl(K_{n+1} - K_n)$ and $T \cap \cl(L_{n+1} - L_n)$ so that they coincide with the model splittings.  These isotopies are supported off $K_{n-1}$ and $L_{n-1}$ respectively. Note that, by the construction of the model splittings, $\cl(J_n - J_{n-1})$ is adapted to both $S$ and $T$ (after performing the isotopies). \newline

We thus obtain an exhausting sequence $\{J_i\}$ for $M$.  The final remarks of the previous paragraph show that there are proper ambient isotopies of $S$ and $T$ so that $\{J_i\}$ is perfectly adapted to both Heegaard surfaces.
\end{proof}

\begin{remark}
So far we have shown that if $S$ and $T$ are end-stabilized splittings and if there are exhausting sequences perfectly adapted to each of them then (after stealing handles from infinity and performing other proper ambient isotopies of $S$ and $T$) there is an exhausting sequence which is perfectly adapted to both of them at the same time and furthermore $S$ and $T$ coincide on the frontiers of the exhausting submanifolds.  Corollary \ref{key cor} then shows that $S$ and $T$ are properly ambient isotopic.  It thus remains to show that an end-stabilized splitting has a perfectly adapted exhausting sequence which is adapted to it.  The following lemmas show how we can achieve this.  This lemma fixes the misstatement in \cite[Proposition 2.2]{FrMe97} mentioned in the introduction.
\end{remark}

\begin{lemma}\label{edge slides again}
Let $M = U \cup_S V$ be an absolute Heegaard splitting of the non-compact 3-manifold $M$ and let $\{K_i\}$ be an exhausting sequence for $M$ adapted to $S$.  Assume that, for each $i$, $V \cap \fr K_i$ consists of discs and that the sequence $\{K_i\}$ has the outer collar property with respect to $U$.  Then after stabilizing $S \cap \cl(K_n - K_{n-1})$, for each $n \geq 3$, a finite number of times, there is a proper ambient isotopy of $S \cap K_n$ with the following properties:
\begin{enumerate}
\item[(i)] The isotopy is fixed on $K_{n-2} \cup \cl(M - K_n)$.
\item[(ii)] $S \cap K_{n-1}$ is the same before and after the isotopy.
\item[(iii)] After the isotopy, $S$ is a relative Heegaard surface for $\cl(K_n - K_{n-1})$.
\end{enumerate}
\end{lemma}

The proof is similar to the proof of Proposition \ref{end-stabilized}.  The reader is referred to Section \ref{edge-slides of reduced spines} for the definitions and properties of edge-slides.

\begin{proof}
Let $N$ be a component of $\cl(K_n - K_{n-1})$.  Let $F_2 = \fr K_n \cap N$ and $F_1 = \fr K_{n-1} \cap N$.  Since $\{K_i\}$ has the outer collar property, there are discs $\delta_1 \subset (U \cap K_{n-1})$ with boundary on $S$ so that $\sigma(U \cap K_{n-1};\delta_1)$ contains a product region $P^U_1 = (F_1 \cap U) \times I \subset U \cap \cl(K_{n-1} - K_{n-2})$ with $F_1 \cap U = (F_1 \cap U) \times \{0\}$.  Let $(F'_1 \cap U)$ signify $(F_1 \cap U) \times \{1\}$; it is a subsurface of $S$ except at the remnants of the discs $\delta_1$.  Similarly, there are discs $\delta_2 \subset U \cap K_{n}$ with boundary on $S$ so that $\sigma(U \cap K_{n};\delta_1)$ contains a product region $P^U_2 = (F_2 \cap U) \times I \subset U \cap \cl(K_{n} - K_{n-1})$ with $(F_2 \cap U) = (F_2 \cap U) \times \{0\}$.  Let $(F'_2 \cap U)$ signify $(F_2 \cap U) \times \{1\}$; it is a subsurface of $S$ except at the remnants of the discs $\delta_2$. The boundaries of the surfaces $F'_1 \cap U$ and $F'_2 \cap U$ are simple closed curves on $S$ which bound discs in $V$.  Let $F'_1$ and $F'_2$ be the surfaces $F'_1 \cap U$ and $F'_2 \cap U$ together with discs in $V$ bound by $\boundary F'_1 \cap U$ and $\boundary F'_2 \cap U$.  Let $P_1$ and $P_2$ be the product regions bounded by $F'_1 \cup F_1$ and $F'_2 \cup F_2$ respectively.  $P^U_1$ and $P^U_2$ are the product regions which are the intersections of $P_1$ with $U$ and $P_2$ with $U$. Let $N' = N \cup P_1$. \newline

Choose a spine for $V$ which intersects each disc of $\delta_1 \cup \delta_2$ exactly once.  We may assume that the spine intersects $P_1$ and $P_2$ in vertical arcs.  Let $\Sigma$ be the intersection of this spine with $N'$.  Corollary \ref{complementary compressionbodies} shows that $U \cap \cl(K_n - K_{n-1})$ and $V \cap \cl(K_n - K_{n-1})$ are compressionbodies.  Since there are not closed components of $\boundary_- \cl(N' - \eta(\boundary N' \cup \Sigma))$, $\cl(N' - \eta(\boundary N' \cup \Sigma))$ is a handlebody and so $\Sigma$ is a reduced spine for $N'$.  (Recall that $\{K_i\}$ is adapted to $S$ and so $U \cap K_{n-1}$ is correctly embedded in $U \cap K_n$.  This is needed to apply Lemma \ref{complementary compressionbodies}.) \newline

We now construct a model splitting of $N'$.  Let $X \cup_W Y$ be any relative Heegaard splitting of $N$ with $Y \cap \fr N = V \cap \fr N$.  Let $\Sigma'$ be a reduced spine for $Y$.  We may assume that $\Sigma' \cap P_2$ consists of vertical arcs. Using the product region $P_1$ we may extend $\Sigma'$ to be a graph in $N'$ whose intersection with $P_1$ consists of vertical arcs.  $\Sigma' \cap \cl(N' - P_1)$ is a reduced spine for $\cl(N' - P_2)$.  \newline

The Reidemeister-Singer theorem and the Marionette Lemma imply that by stabilizing the Heegaard splittings of $N'' = \cl(N' - P_2)$ induced by $\Sigma \cap N''$ and $\Sigma' \cap N''$ they become isotopic.  Perform the necessary stabilizations in such a way that the graphs $\Sigma \cap N''$ and $\Sigma' \cap N''$ still intersect $P_1$ in vertical arcs.  Edge-slides of reduced spines are equivalent to isotopies of the Heegaard surfaces, so there is a sequence of edge-slides which takes (the now stabilized) $\Sigma \cap N''$ to $\Sigma' \cap N''$.  These edge-slides may involve sliding edges of $\Sigma \cap N''$ over other edges or over the surfaces $F'_1 \cup F'_2$.  \newline

These edge-slides define an ambient isotopy of $S \cap N'$ which is fixed off a regular neighborhood of $\cl(N' - P_2)$.  In particular, the isotopy is fixed on $K_{n-2} \cup \cl(M - K_n)$.  After the isotopy, $S \cap K_{n-1}$ is exactly the same as it was before.  Now, however, $S \cap \cl(K_n - K_{n-1})$ is a relative Heegaard surface for $N$ since the model surface was.
\end{proof}

\begin{lemma}\label{perfect adaptations exist}
Suppose that $M = U \cup_S V$ is an end-stabilized absolute Heegaard splitting of $M$.  Then there is an exhausting sequence $\{L_i\}$ which is perfectly adapted to $S$.
\end{lemma}

\begin{proof}
By Section \ref{Balanced Exhausting Sequences} and Corollary \ref{outer collar property 2}, there is an exhausting sequence $\{K_i\}$ which is adapted to $S$, has the outer collar property, and is such that $V \cap \fr K_i$ consists of discs for all $i$.  Recall that, since $S$ is end-stabilized, any time we need to stabilize some $S \cap \cl(K_i - K_j)$ we may do so by a proper ambient isotopy of $S$ in such a way that $K_j$ is fixed throughout the isotopy.  This means the isotopies needed to make each $\cl(K_i - K_j)$ of arbitrarily high genus can be achieved by a single proper ambient isotopy of $S$ in $M$.  \newline

For each $\cl(K_{3i+1} - K_{3i})$, steal handles from infinity and perform the isotopy of $S \cap \cl(K_{3i + 1} - K_{3i})$ needed in order to make $\cl(K_{3i+1} - K_{3i})$ adapted to $S$.  Since each of these isotopies is fixed on $K_{3i - 2}$ their union is a proper ambient isotopy of $S$.  Let $L_i = K_{3i}$ for each $i$.  We claim that $\{L_i\}$ is perfectly adapted to $S$. \newline

It is, of course, adapted to $S$ as each $K_i$ is adapted to $S$ before and after the isotopy.  We need to show that after this isotopy $\cl(K_{3i} - K_{3i -3})$ is adapted to $S$ for $i \geq 2$.  To see this, note that since $V$ intersects each  $\fr K_{3i}$ in discs $V \cap \cl(K_{3i} - K_{3i - 3})$ is a relative compressionbody with preferred surface $S \cap \cl(K_{3i} - K_{3i -3})$ for each $i$. To see that $U \cap \cl(K_{3i} - K_{3i - 3})$ is a relative compressionbody with preferred surface $S \cap \cl(K_{3i} - K_{3i -3})$ note first that $\{L_i\}$ has the outer collar property.  Furthermore, after the isotopy, there are discs $(\delta_1,\boundary \delta_1) \subset (U \cap \cl(K_{3i-2} - K_{3i - 3}),S \cap \cl(K_{3i -2} - K_{3i - 3}))$ which cut off a product region $(U \cap \fr K_{3i - 3}) \times I$ contained in $U \cap (\cl(K_{3i - 2} - K_{3i - 3})) \subset U \cap \cl(K_{3i} - K_{3i - 3})$.  Hence $\{L_i\}$ has both the inner and outer collar properties.  It is easy to see that $\{L_i\}$ is perfectly adapted to $S$ (cf. Section \ref{outer collar prop}). 
\end{proof}

\begin{proof}[Proof of Theorem \ref{Frohman-Meeks 2}]

Suppose, first, that $U_S \cup_S V_S$ and $U_T \cup_T V_T$ are two absolute infinitely stabilized Heegaard splittings of $M$ with the same partition of $\boundary M$.  To show that they are approximately isotopic we will show that given any compact set $C$ there are proper ambient isotopies of $S$ and of $T$ so that after the isotopies, $S$ and $T$ coincide on $C$.  By Section \ref{Balanced Exhausting Sequences} and Corollary \ref{outer collar property 2}, there are exhausting sequences $\{K_i\}$ and $\{L_i\}$ adapted to $S$ and $T$ respectively which have the outer collar property and are such that $V_S \cap \fr K_i$ and $V_T \cap \fr L_i$ consist of discs.  Take subsequences so that $C \subset K_1 \subset L_1 \subset K_2 \subset L_2$.  By Lemma \ref{edge slides again} we may steal handles from infinity for both $S$ and $T$ and then perform further proper ambient isotopies so that $K_1, K_2$ and $\cl(K_2 - K_1)$ are adapted to $S$ and $L_1, L_2,$ and $\cl(L_2 - L_1)$ are adapted to $T$.  By Lemma \ref{mutually adapted} we may steal more handles from infinity and perform more ambient isotopies of $S$ and $T$ so that there is a submanifold $J_1$ containing $K_1$ which is adapted to both $S$ and $T$.  By stealing more handles from infinity, we may stabilize $S \cap J_1$ and $T \cap J_1$ enough times so that they are ambient (in $J_1$) isotopic (Reidemeister-Singer theorem and the Marionette Lemma).  Isotope $S$ and $T$ so that they coincide on $J_1$.  They then also coincide on $C$ and so they are approximately isotopic. \newline

Now suppose that $S$ and $T$ are end-stabilized.  By Lemma \ref{perfect adaptations exist} there are exhausting sequences $\{K_i\}$ and $\{L_i\}$ perfectly adapted to $S$ and $T$ respectively.  By Corollary \ref{mutually adapted exh seq}, we may perform proper ambient isotopies of $S$ and $T$ so that there is an exhausting sequence $\{J_i\}$ perfectly adapted to both of them and is such that, for each $i$, $V_S \cap \fr J_i = V_T \cap \fr J_i$ and the intersections consist of discs.  By Corollary \ref{key cor}, $S$ and $T$ are properly ambient isotopic.
\end{proof}

  \end{document}